\documentclass{amsa}

\newcommand{\editor}{Editors}
\newcommand{\rcptDate}{September 24, 2018}

\usepackage{amsmath}
\usepackage{amssymb}
\usepackage{amsthm}
\usepackage{graphicx}

\newtheorem{thm}{Theorem}[section]

\newtheorem{rem}[thm]{Remark}

\newtheorem{definition}[thm]{Definition}

\renewcommand{\qed}{\ifmmode$\Box$\else{\unskip\nobreak\hfil
	\penalty50\hskip1em\null\nobreak\hfil$\Box$
	\parfillskip=0pt\finalhyphendemerits=0\endgraf}\fi}

\newcommand{\bfu}{{\boldsymbol u}}
\newcommand{\bfv}{{\boldsymbol v}}
\newcommand{\bfw}{{\boldsymbol w}}
\newcommand{\bfF}{{\boldsymbol F}}

\newcommand{\bfz}{{\boldsymbol z}}
\newcommand{\bff}{{\boldsymbol f}}
\newcommand{\bfg}{{\boldsymbol g}}
\newcommand{\bfG}{{\boldsymbol G}}
\newcommand{\bfW}{{\boldsymbol W}}
\newcommand{\bfS}{{\boldsymbol S}}
\newcommand{\bfH}{{\boldsymbol H}}
\newcommand{\bfL}{{\boldsymbol L}}
\newcommand{\bfV}{{\boldsymbol V}}

\begin{document}
\label{page:t}
\vspace*{0.5cm}

\title{VARIATIONAL AND QUASI-VARIATIONAL 
    INEQUALITIES OF NAVIER-STOKES TYPE 
    WITH VELOCITY CONSTRAINTS} \vspace{-0.5cm}

\author{M. Gokieli}
\affiliation{ICM, University of Warsaw, Warsaw, Poland}
\email{maria.gokieli@gmail.com}\vspace{-0.5cm}

\sauthor{N. Kenmochi}
\saffiliation{Faculty of Education, Chiba University, Chiba, Japan}
\semail{nobuyuki.kenmochi@gmail.com}\vspace{-0.6cm}

\tauthor{M. Niezg\'odka}
\taffiliation{CNT Center, Cardinal Stefan Wyszynski University, Warsaw, Poland}
\temail{marekn1506@gmail.com}

\footcomment{
AMS 2010 Subject Classification 34G25, 35G45, 35K51, 35K57, 35K59.\\
}

\maketitle
\noindent
{\bf Abstract.} In this paper we deal with parabolic variational inequalities 
of Navier-Stokes type with time-dependent constraints on velocity fields, 
including gradient constraint case. One of the objectives of this paper is
to propose a weak variational formulation for variational inequalities of 
Navier-Stokes type and to solve them by applying the compactness theorem,
 which was recently developed by the authors (cf. [22]).

Another objective is to approach to a class of quasi-variational inequalities
associated with Stefan/Navier-Stokes problems in which we are taking into 
account the freezing effect of materials in fluids. As is easily understood,
the phase change from liquid into solid gives a great influence to the 
velocity field in the fluid. For instance, in the mussy region, the 
velocity of the fluid
is constrained by some obstacle caused by moving solid. We shall challenge to
the mathematical modeling of Stefan/Navier-Stokes problem as a 
quasi-variational inequality and solve it as an application of parabolic
variational inequalities of Navier-Stokes type.

\vfill

\newpage

\noindent
{\large \bf 1. Introduction}\vspace{0.5cm}

In this paper we study parabolic variational inequalities of
Navier-Stokes type with velocity constraint of the form
 $$ \begin{array}{l}
\displaystyle{
{\bfv}(t) \in K(t),~~0<t<T,~~{\bfv}(0)={\bfv}_0;}\\[0.5cm] 
\displaystyle{\int_Q \{{\bfv}_t\cdot({\bfv}-\boldsymbol {\mathcal \xi}) 
  +\nu \nabla {\bfv}\cdot\nabla({\bfv}-\boldsymbol {\mathcal \xi})+({\bfv}\cdot \nabla){\bfv}
     \cdot({\bfv}-\boldsymbol {\mathcal \xi})\}dxdt }\\[0.5cm]
  \displaystyle{~~~~~~~~~\le \int_Q{\bfg}\cdot({\bfv}-
  \boldsymbol {\mathcal \xi})dxdt, 
~~\forall \boldsymbol {\mathcal \xi}~{\rm with~}\boldsymbol {\mathcal \xi}(t) 
    \in K(t),~~0<t<T,}
\end{array}
\eqno{(1.1)} $$
where $\Omega$ is bounded domain in ${\bf R}^3$ and $Q:=\Omega\times (0,T),~
0<T<\infty$, and $K(t)$ is a prescribed constraint set in the 
3-dimensional solenoidal function space ${\bfW}^{1,2}_{0,\sigma}(\Omega)$;
a positive constant $\nu$, an initial datum ${\bfv}_0$ and a source term
${\bfg}$ are given.
 
We mainly consider the following two cases as $K(t)$:
 $$ K(t)=K^1(\gamma;t):=\{{\bfz}\in{\bfW}^{1,2}_{0,\sigma}(\Omega)~|~|{\bfz}| 
    \le \gamma(\cdot,t)
   ~{\rm a.e.~on~}\Omega\},~0\le t\le T,\eqno{(1.2)}$$
and
 $$ K(t)=K^2(\gamma;t):=\{{\bfz}\in{\bfW}^{1,2}_{0,\sigma}(\Omega)~|~
    |\nabla {\bfz}|
   \le \gamma(\cdot,t)~{\rm a.e.~on~}
    \Omega\},~0\le t \le T,\eqno{(1.3)}$$
where $\gamma$ is a nonnegative continuous function on $\overline Q$ and 
permitted to take $\infty$ somewhere in $\overline Q$; note that the continuity
of $\gamma$ should be understood in the extended sense. Since it is difficult 
to expect the differentiability of the solution ${\bfv}$ in time, we shall
discuss the problem (1.1) in a framework of weak variational inequalities.

There are so many nonlinear dynamical systems in our real world whose mechanism
are still not clear from the theoretical point of view. For instance, we meet
very interesting phase transition phenomena in fluids, such as biofilm growth
and melting ice in the sea or lake. 
Our experiments suggest that in general, biofilm growth is improved 
in fluids, but its mechanism makes complex by 
taking account of freezing or melting effect in fluids. It is a very 
important task for
us to provide some realistic mathematical models to such a complex system. To 
this end the authors started in [20, 21, 22] the development of 
mathematical tools describing various aspects of fluid dynamics; one of them
is the establishment of the theory on variational inequalities of 
Navier-Stokes type with velocity constraints (1.2) and (1.3).

In each of (1.2) and (1.3), there are two cases of $\gamma$, 
the non-degenerate case and degenerate case:
\begin{description}
 \item {(Non-degenerate case)} $~~c_*\le \gamma(x,t)\le \infty$ on 
   $\overline Q$ for a positive constant $c_*$.
 \item {(Degenerate case)} $~~0\le \gamma(x,t)\le \infty$ on 
   $\overline Q$ and $\gamma$ vanishes somewhere in $Q$.
 \end{description}
These cases are separately treated, because in the degenerate case we need
an extended use of Helmholtz decomposition in the solenoidal function spaces 
for the 
construction of a weak solution (see section 4), but in the non-degenerate case
we do not need it and the treatment is much easier. So far as the case
of gradient constraint (1.3) is concerned,
the variational inequality of Navier-Stokes type (1.1) is a new problem and
there has not been any existence result on it in our knowledge.
We shall prove it in the non-degenerate case (see section 3).

In the degenerate case the existence proof of a weak solution of the problem
(1.1) with $K(t)=K^1(\gamma;t)$ was given in the author's paper [20] in which
the crucial step is how to show the strong convergence of regular approximate
solutions in $L^2(Q)^3$. This is quite important to handle well the convergence
of nonlinear terms in approximate variational inequalities for which we needed
the Helmholtz decomposition of solenoidal functions. But,
after the publication of the author's original paper, a gap was found
by the authors in the usage of the 
Helmholtz decomposition. 
We shall make the correction for the gap in section 4 of
this paper under a slight additional assumption on the obstacle function 
$\gamma$ that
  $$ \gamma~{\rm is ~Lipschitz~ continuous~in~a~neighborhood~of~}
       \{(x,t) \in Q~|~\gamma(x,t)=0\}.\eqno{(1.4)}$$
The result in [20] was used in the paper [21] on the biofilm growth
problem, so such an additional assumption as (1.4) should be required in [21],
too.

The second aim of this paper is to give some applications of the results 
obtained in this paper to the Stefan/Navier-Stokes problem, which is
a coupled system of the enthalpy formulation of Stefan problem
with convection,  
 $$ 
  \begin{array}{l}
 \displaystyle{w_t -\Delta \beta(w) +{\bfv}\cdot \nabla w = h(x,t),
   ~~(x,t) \in Q,}\\[0.3cm]
  \displaystyle{w(\cdot,0)=w_0 ~{\rm in~}\Omega,~~\frac {\partial \beta(w)}
   {\partial n}+n_0\beta(w) =0~{\rm on~}
     \Sigma:=\partial \Omega \times (0,T),}
    \end{array} \eqno{(1.5)}$$
and the variational inequality of Navier-Stokes type,
 $$ \begin{array}{l}
\displaystyle{
{\bfv}(t) \in K^i(\gamma(w^{\varepsilon_0});t),~~0<t<T,~~{\bfv}(0)={\bfv}_0;}
\\[0.3cm] 
\displaystyle{\int_Q \{{\bfv}_t\cdot ({\bfv}-\boldsymbol {\mathcal \xi}) 
   +\nu \nabla {\bfv}
   \cdot \nabla({\bfv}-\boldsymbol {\mathcal \xi})+({\bfv}\cdot \nabla){\bfv}
     \cdot({\bfv}-\boldsymbol {\mathcal \xi})\}dxdt}\\[0.5cm]
 \displaystyle{~~~~~~~\le \int_Q{\bfg}\cdot({\bfv}
      -\boldsymbol {\mathcal \xi})dxdt, 
       ~~ \forall \boldsymbol {\mathcal \xi}~{\rm with~}
       \boldsymbol {\mathcal \xi}(t) 
    \in K^i(\gamma(w^{\varepsilon_0});t),~~0<t<T,}
\end{array}
\eqno{(1.6)} $$
where $i=1,~2$, $w_0$ is an initial datum, $h$ is a given source and
$\beta(\cdot)$ is a Lipschitz continuous and increasing function on ${\bf R}$
such that 
$$ \beta(r)=0,~\forall r\in [0,1],~~\beta'(r)>0,~\forall r >1~ {\rm or~} r<0,~~
             \lim_{|r| \to \infty} \frac {\beta(r)}r >0,$$
and $\gamma(\cdot)$ is continuous function from ${\bf R}$ into $[0,\infty]$ 
such that 
  $$ \begin{array}{l}
 \gamma(r)=0,~\forall  r \le 0,~~\gamma~{\rm is~strictly~increasing~on~}
    [0,1),\\[0.3cm]
  \gamma(r) \uparrow \infty~{\rm as~}r\uparrow 1,~~
    \gamma(r)=\infty,~\forall r\geq 1,\\[0.3cm]
  \gamma(r)~{\rm is~Lipschitz~continuous~near~}r=0.
     \end{array} $$ 
We denote by $w^{\varepsilon_0}(x,t)$ the spatial average of $w(x,t)$, namely
  $$ w^{\varepsilon_0}(x,t)=[\rho_{\varepsilon_0}*w(\cdot,t)](x):=
   \int_\Omega \rho_{\varepsilon_0}(x-y)w(y,t)dy,~~
     \forall x \in \Omega,$$
where $\rho_{\varepsilon_0}(\cdot)$ is the usual mollifier on ${\bf R}^3$ with
support in $|x| \le \varepsilon_0$. Throughout this paper, the parameter 
$\varepsilon_0>0$ is fixed, although it is close to $0$, and we do not 
consider the limit $\varepsilon_0 \downarrow 0$. In section 5 we shall prove
the existence of a weak solution $\{w, {\bfv}\}$ to problem (1.5)-(1.6).

Especially, the gradient obstacle problem is a new challenge to the 
Navier-Stokes variational inequalities. From various different motivations
it has been studied by many
researchers so far (cf. [4, 5, 6, 7, 16, 18, 31, 32]), but most of cases were 
treated in the non-degenerate case. It would be expected to generalize to
the degenerate case from some serious physical/mechanical motivations in order
to make more realistic modelings in nonlinear phenomena in fluids 
(cf. [1, 14, 29, 30]).
\vspace{0.5cm}

\noindent
({\bf Notation})

Let 
$\Omega$ be a bounded domain in ${\bf R}^3$ with smooth boundary
$\Gamma:=\partial \Omega$, $Q:=\Omega \times (0,T)$, $0<T<\infty$ and
$\Sigma:=\Gamma \times (0,T)$, and denote by $|\cdot|_X$ the norm in various
function spaces $X$ built on $\Omega$. We consider the usual solenoidal 
function spaces:
\begin{eqnarray*}
 && {\boldsymbol {\mathcal D}}_\sigma(\Omega):=
      \{{\boldsymbol v}=(v^{(1)},v^{(2)},v^{(3)})
       \in {\cal D}(\Omega)^3~|~
  {\rm div}\hspace{0.05cm}{\boldsymbol v}=0~{\rm in~}\Omega\},\\
 &&{\boldsymbol H}_\sigma(\Omega):= {\rm the~ closure~ of~}{\boldsymbol {\mathcal D}}_\sigma(\Omega)~
    {\rm in}~ L^2(\Omega)^3,~{\rm with~norm~}|\cdot|_{0,2},\\
 &&{\boldsymbol V}_\sigma(\Omega):= {\rm the~ closure~ of~}{\boldsymbol {\mathcal D}}_\sigma(\Omega)~
    {\rm in~} H^1_0(\Omega)^3,~{\rm with~norm~} |\cdot|_{1,2},\\
   &&{\boldsymbol W}_\sigma(\Omega):={\rm the~ closure~ of~}
   {\boldsymbol {\mathcal D}}_\sigma(\Omega) ~{\rm in~}W^{1,4}_0(\Omega)^3,\\
    &&~~~~~~~~~~~~{\rm in~ the~ case~of~}(1.2),~
      ~{\rm with~norm~}|\cdot|_{1,4};\\
 &&{\boldsymbol W}_\sigma(\Omega):={\rm the~ closure~ of~}
   {\boldsymbol {\mathcal D}}_\sigma(\Omega)~{\rm in~}W^{2,4}_0(\Omega)^3,\\
 &&~~~~~~~~~~~{\rm in~ the~ case~of~}(1.3);~{\rm with~norm~}|\cdot|_{2,4};
\end{eqnarray*}
in these spaces the norms are given as usual by: for ${\bfv}:=(v^{(1)},v^{(2)},
v^{(3)})$
 $$ |{\boldsymbol v}|_{0,2}:=\left\{ \sum_{k=1}^3 \int_\Omega |v^{(k)}|^2dx 
   \right\}
   ^{\frac 12},~~ |{\boldsymbol v}|_{1,2}:=\left\{ \sum_{k=1}^3 \int_\Omega 
     |\nabla v^{(k)}|^2dx \right\}^{\frac 12}$$
and 
 $$ |{\boldsymbol v}|_{1,4}:=\left\{ \sum_{k=1}^3 \int_\Omega 
  |\nabla v^{(k)}|^4dx  \right\}^{\frac 14},~~
  |{\boldsymbol v}|_{2,4}:=\left\{|{\bfv}|^4_{1,4}+ \sum_{k,j,\ell=1}^3 
   \int_\Omega 
      \left|\frac {\partial^2 v^{(k)}}{\partial x_j\partial x_\ell}\right|^4dx 
   \right\}^{\frac 14}.$$   
For simplicity we denote 
the dual spaces of ${\boldsymbol V}_\sigma(\Omega)$ and 
${\boldsymbol W}_\sigma(\Omega)$ 
by ${\boldsymbol V}^*_\sigma(\Omega)$ and ${\boldsymbol W}^*_\sigma(\Omega)$, 
respectively, which are equipped with their dual norms. 
Also, we denote the inner product in 
${\boldsymbol H}_\sigma(\Omega)$ by $(\cdot,\cdot)_\sigma$ and the duality 
between 
${\boldsymbol V}^*_\sigma(\Omega)$ and ${\boldsymbol V}_\sigma(\Omega)$ by 
$\langle \cdot,\cdot \rangle_\sigma$, namely for 
${\boldsymbol v}_i=(v^{(1)}_i,v^{(2)}_i,v^{(3)}_i),~i=1,2,$ 
$$ ({\boldsymbol v}_1,{\boldsymbol v}_2)_\sigma 
   :=\sum_{k=1}^3\int_\Omega  v_1^{(k)}v_2^{(k)}dx, \qquad
\langle {\bfF}{\boldsymbol v}_1,{\boldsymbol v}_2\rangle_\sigma
    =\sum_{k=1}^3\int_\Omega \nabla v_1^{(k)}\cdot \nabla v_2^{(k)}dx, 
$$
where ${\bfF}$ denotes the duality mapping from ${\bfV}_\sigma(\Omega)$ onto 
${\bfV}^*(\Omega)$.
Then, by identifying the dual of ${\boldsymbol H}_\sigma(\Omega)$ with itself, 
we have:
 $$ {\boldsymbol V}_\sigma(\Omega)
 \subset
 {\boldsymbol H}_\sigma(\Omega)
 \subset {\boldsymbol V}^*_\sigma(\Omega), 
 \quad
 {\boldsymbol W}_\sigma(\Omega)\subset
  C(\overline{\Omega})^3~(C^1(\overline{\Omega})^3~{\rm in~the~case}~(1.3));$$
and all these embeddings are compact.

By the way we introduce the usual simplified notation in the theory 
of Navier-Stokes equations:
  $$ \int_\Omega ({\bfv}\cdot \nabla){\bfw}\cdot {\bfz} dx :=
      \sum_{k,j=1}^3 \int_\Omega v^{(k)}\frac{\partial w^{(j)}}{\partial x_k}
       z^{(j)} dx $$
for ${\bfv}:=(v^{(1)},v^{(2)},v^{(3)}),~{\bfw}:=(w^{(1)},w^{(2)},w^{(3)})$ and
${\bfz}:=(z^{(1)},z^{(2)},z^{(3)})$.

\vspace{0.5cm}

Since $\Omega$ is fixed throughout this paper, the spaces 
${\boldsymbol {\mathcal D}}_\sigma(\Omega),~{\boldsymbol V}_\sigma(\Omega),~
{\boldsymbol H}_\sigma(\Omega)$ and ${\boldsymbol W}_\sigma(\Omega)$ are simply
denoted by ${\boldsymbol {\mathcal D}}_\sigma,~{\boldsymbol V}_\sigma,~
{\boldsymbol H}_\sigma$ and ${\boldsymbol W}_\sigma$, respectively.
When these types of spaces are built on other open set $\Omega'$ in 
${\bf R}^3$,
 we write them as ${\boldsymbol {\mathcal D}}_\sigma(\Omega'),
~{\boldsymbol V}_\sigma(\Omega'),~
{\boldsymbol H}_\sigma(\Omega')$ and ${\boldsymbol W}_\sigma(\Omega')$.
The notation $(\cdot,\cdot)_\sigma$ 
is commonly used for the inner product in ${\bfH}_\sigma$ or in ${\bfH}_\sigma
(\Omega')$ as well as $\langle \cdot,\cdot \rangle_\sigma$ for the duality 
between ${\bfW}^*_\sigma$ and ${\bfW}_\sigma$ or ${\bfW}^*_\sigma(\Omega')$ and
 ${\bfW}_\sigma(\Omega')$ in case of no confusion.

For the general knowledge about solenoidal function spaces and Navier-Stokes 
equation, we refer to the monographs [19, 34].
\vspace{1cm}

\noindent
{\large\bf 2. Time-derivative under constraint and compactness theorem}
\vspace{0.5cm}

\noindent
{\bf 2.1. A compactness theorem}\vspace{0.3cm} 

In this section we recall some results in [22] with the following setup: 
\begin{description}
\item[(h1)] $H$ is a Hilbert space, and its dual $H^*$ is identified with $H$.
\item[(h2)] $V$ is a reflexive
Banach space which is dense and compactly
embedded in $H$, therefore we have
$ V\subset H\subset V^*$ with compact embeddings.
\item[(h3)] $W$ is another reflexive and separable
 Banach space which
is continuously embedded in $V$ and dense in $H$; 
since $H \subset W^*$, we have
   $$ V\subset H\subset W^*~~{\rm with~dense~and~compact
 ~embeddings}. $$
\item[(h4)] 
$V,~V^*,W$ and $W^*$ are strictly convex.
\item[(h5)] The numbers: $p>1$, 
$p':=\frac p{p-1}$, and 
$T>0$ are fixed.
\end{description} \vspace{0.2cm}
	
We begin with the definition of total variation,
which refers here to the time variable. 
For any function $w:[0,T]\to W^*$,
the \textbf{total variation} of $w$, denoted by 
${\rm Var}_{W^*}(w)$,
is defined by
 $$ {\rm Var}_{W^*}(w):=\sup_{
    \begin{array}{l}
    \eta \in C^1_0(0,T;W),\\
   |\eta|_{L^\infty(0,T;W)}\le 1
     \end{array}} \int_0^T
  \langle w,\eta'\rangle_{W^*,W}dt. $$

We refer to [10; Appendice 2] or [15; Chapter 5] 
for the fundamental properties of total
variation functions. Let us now define the set which will be the point of our 
interest in this section.\vspace{0.5cm}

The next lemma is concerned with the compactness property of functions
having bounded total variation from $[0,T]$ into $W^*$.\vspace{0.5cm}

\noindent
{\bf Lemma 2.1.} {\it Let $M$ be any positive number and set
  $$ {\cal X}(M):=\left\{u \in L^p(0,T; H)~\left|~
     \begin{array}{l}
      \displaystyle{|u|_{L^p(0,T;V)} \le M,}\\[0.2cm]
      \displaystyle{ |u|_{L^\infty(0,T;H)}\le M,}\\[0.2cm]
      \displaystyle{{\rm Var}_{W^*}(u) \le M}
     \end{array} \right.\right\}.
  $$
Then we have:
\begin{description}
\item{(1)} Given any sequence $\{u_n\}$ in ${\cal X}(M)$, there is a
subsequence $\{u_{n_k}\}$ of $\{u_n\}$ and a function $u \in {\cal X}(M)$
such that 
  $$u_{n_k}(t) \to u(t)~{\it weakly~ in~} H, ~
      \forall t \in [0,T] ~({\it as~} k \to \infty).$$
Hence, $u_{n_k}(t) \to u(t)$ in $W^*$
for every $t \in [0,T]$.
\item{(2)} ${\cal X}(M)$ is compact and convex in $L^q(0,T;H)$ 
for every $q \in [1,\infty)$.
\end{description}
}\vspace{0.5cm}

See [22; Lemma 3.3] for the proof of Lemma 2.1.\vspace{0.5cm}

\noindent
{\bf Definition 2.1.}
Given $\kappa>0,~M_0 >0$ and $u_0 \in H$, consider the
set $Z_p(\kappa, M_0, u_0)$ in $L^p(0,T;V)\cap L^\infty(0,T;H)$ given by:
$u \in Z_p(\kappa, M_0, u_0)$ if and only if
$ |u|_{L^p(0,T;V)}\le M_0,~|u|_{L^\infty(0,T;H)}\le M_0$ and there exists 
$f \in L^{p'}(0,T;V^*)$ such that
 $$ \int_0^T \langle f,  u\rangle dt \le M_0,
  ~|f|_{L^1(0,T;W^*)} \le M_0$$ and
    $$  \int_0^T \langle \eta'-f, u-\eta \rangle dt \le 
   \frac 12 |u_0-\eta(0)|^2_H,\hspace{3cm}~~~~~~~~~~~~~~~~~~~~~~ \eqno{(2.1)}$$
   $$ \forall \eta 
   \in L^p(0,T;V)~
   {\rm with~}\eta' \in L^{p'}(0,T;V^*),~
  \eta(t) \in \kappa B_W(0),~\forall t \in[0,T],
 $$
where $B_W(0)$ is the closed unit ball in $W$
with center at the origin and $\langle \cdot, \cdot \rangle=
\langle \cdot, \cdot \rangle_{V^*,V}$. \vspace{0.5cm}

\noindent
{\bf Remark 2.1.} In Definition 2.1
the variational inequality (2.1) 
relates $f$ to the time derivative of $u$,
taking into account 
the convex constraint $\kappa B_W(0)$. This is explored as follows.
We note for now that
if $f=u'$ and $u(0)=u_0$, then (2.1) holds for any test function $\eta$. 
Indeed, for
any $u,~\eta \in L^p(0,T;V)$ with $u',~\eta' 
\in L^{p'}(0,T;V^*)$ with $u(0)=u_0$ 
we have by integration by parts
 $$ \int_0^T\langle \eta' -u', u-\eta \rangle dt
 =\frac 12|u(0)-\eta(0)|^2_H
   -\frac 12|u(T)-\eta(T)|^2_H
   \le \frac 12|u_0-\eta(0)|^2_H.$$
Thus, given $u \in L^p(0,T;V)\cap L^\infty(0,T;H)$ and
$u_0 \in H$, the set of all $f$ satisfying (2.1) includes 
$u'$, provided $u'$ exists in $L^{p'}(0,T;V^*)$ and
$u(0)=u_0$. However, in general, it is an extremely
large set; note that in the definition of
$Z_p(\kappa,M_0, u_0)$, any differentiability of
$u$ in time is not required. \vspace{0.5cm}

\noindent
{\bf Lemma 2.2.} {\it Let $Z_p(\kappa, M_0, u_0)$ be the set given by
Definition 2.1. Then there is a positive constant $C^*$ such that
  $$ {\rm Var}_{W^*}(u) \le 
        C^*,~~\forall u \in Z_p(\kappa,M_0,u_0). \eqno{(2.2)}$$
Moreover, we can take $M_0+\frac 1\kappa M_0+\frac 1{2\kappa}|u_0|_H^2$ as
the constant $C^*$ } \vspace{0.5cm}

The uniform estimate (2.2) of total variation for $Z_p(\kappa,M_0,u_0)$ is
directly obtained from variational inequality (2.1). See [20; Lemma 3.2] 
for the proof of Lemma 2.2. 
Combining both of the above lemmas, we arrive at our compactness theorem.
\vspace{0.5cm} 

\noindent
{\bf Theorem 2.1 (cf. [20: Theorem 3.1])} 
{\it Let $\kappa>0,~M_0 >0$ be any numbers and $u_0$ be any element in
$\overline {K(0)}$. Then the set
$Z_p(\kappa,M_0, u_0)$ is relatively compact in
 $L^p(0,T;H)$. Moreover, the convex closure of $Z_p(\kappa,\kappa,M_0)$,
denoted by $\overline{\rm conv}[Z_p(\kappa, M_0,u_0)]$, in $L^p(0,T;V)$
is bounded in $L^p(0,T;V)$ and compact in $L^p(0,T;H)$.} \vspace{0.5cm}

Here we compare Theorem 2.1 with the Aubin compactness theorem [3] 
(or [27; Chapter 1]), saying that for any number $M_0>0$ the set
  $$ \{u~|~
   |u|_{L^p(0,T;V)} \le M_0,~
   |u'|_{L^q(0,T;W^*)} \le M_0 \}, ~1<p<\infty,~ 1<q<\infty, 
  $$
is relatively compact in $L^p(0,T;H)$. We can say in rough that
our compactness theorem is the one obtained by replacing 
$^^ ^^ |u'|_{L^q(0,T;W^*)}\le M_0"$ by the
total variation estimate $^^ ^^ {\rm Var}_{W^*}(u) \le M_0"$ in the Aubin
compactness theorem.\vspace{0.5cm}

\noindent
{\bf Remark 2.2.}   A compactness theorem
of the Aubin type was extended to various directions, for instance
[13] and [23], and further to a quite general setup [33].
\vspace{0.5cm}

\noindent
{\bf 2.2. Time-derivative under convex constraints} \vspace{0.3cm}

We assume here \textbf{(h1)}, \textbf{(h2)}, \textbf{(h4)} 
except for $W$ and \textbf{(h5)}:  
we will not be using the space $W$ here. Again, 
for the sake of simplicity of notation, we
write $\langle\cdot,\cdot\rangle$ for
$\langle\cdot,\cdot\rangle_{V^*,V}$.

As $V^*$ is strictly convex, the duality mapping $F$ from $V$ into $V^*$, 
associated with the gauge function 
$r \to  |r|^{p-1}$, namely $F: V \to V^*$ is the subdifferential of 
$u\to \frac 1p |u|_V^p$, is singlevalued and
demicontinuous from $V$ into $V^*$. \vspace{0.5cm}

\noindent
{\bf Definition 2.2.}
Let $\{K(t)\}_{t \in [0,T]}$ be a family of non-empty, closed and
convex sets in $V$ such that there are functions 
$\alpha \in W^{1,2}(0,T)$ and
$\beta \in W^{1,1}(0,T)$ satisfying the following property:
for any $s, t \in [0,T]$ and any
$z \in K(s)$ there is $\tilde z \in K(t)$ such that
$$
|\tilde z -z|_H \le |\alpha(t)-\alpha(s)|
 (1+|z|^{\frac p2}_V),~~
 |\tilde z|^p_V - |z|^p_V \le |\beta(t)-\beta(s)|
 (1+|z|^p_V).
 \eqno{(2.3)} $$
We denote by $\Phi(\alpha,\beta)$ the set of all such families $\{K(t)\}$, and 
put
 $$\Phi_S:=\bigcup_{\alpha \in W^{1,2}(0,T),~
 \beta\in W^{1,1}(0,T)} \Phi(\alpha,\beta),
  $$
which is called the {\bf strong class of 
time-dependent convex sets}. \vspace{0.5cm}

Given $\{K(t)\} \in \Phi_S$, we consider the following time-dependent 
convex function on $H$:
$$
  \varphi_K^t(z):= \left \{
        \begin{array}{l}
 \displaystyle {\frac 1p |z|^p_V+I_{K(t)}(z), ~~~{\rm if~}z \in K(t),}\\[0.3cm]
 \displaystyle {\infty,~~~~~~~~~~~~~~~~~~~~~{\rm otherwise},}
        \end{array} \right. 
$$
where $I_{K(t)}(\cdot)$ is the indicator function of
$K(t)$ on $H$. For each $t \in [0,T]$, 
$\varphi^t_K(\cdot)$
is proper, l.s.c. and strictly convex on $H$ and on $V$. 
By the general theory on nonlinear
evolution equations generated by time-dependent 
subdifferentials, see [24], 
condition (2.3) is a sufficient condition in order that
for any $f \in L^2(0,T;H)$ and $u_0 \in \overline{K(0)}$ (the closure of 
$K(0)$ in $H$), the Cauchy problem
 $$ u'(t)+\partial \varphi^t_K(u(t))
 \ni f(t),~~u(0)=u_0,~{\rm in~}H,$$
admits a unique solution $u$ such that
$u \in C([0,T];H)\cap L^p(0,T;V)$ with $u(0)=u_0$, 
$t^{\frac 12}u'\in L^2(0,T;H)$ and 
$t \to t \varphi^t_K(u(t))$ is bounded on 
$(0,T]$ and absolutely continuous on any
compact interval in $(0,T]$.

Next, taking constraints of obstacle type into account, we introduce a weak 
class of time-dependent convex sets. To this end, we first 
recall a  notion of convergence of time--dependent convex sets introduced in
[17].
This convergence is defined by means of admissible geometrical perturbations:
we choose them to be homothetic and parallel transformations. 
To set this up, we define a perturbation operator 
${\cal F}_\varepsilon(t):V\to V$, which can be a sum of 
an expansion / contraction with an $\varepsilon$--dependent modulus 
(close to $1$) and an $(\varepsilon,t)$--dependent (small) 
parallel transformation.
Roughly, the sets will be considered to be close one to the other if 
after  this kind of perturbation, at any time moment, they 'fit'
one to the other.
Note that we do not include rotations in our perturbation operator 
in order to avoid complexity which would be irrelevant 
from the point of view of applications. 
However, this can be done; see [26]. \vspace{0.5cm}

\noindent
{\bf Definition 2.3.} Let $c_0$
be a fixed constant and $\sigma_0$ be a fixed 
function in $C([0,T];V)$ with $\sigma'_0 \in L^{p'}(0,T;V^*)$. Associated
with these $c_0$ and $\sigma_0$,  for any
small positive number $\varepsilon$,
the mapping ${\cal F}_\varepsilon: [0,T]\times V \to V$ is defined by
$$
{\cal F}_\varepsilon(t)z 
=(1+ \varepsilon c_0)z +\varepsilon \sigma_0(t),
  ~\forall t \in [0,T],~\forall z \in V.
 \eqno{(2.4)} $$
Let $\{K(t)\} _{t \in [0,T]}$ be a family of non-empty, closed and
convex sets in $V$ 
and $\{K_n(t)\} _{t \in [0,T]}$ a sequence of such families. 
We say that \textbf{$\{K_n(t)\}$
converges to $\{K(t)\}$} as $n \to\infty$,
which is denoted by 
 $$ K_n(t)\Longrightarrow K(t)~~\text{ on~}[0,T]~(\text{ as~}n \to \infty),
 $$ 
if for any 
$\varepsilon \in (0,\varepsilon_1]~(0<\varepsilon_1 <1) $ 
there is a positive integer 
$N_\varepsilon$ satisfying
 $$ {\cal F}_\varepsilon(t) (K_n(t)) \subset K(t),~~
{\cal F}_\varepsilon(t) (K(t)) \subset K_n(t),~~ 
      \forall t \in [0,T],~~\forall n\geq 
   N_\varepsilon.$$

Note that this notion of convergence depends on 
the choice of the perturbation operator ${\cal F}_\varepsilon(t)$, 
which depends itself on 
the constant $c_0$ and on the function $\sigma_0$. 
The operator's form defines the perturbations that we allow, and that we can 
further restrict by choosing concrete 
$c_0$ and $\sigma_0$. 
As we are going to see in the examples, 
it is often enough to take them as  equal to $0$ or $\pm 1$.
We are now ready to define the weak class of constraints, 
which is the closure of the strong 
class with respect to this convergence. \vspace{0.5cm}

\noindent
{\bf Definition 2.4.}
$\{K(t)\}\in \Phi_W$, the {\bf weak class of time-dependent
convex sets}, if and only if the following two conditions are satisfied:
\begin{description}
\item {(a)} $K(t)$ is a closed and convex set in $V$ for all 
$t \in[0,T]$,
\item {(b)} there exists a
sequence $\{\{K_n(t)\}\}_{n \in {\bf N}} \subset \Phi_S $ such that 
$ K_n(t)\Longrightarrow K(t)$ on $[0,T]$, as $n \to \infty$, according to 
Definition 2.3.
\end{description}\vspace{0.5cm}

We give typical examples of $\{K(t)\}$ in the weak class $\Phi_W$.
\vspace{0.5cm}

\noindent{\bf Example 2.1.}
 Let $\Omega$ be a bounded smooth domain in 
${\bf R}^3$ and $Q:=\Omega \times (0,T)$. Let
 $$H:={\bfH}_\sigma(\Omega),~~~V:={\bfV}_\sigma(\Omega),~ 
 $$ 
Moreover, let 
$\gamma:=\gamma(x,t) \in C(\overline Q)$ such that 
   $$\gamma \geq c_* ~~{\rm on~} \overline Q $$
for a positive constant $c_*$ with $0<c_*<1$, and choose a sequence 
$\{\gamma_n\}$ in 
$C^2(\overline Q)$ such 
that $\gamma_n\geq c_*$ and $\gamma_n \to \gamma$ in  $C(\overline Q)$. 
Now, constraint sets $K(t)$ and $K_n(t)$ are defined by
 $$ K^1(t):=\{{\bfz} \in {\bfV}_\sigma~|~|{\bfz}(x)| \le \gamma(x,t)~{\rm for~a.e.~}x \in \Omega\},~
   \forall t \in [0,T],$$
and 
$$ K^1_n(t):=\{{\bfz} \in {\bfV}_\sigma~|~|{\bfz}(x)| \le \gamma_n(x,t)~{\rm for~a.e.~}x \in 
\Omega\},~
   \forall t \in [0,T].$$
Given $\varepsilon >0$, take a positive integer
$N_\varepsilon$ so that
    $$ |\gamma_n-\gamma |\le \varepsilon c_* 
     ~~{\rm on~}\overline Q,~\forall n \geq 
     N_\varepsilon. \eqno{(2.5)}$$

In this case, with the choice of $c_0=-\frac 1{c_*}$ and $\sigma\equiv 0$, 
consider the mapping ${\cal F}_\varepsilon(t)$ of the form
${\cal F}_\varepsilon(t){\bfz} = (1-\frac \varepsilon{c_*}){\bfz}$, which 
maps $V$ into itself for all small $\varepsilon>0$. 
Then we have:

(i) We show that $\{K_n(t)\} \in \Phi_S$. Fix $n$ and note that $\gamma_n\in
C^2(\overline Q)$. Therefore it is possible to take a partition $0=t_0
<t_1<t_2<\cdots <t_k(n)=T$ of $[0,T]$ so that
  $$ |\gamma_n(s)-\gamma_n(t)|_{C(\overline\Omega)} \le L(\gamma_n)|s-t|
     < c_*,~
     \forall s,~t \in [t_{k-1},t_k],~k=1,2,\cdots, k(n), \eqno{(2.6)}$$
where $L(\gamma_n)$ is the Lipschitz constant of $\gamma_n$.
Now, given ${\bfz}\in K^1_n(s)$ and $s,~t \in [t_{k-1},t_k]$, 
the function
$\tilde {\bfz}=(1-\frac {|\gamma_n(s)-\gamma_n(t)|_{C(\overline \Omega)}}
{c_*}){\bfz}$ satisfies by (2.6)
  \begin{eqnarray*}
  |\tilde{\bfz}(x)| 
    &=& \left(1-\frac {|\gamma_n(s)-\gamma_n(t)|
       _{C(\overline \Omega)}}
   {c_*}\right)|{\bfz}(x)|\\
   &\le& \left(1-\frac {|\gamma_n(s)-\gamma_n(t)|_{C(\overline \Omega)}}
    {c_*} \right) \gamma_n(x,s)\\
   &\le& \gamma_n(x,s)-|\gamma_n(s)-\gamma_n(t)|_{C(\overline\Omega)}\\
    &\le& \gamma_n(x,s)-|\gamma_n(x,s)-\gamma_n(x,t)|\le \gamma_n(x,t)
  \end{eqnarray*}
Thus $\tilde{\bfz} \in K^1_n(t)$ and
 $$ |\tilde{\bfz}-{\bfz}|_{0,2} = \frac{
      |\gamma_n(s)-\gamma_n(t)|_{C(\overline\Omega)}}{c_*}
   |{\bfz}|_{0,2} \le \frac{L(\gamma_n)}{c_*}|s-t|
     |{\bfz}|_{0,2}.   $$ 

Generally, for any $s,~t \in [0,T]$ with $s <t$ and 
${\bfz}\in K^1_n(s)$, by 
repeating the above argument we can find $\tilde {\bfz} \in K^1_n(t)$
such that 
  $$ |{\bfz}-\tilde{\bfz}|_{0,2} \le L_n|s-t|
          (1+|{\bfz}|_{0,2}),$$
for a positive constant $L_n$ depending only on $n$.
Also, we have
$|\nabla \tilde{\bfz}|^2_{0,2} \le 
         |\nabla{\bfz}|^2_{0,2}$, too. 
Therefore $\{K^1_n(t)\}$ belongs to $\Phi(\alpha,\beta)$ with 
$p=2,~\alpha(t)=L_nt$ and $\beta(t)=0$. Thus $\{K^1_n(t)\} \in \Phi_S$. 
\vspace{0.5cm}

(ii) $\{K^1(t)\} \in \Phi_W$. 
In fact, for any ${\bfz}_n \in K^1_n(t)$, we have by (2.5)
 \begin{eqnarray*}
   |{\cal F}_\varepsilon {\bfz}_n| &=&\left(1-\frac \varepsilon{c_*}\right)
  |{\bfz}_n|\\
   &\le& \left(1-\frac \varepsilon{c_*}\right)\gamma_n(\cdot,t) 
   \le \left(1-\frac \varepsilon{c_*}\right)(\gamma(\cdot,t)
    + \varepsilon c_*)\\
   &\le& \gamma(\cdot,t)+\varepsilon(c_*-1-\varepsilon) \le \gamma(\cdot,t)
   ~{\rm a.e.~on~} \Omega,
   \end{eqnarray*}
which implies ${\cal F}_\varepsilon (K^1_n(t))\subset K^1(t)$. 
Similarly, 
${\cal F}_\varepsilon (K^1(t)) \subset K^1_n(t)$. Hence 
$K^1_n(t) \Longrightarrow K^1(t)$ on $[0,T]$, and thus
$\{K^1(t)\} \in \Phi_W$. \vspace{0.5cm}

\noindent
{\bf Example 2.2.} Under the same situation as in Example 2.1, let us consider
the gradient constraint case:
 $$ K^2(t):=\{{\bfz}\in {\bfV}_\sigma~|~|\nabla {\bfz}| \le
     \gamma(\cdot,t)~{\rm a.e.~on~}\Omega\},~~0\le t\le T.$$
and
$$ K^2_n(t):=\{{\bfz}\in {\bfV}_\sigma~|~|\nabla {\bfz}| \le
     \gamma_n(\cdot,t)~{\rm a.e.~on~}\Omega\},~~0\le t\le T.$$
Then we see just in the same way as Example 2.1 that $\{K^2_n(t)\} 
\in \Phi_S$ and $K^2_n(t) \Longrightarrow 
K^2(t)$ on $[0,T]$ (as $n \to \infty$). Hence 
$\{K^2(t)\} \in \Phi_W$. 
\vspace{0.5cm}

Next, we introduce the time-derivative under constraint $\{K(t)\} \in \Phi_W$.
Put
  $$ {\mathcal K}:=\{v \in L^p(0,T;V)~|~ v(t) \in K(t)
  ~{\rm for~a.e.~}t \in [0,T]\} $$
and
 $$ {\mathcal K}_0:=\{ \eta
 \in{\mathcal K} ~|~ \eta'\in 
  L^{p'}(0,T;V^*)\}.$$

\vspace{0.5cm}

\noindent
{\bf Definition 2.5.} Let $\{K(t)\} \in \Phi_W$
and $u_0 \in \overline {K(0)}$.
Then we define an operator $L_{u_0}$ whose graph is given as follows:
$ f \in L_{u_0}u$ if and only if
 $$u \in {\cal K},~f \in L^{p'}(0,T; V^*),~~
 \int_0^T \langle \eta'-f, u-\eta\rangle dt 
  \le \frac 12 |u_0-\eta(0)|^2_H, 
  ~~\forall \eta \in 
     {\mathcal K}_0. $$
\vspace{0.2cm}

In the next theorems we mention some of important properties of 
$L_{u_0}$.\vspace{0.5cm}

\noindent
{\bf Theorem 2.2.} {\it Let $\{K(t)\} \in \Phi_W$ and $u_0 \in 
\overline{K(0)}$. Then, $L_{u_0}$ is maximal monotone from 
$D(L_{u_0})\subset L^p(0,T;V)$
 into 
$L^{p'}(0,T;V^*)$, and the domain $D(L_{u_0})$
is included in the set $\{u\in C([0,T];H)\cap
{\cal K} ~|~u(0)=u_0\}$.}
\vspace{0.5cm}

The characterization and fundamental properties of the mapping $L_{u_0}$
are given in the following theorem.\vspace{0.5cm}

\noindent
{\bf Theorem 2.3.} {\it Let 
$\{K(t)\} \in \Phi_W$. Then we have: 
\begin{description}
\item{(1)} Let $u_0 \in \overline{K(0)}$.
Then $f \in L_{u_0}u$ if and only if there are
$\{\{K_n(t)\}\} \subset \Phi_S$, $\{u_n\} \subset L^p(0,T;V)$ with 
$u_n \in {\cal K}_n:=\{v \in L^p(0,T;V)~|~v(t) \in K_n(t)~for~a.e~t
\in [0,T]\}$ and $u'_n \in 
L^{p'}(0,T;V^*)$, $\{f_n\} \subset 
L^{p'}(0,T;V^*)$ such that 
      $$ K_n(t) \Longrightarrow K(t)~{\rm on~}
      [0,T],$$
 $$ u_n \to u~{\rm in~}C([0,T];H)~and~weakly~in~
     L^p(0,T;V), $$
 $$ f_n\to f~weakly~in~L^{p'}(0,T;V^*), $$
  $$  \int_0^T \langle u'_n-f_n,u_n -v\rangle dt \le 0, ~\forall v \in
      {\mathcal K}_n,~\forall n,$$
  $$\limsup_{n\to \infty}\int_{t_1}^{t_2} \langle f_n,u_n\rangle dt
 \le \int_{t_1}^{t_2} \langle f,u \rangle dt, ~\forall t_1,~t_2~with
~0\le t_1\le t_2 \le T.
  $$
\item{(2)} 
Let $u_0 \in \overline{K(0)}$
and $f \in L_{u_0}u$. Then, for any $t_1,~t_2 \in 
[0,T]$ with $t_1 \le t_2$,
$$ \int_{t_1}^{t_2}\langle \eta'-f,u-\eta\rangle dt
  + \frac 12|u(t_2)-\eta(t_2)|^2_H \le 
 \frac 12|u(t_1)-\eta(t_2)|^2_H,~~\forall
 \eta \in {\cal K}_0.$$
\item{(3)} Let $u_{i0} \in \overline {K(0)},$ and 
$f_i \in L_{u_{i0}}u_i$ for $i=1,2$. Then,
for any $t_1,~t_2 \in [0,T]$ with $t_1\le t_2$, 
 $$ \frac 12 |u_1(t_2)-u_2(t_2)|^2_H \le 
 \frac 12 |u_1(t_1)-u_2(t_2)|^2_H +
  \int_{t_1}^{t_2} \langle f_1-f_2, u_1-u_2 \rangle dt.
  \eqno{(2.7)}$$
\end{description}
}\vspace{0.5cm}
We refer to [22; Theorems 5.1 and 5.2] for the precise proof of the above 
theorems. \vspace{0.5cm}

\noindent
{\bf Remark 2.3.} If $0 \in K(t)$ for all $t \in [0,T]$, then 
the relation $f \in L_{u_0}u$ gives
 $$ \frac 12|u(s)|^2_H \le \frac 12|u_0|^2_H +\int_0^s \langle f, u\rangle 
    dt,~~\forall s \in [0,T]; $$
in fact, noting $0 \in L_00$ and using (2.7), we get the above
inequality.\vspace{0.5cm}

Once the maximal monotonicity of $L_{u_0}$ is proved, it is quite useful for
the weak solvability of parabolic variational inequalities with time-dependent
constraint $K(t)$. In fact, for any coercive maximal monotone 
or pseudomonotne operator $A: D(A)=
L^p(0,T;V)  \to L^{p'}(0,T;V^*)$ and $f \in L^{p'}(0,T;V^*)$ we see by the
general theory on nonlinear operators of monotone type (cf. [25; section 5])
the functional inclusion
  $$ L_{u_0}u+ Au \ni f~~~{\rm in ~}L^{p'}(0,T; V^*) \eqno{(2.8)}$$
admits a solution $u$.
\vspace{1cm}

\noindent
{\large\bf 3. Variational inequalities of Navier-Stokes 
type}\vspace{0.5cm}

\noindent
{\bf 3.1. Weak formulation of variational inequalities of Navier-Stokes type}
\vspace{0.3cm}

We are given a nonnegative function $\gamma=\gamma(x,t)$ on $\overline Q$ 
as an obstacle function such that
$ 0\le \gamma(x,t) \le \infty$ for all $(x,t) \in \overline Q$. For simplicity
we use the following notation: for any constant $c \in [0,\infty]$,
  $$\overline Q(\gamma = c):=\{(x,t)\in \overline Q~|~\gamma(x,t)= c\},$$
  $$ \overline Q(\gamma \geq c):=\{(x,t)\in \overline Q~|~\gamma(x,t)\geq c\},~
     \overline Q(\gamma \le c):=\{(x,t)\in \overline Q~|~\gamma(x,t)\le c\},$$
and similarly $\overline Q(\gamma > c)$ and $\overline Q(\gamma < c)$ are 
defined. Besides, for the set $\hat Q:=\Omega\times [0,T]$, $\hat Q(\gamma=c),~
\hat Q(\gamma\geq c),~{\rm etc.}$ are similarly defined, too. \vspace{0.5cm}

We assume that $\gamma$ is 
continuous from $\overline Q$ into $[0,\infty]$, namely,
 $$ \left\{
    \begin{array}{l}
 {\rm the~set~}\overline Q(\gamma=\infty)
  ~{\rm is~closed~in~} \overline Q,\\[0.18cm] 
     \forall \kappa \in (0,\infty),~\gamma~{\rm is~continuous~on~} 
      \overline Q(\gamma \le \kappa),
      \\[0.18cm]
     \forall M \in (0,\infty), {\rm there~is~an~open~set~}U_M~(\subset
      {\bf R}^4)
     ~{\rm containing~}
      \overline Q(\gamma =\infty) \\
     {\rm ~~~such ~that~} \gamma \geq M~{\rm on~}U_M\cap \overline Q.
    \end{array} 
    \right. \eqno{(3.1)} $$
It is easily seen that (3.1) is equivalent to the  continuity on $\overline Q$ 
in the usual sense, of the function
  $$\alpha(x,t):=\left \{
        \begin{array}{l}
       \displaystyle{  \frac {\gamma(x,t)}{1+\gamma(x,t)},
     ~~~{\rm if~}0\le \gamma(x,t) < \infty,}\\[0.5cm]
         1.~~~~~~~~~~~~~~~{\rm if~}\gamma(x,t)=\infty,
        \end{array}
     \right. $$

We are now ready to formulate our problem with
the constraint sets $K(t)$ given by
  $$K^1(\gamma;t):=\{{\bfz}\in {\bfV}_\sigma~|~|{\bfz}|\le 
   \gamma(\cdot,t)~{\rm a.e.~on~}\Omega\},~ ~0\le t\le T, $$
or
  $$K^2(\gamma;t):=\{{\bfz}\in {\bfV}_\sigma~|~|\nabla{\bfz}|
    \le \gamma(\cdot,t)~{\rm a.e.~on~}\Omega\},~~0\le t\le T, $$
and the classes of test functions
 $$\boldsymbol{\cal K}^i(\gamma):=\{\boldsymbol {\mathcal \xi}\in
   L^2(0,T;{\bfV}_\sigma)
 ~|~\boldsymbol {\mathcal \xi}(t) \in K^i(\gamma;t),~{\rm a.e.~}t \in (0,T)\},~
 i=1,2, $$
 $$ \boldsymbol{\cal K}^1_0(\gamma):=\left\{{\boldsymbol {\mathcal \xi}} \in 
C^1([0,T];{\boldsymbol W}_\sigma)~\left|~
      \begin{array}{l}
    {\boldsymbol {\mathcal \xi}}\in K^1(\gamma;t),~\forall t\in [0,T],\\[0.2cm]
   {\rm supp}(|{\boldsymbol {\mathcal \xi}}|) \subset \hat Q(\gamma>0),
      \end{array} \right.\right\} $$
and
$$ \boldsymbol{\cal K}^2_0(\gamma):=\left\{{\boldsymbol {\mathcal \xi}} \in 
C^1([0,T];{\boldsymbol W}_\sigma)~\left|~
      \begin{array}{l}
      {\boldsymbol {\mathcal \xi}}\in K^2(\gamma;t),~\forall t\in [0,T],
     \\[0.2cm]
      {\rm supp}(|\nabla{\boldsymbol {\mathcal \xi}}|) \subset \hat Q(\gamma>0)
      \end{array} \right. \right\}
       $$
and ${\rm supp}(|\boldsymbol {\mathcal \xi}|)$ and ${\rm supp}(|\nabla
   \boldsymbol {\mathcal \xi}|)$ denote the 
supports of $|\boldsymbol {\mathcal \xi}|$ and $|\nabla 
\boldsymbol {\mathcal \xi}|$, respectively. 
 \vspace{0.5cm}

\noindent
{\bf Definition 3.1.} 
For given data
 $$  \nu >0~(\text{constant}),
     ~ {\boldsymbol g}\in L^2(0,T; {\boldsymbol H}_\sigma),~
   {\boldsymbol v}_0 \in {\boldsymbol H}_\sigma,
 $$
our problem, referred to $NS^i(\gamma;{\boldsymbol g},{\boldsymbol v_0})$, 
is to find a function
${\boldsymbol v}:=(v^{(1)},v^{(2)},v^{(3)})$ from $[0,T]$ into 
${\boldsymbol H}_\sigma$ satisfying the following (i) and (ii).
\begin{description}
\item{(i)} ${\boldsymbol v}(0)
={\boldsymbol v}_0$ in~${\boldsymbol H}_\sigma$, and
$t\mapsto ({\boldsymbol v}(t), \boldsymbol {\mathcal \xi}(t))_\sigma$ is of  
bounded variation on
$[0,T]$ for all ${\boldsymbol {\mathcal \xi}} \in \boldsymbol{\cal K}^i_0
(\gamma)$.
 \item{(ii)} $ \sup_{t \in [0,T]}|{\bfv}(t)|_{0,2}< \infty$, 
   ${\bfv}\in \boldsymbol{\cal K}^i(\gamma)$
and 
   $$\int_0^t({\boldsymbol {\mathcal \xi}}'(\tau), {\boldsymbol v}(\tau)
  -{\boldsymbol {\mathcal \xi}}(\tau))_\sigma d\tau
    +\nu \int_0^t \langle {\bfF}{\boldsymbol v}(\tau), {\boldsymbol v}(\tau)
  - {\boldsymbol {\mathcal \xi}}(\tau)\rangle
      _\sigma d\tau ~~~~~~~~~~~~~~~~~~~~~~~~~~$$
  $$~~+ \int_0^t \int_\Omega ({\boldsymbol v}(x,\tau)\cdot\nabla){\bfv}
(x,\tau)\cdot
      ({\boldsymbol v}(x,\tau)-{\boldsymbol {\mathcal \xi}}(x,\tau))dx d\tau 
     +\frac 12 |{\boldsymbol v}(t)-{\boldsymbol {\mathcal \xi}}(t)|^2_{0,2}  
   \eqno{(3.2)}$$
  $$ ~~~~~~~~~~~~~~~\le \int_0^t ({\boldsymbol g}(\tau),{\boldsymbol v}(\tau)
    -{\boldsymbol {\mathcal \xi}}
(\tau))_\sigma d\tau
     + \frac 12 |{\boldsymbol v}_0-{\boldsymbol {\mathcal \xi}}(0)|^2_{0,2},
~~\forall t \in [0,T],
    ~\forall {\boldsymbol {\mathcal \xi}} \in \boldsymbol{\cal K}^i_0(\gamma).~~~~~~~~~~ $$
\end{description}
Such a function ${\bfv}$ is called a weak solution of 
$NS^i(\gamma;{\bfg},{\bfv}_0),~i=1,2$.
 \vspace{0.5cm}

\noindent
{\bf Remark 3.1.} In the non-degenerate case of $\gamma$, namely 
$ \gamma \geq c_* ~(>0)$ on $\overline Q$, $\boldsymbol{\cal K}^i_0(\gamma)$ 
is simply described as
 $$\boldsymbol{\mathcal K}^i_0(\gamma):=\{{\boldsymbol v} \in 
C^1([0,T];{\boldsymbol W}_\sigma(\Omega))~|~
      {\boldsymbol v}\in K^i(\gamma;t),~\forall t\in [0,T]\}.\eqno{(3.3)}$$
Hence, (i) of Definition 3.1 implies that $t\mapsto ({\boldsymbol v}(t), 
\boldsymbol {\mathcal \xi}(t))_\sigma$ is of  
bounded variation on
$[0,T]$ for all ${\boldsymbol {\mathcal \xi}} \in C^1([0,T];{\bfW}_\sigma)$.
Similary, in the degenerate case of $\gamma$, it implies in the case $i=1$ 
that 
$t\mapsto ({\boldsymbol v}(t), \boldsymbol {\mathcal \xi}(t))_\sigma$ is of  
bounded variation on
$[0,T]$ for all ${\boldsymbol {\mathcal \xi}} \in C^1([0,T];{\bfW}_\sigma)$
with ${\rm supp}(|\boldsymbol {\mathcal \xi}|) 
\subset \hat Q(\gamma>0)$.\vspace{0.5cm}

We note in Definition 3.1 that $\bfv$ is defined for \emph{every} $t\in[0,T]$
 according to the given $\bfv_0$,
even if we do not require it to be continuous in time:
our definition permits jumps in time, including 
the initial time $t=0$.
What we will prove, is that $\bfv$ is a limit of continuous
approximate solutions. \vspace{0.5cm}

The first main result of this paper is stated as follows.\vspace{0.5cm}

\noindent
{\bf Theorem 3.1 (Non-degenerate case).} {\it In addition to (3.1), assume 
that 
 $$\gamma \geq c_*~~{\rm on~}\overline Q ~~{\rm for~a~positive~constant~}c_*.
    \eqno{(3.4)}$$
Let $i=1$ or $2$. Let ${\bfg}\in L^2(0,T;{\bfH}_\sigma)$ and
${\bfv}_0 \in K^i(\gamma; 0)$ and 
  $$ \left \{ 
     \begin{array}{l}
     |{\bfv}_0| \in L^\infty(\Omega) ~~~{\it in~the~case~of~}i=1,
       \\[0.2cm]
      |\nabla {\bfv}_0| \in L^\infty(\Omega)~~~{\it in~the~case~of~}i=2.
     \end{array} \right. \eqno{(3.5)}
  $$
Then there exists at least one weak solution of 
$NS^i(\gamma;{\bfg},{\bfv}_0)$; namely
\begin{description}
\item{(i)} ${\bfv}(0)
={\bfv}_0$ in~${\bfH}_\sigma$, and
$t\mapsto ({\bfv}(t), \boldsymbol {\mathcal \xi}(t))_\sigma$ is of  
bounded variation on
$[0,T]$ for all $\boldsymbol {\mathcal \xi} \in C^1([0,T];{\bfW}_\sigma)$,
\item{(ii)} $\sup_{t \in [0,T]}|{\bfv}(t)|_{0,2}< \infty$, 
${\bfv}\in \boldsymbol{\cal K}^i(\gamma)$
and (3.2) holds.
\end{description}}
\vspace{0.5cm}

As to the weak solvability of $NS^1(\gamma;{\bfg},{\bfv}_0)$ in the degenerate
case we have:
\vspace{0.5cm}

\noindent
{\bf Theorem 3.2 (Degenerate case).} {\it Assume (3.1) holds, ${\bfg}\in
L^2(0,T;{\bfH}_\sigma)$ and ${\bfv}_0 \in K^1(\gamma;0) \cap 
L^\infty(\Omega)^3$ with 
 $${\rm supp}(|{\bfv}_0|) \subset \{x \in \Omega~|~\gamma(x,0)>0\}. 
    $$
Moreover assume that for each $t \in [0,T]$ and all small $\kappa>0$,
$\gamma(x,t)$ is Lipschitz continuous on
$\overline \Omega(\gamma(\cdot,t)\le \kappa):=\{x \in \overline \Omega~|~\gamma(x,t) \le \kappa\}$, namely
   $$|\gamma(x,t)-\gamma(x',t)|\le L_\gamma(t,\kappa)|x-x'|,
    ~~\forall x,~\forall x'
     \in \overline \Omega(\gamma(\cdot,t)\le \kappa), \eqno{(3.6)}$$ 
where $L_\gamma(t,\kappa)$ is a positive constant depending on $t, \kappa$.
Then $NS^1(\gamma;{\bfg},{\bfv}_0)$ has at least one weak solution ${\bfv}$
in the sense of Definition 3.1.}
\vspace{0.5cm}

In the degenerate case of $\gamma$, 
we need the Helmholtz decomposition of solenoidal functions (cf. [19]) in the 
construction of weak solution of $NS^1(\gamma;{\bfg},{\bfv}_0)$; regarding
problem $NS^2(\gamma; {\bfg},{\bfv}_0)$ the degenerate case
is still open question except some special cases of $\gamma$. \vspace{0.5cm}

\noindent
{\bf Remark 3.2.} The degenerate case of $\gamma$, namely 
$\gamma=0$ somewhere in $Q$, problem 
$NS^1(\gamma;{\bfg},{\bfu}_0)$ was earlier discussed without 
condition (3.6) in the statement of the main result 
([20; Theorem 1.1]). However, after the publication of this result, 
unfortunately
a gap was found in the proof by the authors. In this paper, we shall make the
precise correction in section 4 of this paper.
\vspace{0.5cm}

\noindent
{\bf 3.2. Approximation of $NS^i(\gamma; {\bfg},{\bfv}_0),~i=1,2$ } 
\vspace{0.3cm}

Our main theorems will be proved in two steps of 
\begin{itemize} 
  \item Approximation of $NS^i(\gamma;{\bfg},{\bfv}_0)$,
  \item Convergence of approximate solutions.
\end{itemize}
We begin with the approximation of $\gamma$ given by
  $$ \gamma_{\delta,N}(r):= (\gamma(r) \lor \delta) \land N,
  ~r \in {\bf R},   ~0<\delta < 1, ~N>0,$$
where $a \lor b=\max\{a,b\},~a\land b=\min\{a,b\}$ for any real numbers 
$a,~b$.
Clearly, $\gamma_{\delta,N}$ is everywhere bounded, strictly positive and 
continuous on $\overline Q$, and $\gamma_{\delta,N} \to \gamma$ as 
$\delta \downarrow 0$ and $N\uparrow \infty$ in the sense that
  $$ \left \{ 
     \begin{array}{l}
  \forall M>0,~\exists \delta_M>0,~\exists N_M,
     ~\exists M'_M >0~{\rm such~ that~}\\[0.1cm]
    \hspace{2cm} \overline Q(\gamma >M) \subset \overline 
      Q(\gamma_{\delta,N} >M'_M),~
           \forall \delta \in (0,\delta_M),~\forall N >N_M,\\[0.3cm]  
    \forall k >0, ~\gamma_{\delta,N} \to \gamma
      ~{\rm uniformly~on~}\overline Q(\gamma \le k)~{\rm as~}\delta \to 0
       ~{\rm and}~ N\to \infty.
     \end{array} \right. \eqno{(3.7)}
  $$

The approximate problem $NS^i(\gamma_{\delta,N};{\bfg},{\bfv}_0)$ to 
$NS^i(\gamma; {\bfg}, {\bfv}_0)$ is formulated as follows: Find 
a function ${\bfv}$ such that
  $$ {\bfv} \in C([0,T];{\bfH}_\sigma)\cap \boldsymbol{\mathcal K}^i
(\gamma_{\delta,N}),
   ~~ {\bfv}(0)={\bfv}_0,
  $$
and
 $$ \int_0^t \langle {\boldsymbol {\mathcal \xi}}', 
    {\bfv}-\boldsymbol {\mathcal \xi}\rangle_\sigma d\tau 
    +\nu \int_0^t\langle {\bfF}{\bfv},{\bfv}-\boldsymbol {\mathcal \xi} 
   \rangle_\sigma d\tau ~~~~~~~~~~~~~~~~~~~~~~~~$$ 
  $$ +\int_0^t\int_\Omega ({\bfv}\cdot \nabla){\bfv}\cdot
    ({\bfv}-\boldsymbol {\mathcal \xi})dx d\tau
   +\frac 12|{\bfv}(t)-\boldsymbol {\mathcal \xi}(t)|^2_{0,2} \eqno{(3.8)}$$
 $$ \le \int_0^t \langle {\bfg},{\bfv}-\boldsymbol {\mathcal \xi}\rangle_\sigma d\tau  +\frac 12|{\bfv}_0-\boldsymbol {\mathcal \xi}(0)|^2_{0,2},
   ~~\forall \boldsymbol {\mathcal \xi}\in \boldsymbol{\mathcal K}^i_0
   (\gamma_{\delta,N}),
    \forall t \in [0,T],$$
where
 $$ \boldsymbol{\mathcal K}^i(\gamma_{\delta,N}):=\{\boldsymbol {\mathcal \xi}
     \in L^2(0,T;{\bfV}_\sigma)~|~
     \boldsymbol {\mathcal \xi}(t) \in K^i(\gamma_{\delta,N};t), 
     ~{\rm a.e.~}t \in [0,T]\},$$
$$ \boldsymbol{\mathcal K}^i_0(\gamma_{\delta,N}):=\{\boldsymbol {\mathcal \xi}
    \in C^1([0,T];{\bfW}_\sigma)~|~ \boldsymbol {\mathcal \xi}(t) 
    \in K^i(\gamma_{\delta,N};t),~\forall t \in [0,T]\}$$
with
  $$K^1(\gamma_{\delta,N};t):=\{{\bfz}\in {\bfV}_\sigma~|~|{\bfz}|
 \le \gamma_{\delta,N}(\cdot,t)~{\rm a.e.~on~}\Omega\},~\forall t \in [0,T],$$
  $$K^2(\gamma_{\delta,N};t):=\{{\bfz}\in {\bfV}_\sigma~|~|\nabla{\bfz}|
 \le \gamma_{\delta,N}(\cdot,t)~{\rm a.e.~on~}\Omega\},~\forall t \in [0,T].$$
 
As was seen in Examples 2.1 and 2.2
these classes $\{K^i(\gamma_{\delta,N};t)\}$, $i=1,2$, belong to the weak class $\Phi_W$, with $p=2$, of time-dependent convex sets in ${\bfV}_\sigma$ and 
the time derivative ${\bfL}^i_{{\bfv}_0}(\delta,N)$ associated with 
constraints $\{K^i(\gamma_{\delta,N};t)\}$ is well defined by Theorem 2.2.
By definition, ${\bff} \in {\bfL}^i_{{\bfv}_0}(\delta,N){\bfv}$
if and only if ${\bff} \in L^2(0,T;{\bfV}^*_\sigma),~{\bfv} \in 
\boldsymbol{\mathcal K}^i(\gamma_{\delta,N})$ and 
  $$ \int_0^T \langle \boldsymbol {\mathcal \xi}' -{\bff}, {\bfv}
    -\boldsymbol {\mathcal \xi} \rangle_\sigma dt \le 
    \frac 12|{\bfv}_0-\boldsymbol {\mathcal \xi}(0)|^2_{0,2},~~
  \forall \boldsymbol {\mathcal \xi} \in 
    \boldsymbol{\mathcal K}^i_0(\gamma_{\delta,N}).
   \eqno{(3.9)}$$

Since $\gamma_{\delta,N} \le N$ on ${\bf R}$, there is a positive 
constant $C_N$, depending only on $N$, such that
  $$  K^i(\gamma_{\delta,N};t) \subset K^*:=\{{\bfz}\in {\bfV}_\sigma~|~
        |{\bfz}| \le C_N~{\rm a.e.~on~}\Omega\},~\forall t \in [0,T],~
     i=1,~2.$$
With this set $K^*$,
we introduce a mapping ${\bfG}(\cdot, \cdot): 
 K^*\times {\bfV}_\sigma \to {\bfV}^*_\sigma$ defined by
 $$ \langle {\bfG}({\bfw},{\bfv}), {\bfz}\rangle_\sigma 
     := \int_\Omega ({\bfw}\cdot \nabla){\bfv}\cdot{\bfz} dx,~
  \forall {\bfw} \in K^*,~\forall {\bfv},~{\bfz}\in {\bfV}_\sigma.$$
For any fixed ${\bfw} \in K^*$, the mapping 
${\bfz} \mapsto {\bfG}({\bfw},{\bfz})$ is bounded, linear and monotone 
from ${\bfV}_\sigma$ into ${\bfV}^*_\sigma$, because
  $$ \langle {\bfG}({\bfw},{\bfz}),{\bfz}\rangle_\sigma =0,~
      ~\forall {\bfz} \in {\bfV}_\sigma, \eqno{(3.10)}$$
by the divergencefreeness of ${\bfw}$.\vspace{0.5cm}

\noindent
{\bf Proposition 3.1.} {\it Let ${\bfg} \in L^2(0,T;{\bfH}_\sigma)$, 
${\bfv}_0 \in K^i(\gamma,0)$ 
and assume
(3.5) holds for the initial datum ${\bfv}_0$. Then, for any small positive 
$\delta>0$ and large $N>0$ there is one and only one solution ${\bfv}$ of 
$NS^i(\gamma_{\delta,N};{\bfg},{\bfv}_0)$ and it is given by the solution of
 $$ {\bfg} \in {\bfL}^i_{{\bfv}_0}(\delta,N){\bfv}+\nu {\bfF}{\bfv}
   +{\bfG}({\bfv},{\bfv})~~{\it in~}L^2(0,T;{\bfV}^*_\sigma).$$
} 

Prior to the proof of Proposition 3.1 we prepare two lemmas.\vspace{0.5cm}

Since
${\bfW}_\sigma$ is compactly embedded in $C(\overline \Omega)^3$ for 
$i=1$ and in $C^1(\overline \Omega)^3$ for $i=2$,
there is a 
positive constant $\kappa_{\delta,N}$, depending on $\delta$ and $N$, 
such that
 $$  \kappa_{\delta,N} B_{{\bfW}_\sigma}(0) \subset K^i(\gamma_{\delta,N};t), 
  ~\forall t \in [0,T],     \eqno{(3.11)}$$
where $B_{{\bfW}_\sigma}(0)$ is the closed unit ball around the origin 
in ${\bfW}_\sigma$. \vspace{0.5cm}

The first lemma follows immediately from the definition of ${\bfG}:
K^*\times {\bfV}_\sigma \to {\bfV}_\sigma^*.$ \vspace{0.5cm}

\noindent
{\bf Lemma 3.1.} {\it We have that 
  $$|\langle {\bfG}({\bfw},{\bfv}),{\bfz} \rangle_\sigma|
     \le  C_N |{\bfv}|_{1,2}|{\bfz}|_{0,2},~~
   \forall {\bfw}\in K^*,~\forall {\bfv},~{\bfz} \in {\bfV}_\sigma. 
     $$ }

By condition (3.5), ${\bfv}_0 \in K^i(\gamma_{\delta,N};0)$ for all $\delta>0$
and large $N$. 
Let us consider a functional inclusion of the form:
 $$ {\bfL}^i_{{\bfv}_0}(\delta,N){\bfv}+\nu{\bfF}{\bfv}+{\bfG}({\bfw},{\bfv})
  \ni {\bfg}   ~~{\rm in ~}L^2(0,T;{\bfV}_\sigma). \eqno{(3.12)}$$
We observe easily that the mapping ${\bfv}\mapsto \nu {\bfF}{\bfv}
+{\bfG}({\bfw},{\bfv})$ is 
maximal monotone, strictly monotone and coercive from $L^2(0,T;{\bfV}_\sigma)$ 
into $L^2(0,T;{\bfV}^*_\sigma)$. Hence, from the general theory on monotone 
operators in Banach spaces (cf. [25]) it follows that the range of the sum
${\bfL}^i_{{\bfv}_0}(\delta,N)+\nu{\bfF}+{\bfG}({\bfw},\cdot)$ is the 
whole of 
$L^2(0,T;{\bfV}^*_\sigma)$, namely there is one and only one 
${\bfv} \in L^2(0,T;{\bfV}_\sigma)$ which satisfies (3.12); we denote by
${\bfS}^i{\bfw}$ the solution ${\bfv}$ and obtain from the energy estimate in
Remark 2.3 that
 $$ \frac 12|{\bfv}(t)|^2_{0,2}+\nu \int_0^t 
   |{\bfv}(\tau)|^2_{{\bfV}_\sigma}
      d\tau \le \frac 12 |{\bfv}_0|^2_{0,2}+\int_0^t ({\bfg}(\tau),
    {\bfv}(\tau))_\sigma dt,
     ~~\forall t \in [0,T].$$
Hence 
 $$ \sup_{t \in [0,T]}|{\bfv}(t)|^2_{0,2} 
   + \nu|{\bfv}|^2_{L^2(0,T;{\bfV}_\sigma)}
   \le |{\bfv}_0|^2_{0,2}+\frac {C_0^2}\nu 
          |{\bfg}|^2_{L^2(0,T;{\bfH}_\sigma)}=:
  M_0({\bfu}_0,{\bfg}) \eqno{(3.13)} $$
for a positive constant $C_0$ satisfying $|{\bfz}|_{0,2}\le C_0|{\bfz}|_{1,2}$
for all ${\bfz} \in {\bfV}_\sigma$. 

On account of this result,
with
  $$\boldsymbol{\mathcal K}^*:=\left\{{\bfw} \in L^2(0,T;{\bfV}_\sigma)~\left|~
  \sup_{t \in [0,T]}
   |{\bfw}(t)|^2_{0,2}+
  \nu|{\bfw}|^2_{L^2(0,T;{\bfV}_\sigma)}\le M_0({\bfu}_0,{\bfg}) \right.
   \right\},$$
we can define an operator 
$\boldsymbol{\mathcal S}^i: \boldsymbol{\mathcal K}^i(\gamma_{\delta,N})\cap 
\boldsymbol{\mathcal K}^*
 \to \boldsymbol{\mathcal K}^i(\gamma_{\delta,N})\cap \boldsymbol{\mathcal K}^*$ by putting 
   $${\mathcal S}^i{\bfw}={\bfv}.$$  

\noindent
{\bf Lemma 3.2.} {\it For $i=1,~2$, the operator
${\mathcal S}^i$ is compact in 
$\boldsymbol{\mathcal K}^i(\gamma_{\delta,N})\cap \boldsymbol{\mathcal K}^*$ 
with respect to the 
topology of $L^2(0,T; {\bfH}_\sigma)$. }\vspace{0.3cm}

\noindent
{\bf Proof.} 
Let ${\bfw}_n$ be any sequence in $\boldsymbol{\mathcal K}^i(\gamma_{\delta,N})
\cap \boldsymbol{\mathcal K}^*$
such that ${\bfw}_n \to {\bfw}$ weakly in $L^2(0,T;{\bfH}_\sigma)$ (as
$n\to \infty$). Clearly, ${\bfw} \in \boldsymbol{\mathcal K}^i
(\gamma_{\delta,N}) \cap
\boldsymbol{\mathcal K}^*$, since 
it is closed and convex in $L^2(0,T;{\bfH}_\sigma)$. 
Also, putting ${\bfv}_n:={\mathcal S}^i{\bfw}_n$, we observe from (3.13) that
$ {\bff}_n:={\bfg}-\nu {\bfF}{\bfv}_n-{\bfG}^i({\bfw}_n,{\bfv}_n) \in 
{\bfL}^i_{{\bfv}_0}(\delta,N){\bfv}_n$ and $\{{\bff}_n\}$ is
bounded in $L^2(0,T;{\bfV}^*_\sigma)$. Noting (3.11) and applying Theorem 2.1 
for the triplet
     $${\bfV}_\sigma \subset {\bfH}_\sigma \subset
          {\bfV}^*_\sigma, $$
we see that $\{{\bfv}_n\}$ is relatively compact in $L^2(0,T;{\bfH}_\sigma)$
and hence in $L^2(Q)$. 

Now, we extract a subsequence $\{{\bfv}_{n_k}\}$ so that
${\bfv}_{n_k} \to {\bfv}$ in $L^2(0,T;{\bfH}_\sigma)$ as $k \to \infty$ for
some ${\bfv}$; note that ${\bfv} \in \boldsymbol{\mathcal K}^i
(\gamma_{\delta,N})\cap 
\boldsymbol{\mathcal K}^*$. 
Also, from (3) of Theorem 2.3 it follows that for all $k,~j$
 $$ \frac 12 |{\bfv}_{n_k}(t)-{\bfv}_{n_j}(t)|^2_{0,2} +\nu \int_0^t
     |{\bfv}_{n_k}-{\bfv}_{n_j}|^2_{1,2}d\tau~~~~~~~~~~~~~~~~~$$ 
  $$~~~~~~~~~~~~~~~~\le
    - \int_0^t \langle {\bfG}({\bfw}_{n_k},{\bfv}_{n_k})
     -{\bfG}({\bfw}_{n_j},{\bfv}_{n_j}),{\bfv}_{n_k}
       -{\bfv}_{n_j}\rangle_\sigma d\tau =:I_{k,j}(t). $$
Here we estimate $I_{k,j}(t)$ by Lemma 3.1 and (3.10) as follows:
  \begin{eqnarray*}
  I_{k,j}(t) 
   &\le& -\int_0^t \langle {\bfG}({\bfw}_{n_k},{\bfv}_{n_k}-{\bfv}_{n_j})
     ,{\bfv}_{n_k}-{\bfv}_{n_j}\rangle_\sigma d\tau \\
   & &~~~~~~~~~~ -\int_0^t \langle {\bfG}({\bfw}_{n_k}-{\bfw}_{n_j},
       {\bfv}_{n_j}),{\bfv}_{n_k}-{\bfv}_{n_j}\rangle_\sigma d\tau \\
   &=&-\int_0^t \langle {\bfG}({\bfw}_{n_k}-{\bfw}_{n_j},{\bfv}_{n_j})
     ,{\bfv}_{n_k}-{\bfv}_{n_j}\rangle_\sigma d\tau \\
   &\le& 2 C_N \int_0^T |{\bfv}_{n_j}|_{1,2}
            |{\bfv}_{n_k}-{\bfv}_{n_j}|_{0,2} d\tau \to 0 ~({\rm as~}k,~j 
          \to \infty). 
   \end{eqnarray*} 
As a consequence, $ {\bfv}_{n_k}\to {\bfv}$ in $C([0,T];
{\bfH}_\sigma)\cap L^2(0,T; {\bfV}_\sigma)$ and ${\bfG}({\bfw}_{n_k},
{\bfv}_{n_k}) \to {\bfG}({\bfw},{\bfv})$ weakly in $L^2(0,T;{\bfV}^*_\sigma)$
as $k\to \infty$, whence
${\bff}_{n_k} \to {\bff}:= {\bfg}-\nu{\bfF}{\bfv}-{\bfG}({\bfw},{\bfv})$
weakly in $L^2(0,T;{\bfV}^*_\sigma)$ as $k \to \infty$. Hence, by the 
demiclosedness of maximal monotone operators, ${\bff} \in 
{\bfL}^i_{{\bfv}_0}(\delta,N){\bfv}$, namely
  $$ {\bfL}^i_{{\bfv}_0}(\delta,N){\bfv} +\nu{\bfF}{\bfv}+{\bfG}({\bfw},
   {\bfv}) \ni  {\bfg},$$
which shows that ${\mathcal S}^i{\bfw}={\bfv}$. By the uniqueness of solution
of this functional inclusion, it is concluded that ${\mathcal S}^i{\bfw}_n
={\bfv}_n\to {\bfv}={\mathcal S}^i{\bfw}$ in $L^2(0,T;{\bfV}_\sigma)$ without 
extracting any subsequence from $\{{\bfw}_n\}$. Thus ${\mathcal S}^i$ is 
compact
in $\boldsymbol{\mathcal K}^i(\gamma_{\delta,N})\cap \boldsymbol{\mathcal K}^*$  with respect to the 
topology of $L^2(0,T;{\bfH}_\sigma)$. \hfill $\diamondsuit$ \vspace{0.5cm}

\noindent
{\bf Proof of Proposition 3.1:} By virtue of Lemma 3.2, the operator
${\mathcal S}^i$ admits at least one fixed point, 
${\bfv}={\mathcal S}^i{\bfv}$ in
$\boldsymbol{\mathcal K}^i(\gamma_{\delta,N})$. Besides, by the 
definition (3.9) of 
${\bfL}^i_{{\bfv}_0}(\delta,N)$,
this fixed point ${\bfv}$ satisfies (3.8),
namely ${\bfv}$ is a solution of 
$NS^i(\gamma_{\delta,N};{\bfg},{\bfv}_0)$. Finally we prove the uniqueness of 
solution.
Let ${\bfv}$ and $\bar {\bfv}$ be two solutions 
$NS^i(\gamma_{\delta,N}; {\bfg},{\bfv}_0)$. Then, by (3) of Theorem 2.3, 
we have
 $$ \frac 12 |{\bfv}(t)-\bar{\bfv}(t)|^2_{0,2} 
      +\nu \int_0^t |{\bfv}-\bar{\bfv}|^2_{1,2} d\tau 
    \le -\int_0^t\langle {\bfG}({\bfv},{\bfv})
     -{\bfG}(\bar{\bfv},\bar{\bfv}), {\bfv}-\bar {\bfv}\rangle_\sigma d\tau 
     =:I(t).
   \eqno{(3.14)}$$
Just as in the proof of Lemma 3.2, the right hand side of (3.14) is dominated
by
 \begin{eqnarray*} 
 I(t) &\le&
 -\int_0^t \langle {\bfG}({\bfv}-\bar {\bfv},{\bfv}),{\bfv}-\bar{\bfv}\rangle
    _\sigma d\tau 
  = \int_0^t \langle {\bfG}({\bfv}-\bar {\bfv},{\bfv}-\bar{\bfv}),{\bfv}
     \rangle_\sigma d\tau \\
  & \le& 2C_N \int_0^T |{\bfv}-\bar{\bfv}|_{1,2}
            |{\bfv}-\bar{\bfv}|_{0,2} d\tau \\
   &\le& 2C_N \left\{\varepsilon \int_0^t 
     |{\bfv}-\bar{\bfv}|^2_{1,2}d \tau
           + c_\varepsilon \int_0^t |{\bfu}-\bar{\bfu}|^2_{0,2} \tau \right\},
 \end{eqnarray*}
where $\varepsilon$ is any small positive number and $c_\varepsilon$ is a 
positive number depending only on $\varepsilon$. Now, taking $\varepsilon>0$
so as to satisfy $2C_N \varepsilon <\frac \nu 2$, we get from (3.14) that
 $$  |{\bfv}(t)-\bar{\bfv}(t)|^2_{0,2} 
      +\nu \int_0^t |{\bfv}-\bar{\bfv}|^2_{1,2} d\tau 
    \le C' \int_0^t |{\bfv}-\bar{\bfv}|^2_{0,2} d\tau,
    ~\forall t \in [0,T].$$
where $C'$ is a positive constant. According to the Gronwall inequality, we
see that ${\bfv}=\bar{\bfv}$ on $[0,T]$. Thus the solution of 
$NS^i(\gamma_{\delta,N};{\bfg}, {\bfu}_0)$ is unique and the proof of 
Proposition 3.1
is now complete. \hfill $\diamondsuit$ 
\vspace{0.5cm}

\noindent
{\bf Remark 3.3.} In the case when $\gamma$ is everywhere bounded and strictly 
positive on $\overline Q$, we have $\gamma=\gamma_{\delta,N}$ for all small 
$\delta>0$ and large $N>0$, so that Proposition 3.1 gives that 
$NS^i(\gamma;{\bfg},{\bfv}_0)$ has
one and only one weak solution ${\bfv}$ in $C([0,T];{\bfH}_\sigma)\cap
\boldsymbol{\cal K}^i(\gamma)$.\vspace{0.5cm}

\noindent
{\bf 3.3. Proof of Theorem 3.1 (Non-degenerate case)}\vspace{0.3cm}

In this subsection, we accomplish the proof of Theorem 3.1 by proving the 
convergence of approximate solutions constructed by Proposition 3.1.

By assumption (3.4), we see that
    $$\gamma_{\delta,N}(x,t) = \gamma(x,t)\land N=:\gamma_N(x,t),~~\forall 
       \delta ~{\rm with~}0<\delta \le c_*,~\forall {\rm large}~N>0. 
     $$

Let ${\bfv}_N$ be the approximate solution of $NS^i(\gamma_N;{\bfg},
{\bfv}_0), ~i=1,~2,$ for large $N>0$. Then, by (3.13) 
there exists sequences 
$\{N_n\}$ tending to $\infty$ (as $n \to \infty $) such that
  $${\bfv}_n:={\bfv}_{N_n} \to {\bfv}~{\rm weakly~in~}L^2(0,T;
    {\bfV}_\sigma)~~{\rm and ~weakly}^*~{\rm in~}L^\infty(0,T;{\bfH}_\sigma)
    \eqno{(3.15)} $$
for a certain function ${\bfv} \in L^2(0,T;{\bfV}_\sigma)\cap
L^\infty(0,T;{\bfH}_\sigma)$, satisfying
  $$ {\bfv} \in \boldsymbol{\mathcal K}^i(\gamma).  \eqno{(3.16)}$$
We note that ${\bfv}_n$ is the solution of
   $$ {\bff}_n :={\bfg}-\nu {\bfF}{\bfv}_n-{\bfG}({\bfv}_n,{\bfv}_n)
      \in {\bfL}^i_{{\bfv}_0}(c_*,N_n){\bfv}_n,$$
which is equivalent to the variational inequality (cf. (3.8)):
  $$ {\bfv}_n \in C([0,T];{\bfH}_\sigma)\cap \boldsymbol{\cal K}^i(\gamma_{N_n}),~
   ~{\bfv}_n(0)={\bfv}_0,$$
 $$ \int_0^t \langle \boldsymbol {\mathcal \xi}'-{\bff}_n, 
    {\bfv}_n-\boldsymbol {\mathcal \xi}\rangle_\sigma d\tau 
   +\frac 12|{\bfv}_n(t)-\boldsymbol {\mathcal \xi}(t)|^2_{0,2} 
    \le \frac 12|{\bfv}_0-\boldsymbol {\mathcal \xi}(0)|^2_{0,2} 
     \eqno{(3.17)} $$
   $$\forall \boldsymbol {\mathcal \xi}\in \boldsymbol{\mathcal K}^i_0
   (\gamma_{N_n}),
    \forall t \in [0,T].$$
Here we note that for any $\boldsymbol {\mathcal \xi} \in 
C^1([0,T];{\bfW}_\sigma) \subset 
C(\overline Q)^3$
  \begin{eqnarray*}
    \int_0^T | \langle {\bfG}({\bfv}_n,{\bfv}_n),\boldsymbol {\mathcal \xi}
     \rangle |dt
     &\le& C\int_0^T|{\bfv}_n|_{0,2}|{\bfv}_n|_{1,2}
           |\boldsymbol {\mathcal \xi}|_{{\bfW}_\sigma} dt\\
     &\le & CM_0({\bfg},{\bfv}_0)|\boldsymbol {\mathcal \xi}|_{L^\infty(0,T;
     {\bfW}_\sigma)}
   \end{eqnarray*}
for some positive constant $C$ independent of $n$. 
This shows that
$\{{\bfG}({\bfv}_n,{\bfv}_n)\}$ is bounded in $L^1(0,T;{\bfW}^*_\sigma)$,
so is $\{{\bff}_n\}$ in $L^1(0,T;{\bfW}^*_\sigma)$. 
Also, by non-degenerate
condition $c_*>0$, we can choose a positive number $\kappa$ so that
    $$ \kappa B_{{\bfW}_\sigma}(0) \subset K^i(\gamma_{N_n};t),~
       \forall t \in [0,T],~~\forall n, ~i=1,2. $$
Hence, for a sufficient large constant $M'_0 >0$ it follows that 
$${\bfv}_n \in Z_2(\kappa,M'_0, {\bfv}_0), \forall n,$$
so that $\{{\rm Var}_{W^*_\sigma}({\bfv}_n)\}$ is bounded by Lemma 2.2.
Moreover, by virtue of Theorem 2.1 and  Lemma 2.1, $\{{\bfv}_n\}$ is 
relatively compact in $L^2(0,T;{\bfH}_\sigma)$ and  
there is a subsequence 
$\{{\bfv}_{n_k}\}$ such that 
${\bfv}_{n_k} \to {\bfv}~{\rm~ in~} L^2(0,T;{\bfH}_\sigma) ~{\rm as~}
     k\to \infty$ and ${\bfv}_{n_k}(t) \to {\bfv}(t)$ weakly in 
${\bfH}_\sigma$ for every $t \in [0,T]$; we may assume that
the limit function ${\bfv}$ is the same one
as in (3.15). For simplicity we write this subsequence by $\{{\bfv}_n\}$,
again; we have together with (3.15) and (3.16)
   $$ \begin{array}{l}
   {\bfv}_n \to {\bfv}~{\rm in~}L^2(0,T;{\bfH}_\sigma)~{\rm and~weakly~in~}
        L^2(0,T;{\bfV}_\sigma),\\[0.3cm]
    {\bfv}_n(t) \to {\bfv}(t) ~{\rm weakly~in~} {\bfH}_\sigma,~~\forall
     t \in [0,T]. 
      \end{array} 
   \eqno{(3.18)}$$
This implies that 
  $$\int_0^t \langle {\bfG}({\bfv}_n,{\bfv}_n), \boldsymbol {\mathcal \xi}
     \rangle_\sigma d\tau
    \to \int_0^t \langle {\bfG}({\bfv},{\bfv}), \boldsymbol {\mathcal \xi}
     \rangle_\sigma d\tau,~ 
  \forall \boldsymbol {\mathcal \xi} \in C^1([0,T];{\bfW}_\sigma),
    ~\forall t\in [0,T].
     \eqno{(3.19)}$$
Taking any function $\boldsymbol {\mathcal \xi} \in
\boldsymbol{\cal K}^i_0(\gamma)$, we see that 
$\boldsymbol {\mathcal \xi}$ is a test function for (3.17), 
since $ \boldsymbol {\mathcal \xi} \in \boldsymbol{\cal K}^i_0
( \gamma_{N_n})$ for all large $n$.
Now, substitute $\boldsymbol {\mathcal \xi}$ in (3.17) and pass to the limit as $n \to \infty$ to
have by (3.18) and (3.19) that
  $$ \begin{array}{l}
  \displaystyle{
 \int_0^t \langle \boldsymbol {\mathcal \xi}', {\bfv}
    -\boldsymbol {\mathcal \xi}\rangle_\sigma d\tau
    +\nu \int_0^t \langle {\bfF}{\bfv}, {\bfv}-\boldsymbol {\mathcal \xi}
    \rangle_\sigma
    +\int_0^t \langle {\bfG}({\bfv},{\bfv}),{\bfv}-\boldsymbol {\mathcal \xi}
      \rangle_\sigma d\tau }\\[0.3cm]
 \displaystyle{~~~~~~~
   +\frac 12|{\bfv}(t)-\boldsymbol {\mathcal \xi}(t)|^2_{0,2} \le 
 \frac 12|{\bfv}_0-\boldsymbol {\mathcal \xi}(0)|^2_{0,2}+\int_0^t \langle 
   {\bfg}, 
      {\bfv}-\boldsymbol {\mathcal \xi}\rangle_\sigma d\tau } 
  \end{array}  $$
and ${\bfv}$ is of bounded variation as a function from $[0,T]$ into 
${\bfW}^*_\sigma$, whence the function $t \mapsto ({\bfv}(t),\boldsymbol 
{\mathcal \xi}(t))_\sigma$ 
is of bounded variation on $[0,T]$ for each $\boldsymbol {\mathcal \xi} \in 
\boldsymbol{\mathcal K}^i_0(\gamma)$ (cf. (3.3)
in Remark 3.1). Clearly ${\bfv}(0)={\bfv}_0$, because 
${\bfv}_n(0)={\bfv}_0 \to {\bfv}(0)$ weakly in ${\bfH}_\sigma$ by (3.18). 
Thus ${\bfv}$
is a weak solution of $NS^i(\gamma;{\bfg},{\bfv}_0)$. \hfill $\diamondsuit$
\vspace{0.5cm}

\noindent
{\bf Remark 3.4.} The obstacle function
$\gamma$ was approximated by $\gamma_{\delta,N}$ and this approximation
satisfies (3.7), which was used in the proof of Theorem 3.1. As is
easily checked, Theorem 3.1 can be proved by means of any approximation 
of $\gamma$, as long as (3.7) is fulfilled, although we need some easy
modifications in the proof. \vspace{0.5cm}

\noindent
{\bf Remark 3.5.} When $\gamma$ is so close to $0$ on some region $Q'$,
for the solution ${\bfv}$ of $NS^i(\gamma;{\bfg},{\bfv}_0)$ we see
that
  $$|{\bfv}|~{\rm is~ close~ to~} 0~{\rm on~}Q'
    ~~{\rm in~the~ case~ of~}i=1, $$ 
or 
  $$|\nabla {\bfv}| ~{\rm is~ close~ to~} 0~{\rm on~}Q' 
       ~~{\rm in~the~ case~ of~}i=2. $$
The former means that the velocity ${\bfv}$ is close to $0$ on such a 
region $Q'$, and the 
latter that ${\bfv}$ is close to a vector field independent of space 
variable $x$, namely it depend almost only on time, but ${\bfv}$ itself 
is not necessarily close to $0$ on $Q'$.
\vspace{1cm}

\noindent
{\large\bf 4.  Degenerate case of $NS^1(\gamma;{\bfg},{\bfv}_0)$}
\vspace{0.5cm}

\noindent
{\bf 4.1. Helmholtz decomposition}\vspace{0.3cm}

In this subsection we suppose in addition to (3.1) that $\gamma$ is 
nonnegative and satisfies (3.6).

Given sequences $\{\delta_n\}$ with $\delta_n \downarrow 0$ and 
$\{N_n\}$ with $N_n \uparrow
\infty$, we put
   $$\gamma_n(x,t)= (\gamma(x,t)\lor \delta_n)\land N_n,~\forall
    x \in \Omega,~\forall t \in [0,T],~\forall n \in {\rm N}. $$

Let $t_0$ be fixed in $[0,T]$ with $\Omega_0:=\{x \in \Omega~|~
\gamma(x,t_0)>0\} \ne \emptyset$ and in this subsection denote $\gamma(x,t_0)$
simply by $\gamma(x)$ for $x \in \Omega$.

For any function ${\bfz}$
in $L^2(\Omega_0)^3$, consider the Helmholtz-decomposition (cf. [19, 34])
 $$ {\bfz}=\tilde{\bfz} +\nabla q,~~\tilde{\bfz} \in {\bfH}_\sigma(\Omega_0),~~
  q \in L^2_{loc}(\Omega_0),~\nabla q \in L^2(\Omega_0)^3;
 $$
we note that the correspondence ${\bfz}\to \tilde {\bfz}$ is the projection of 
${\bfz}$ onto ${\bfH}_\sigma(\Omega_0)$ and $\Omega_0$ is an open set in 
$\Omega$ which is possibly not Lipschitz.
\vspace{0.5cm}

\noindent
{\bf Lemma 4.1.} {\it Let ${\bfz},~ \tilde{\bfz},~q$ be as above, 
$M >0$ be any number and put
    $$ q^M(x) :=(q(x)\land M)\lor (-M),~~x \in \Omega_0.$$
Then $q^M \in H^1(\Omega_0)$ and ${\bfz}^M :=\tilde{\bfz}+\nabla q^M
\to {\bfz}$ in $L^2(\Omega_0)^3$ as $M \uparrow \infty.$}\vspace{0.3cm}

\noindent
{\bf Proof.} Note that $|q^M(x)| \le M$ and $|\nabla q^M(x)|\le
|\nabla q(x)|$ for a.e. $x \in \Omega_0$. Besides, since $q^M \to q$
a.e. on $\Omega_0$ as $M \uparrow \infty$, we have $ |\nabla q^M| \uparrow
|\nabla q|$ a.e. on $\Omega_0$ as $M\uparrow \infty$. Hence we have the 
conclusion of the lemma. \hfill $\diamondsuit$
\vspace{0.5cm}

Next consider an approximation of $\Omega_0$ by smooth (or Lipschitz) open 
sets. Given any positive number $\varepsilon >0$ and any set $\Omega'$ in 
$\Omega$ we use the notation:
 $$ U_\varepsilon (\Omega'):=\{x \in \Omega~|~\exists x' \in \Omega'~
  {\rm s.t.~} |x-x'| < \varepsilon\}~~(\varepsilon{\rm -neighborhood~ of~}
   \Omega').$$

\noindent
{\bf Lemma 4.2.} {\it Let 
$\varepsilon >0$ be any small number. Then there is
a bounded open set $\omega_\varepsilon$, with a Lipschitz boundary,
in $\Omega_0$ such that
 $$ \overline {U_\varepsilon(\omega_\varepsilon)} \subset \Omega_0,$$
and
 $$ \forall x \in \Omega_0-\omega_\varepsilon,~\exists x' \in
   \Omega-\Omega_0~{\it such~that~}|x-x'|\le c_o\varepsilon,
            $$
where $c_0>1$ is a positive constant independent of $\varepsilon$, 
and moreover, $\omega_\varepsilon$ is increasing as
$\varepsilon \downarrow 0$ and
  $$ \bigcup_{\varepsilon >0} \omega_\varepsilon =\Omega_0,~~{\it hence~meas}
   (\Omega_0-\omega_\varepsilon) \to 0~{\it as~}\varepsilon \to 0.
       $$
   } 

The proof of Lemma 4.2 is elementary, so it is omitted.\vspace{0.5cm}

\noindent
{\bf Lemma 4.3.} {\it Let $\varepsilon >0$ and $\omega_\varepsilon$ be the same
as in Lemma 4.2. Then there is a cut-off function $\alpha_\varepsilon \in 
{\cal D}(\Omega)$ such that
  $$ 0\le \alpha_\varepsilon \le 1~~{\it on ~}\Omega,
       ~~{\rm supp}(\alpha_\varepsilon) \subset \Omega_0,\eqno{(4.1)}$$
  $$ \alpha_\varepsilon =1~~{\it on ~}\omega_\varepsilon,~~  
     \alpha_\varepsilon =0~~{\it on ~}\Omega-\Omega_0, \eqno{(4.2)} $$
  $$ |\nabla \alpha_\varepsilon| \le \frac {C(\Omega_0)}\varepsilon~~{\it on~}
    \Omega,\eqno{(4.3)}$$
where $C(\Omega_0)$ is a positive constant depending on $\Omega_0$ (but 
independent of $\varepsilon$). }\vspace{0.3cm}

\noindent
{\bf Proof.} As $\alpha_\varepsilon$ we can choose the convolution
  $$  \rho_{\frac \varepsilon 4}*\chi_{{\overline {U_{\frac \varepsilon 2}
                     (\omega_\varepsilon)}}}(x) := \int_\Omega 
  \rho_{\frac{\varepsilon}4}(x-y)\chi_{{\overline {U_{\frac \varepsilon 2}
                     (\omega_\varepsilon)}}}(y)dy,~~x \in \Omega,$$
where 
$\chi_{{\overline {U_{\frac \varepsilon 2}(\omega_\varepsilon)}}}$ is 
the characteristic function of $\overline {U_{\frac \varepsilon 2}
(\omega_\varepsilon)}$ and $\rho_{\frac \varepsilon 4}$ is the 3-D 
mollifier with
support $|x| \le \frac {\varepsilon}4$.  
The proof is elementary. See [2; {\bf 4.19}] for the statement.
\hfill $\diamondsuit$ \vspace{0.5cm}

\noindent
{\bf Lemma 4.4.} {\it Let $\{{\bfv}_n\}$ be a bounded
sequence in ${\bfH}_\sigma$ such that 
 $$|{\bfv}_n(x)| \le \gamma_n(x),~{\it ~a.e.~} x \in \Omega,~
         \forall n, \eqno{(4.4)}$$
and
$$ \int_\Omega ({\bfv}_n-{\bfv}_m)\cdot {\bfz}dx  \to 0
    ~~~{\it as~}n,~m \to \infty,  \eqno{(4.5)} $$
for any ${\bfz} \in {\bfW}_\sigma$ with supp$(|{\bfz}|) \subset \Omega_0$.
Then, there is ${\bfv} \in {\bfH}_\sigma$ such that
${\bfv}_n \to {\bfv}$ weakly in $L^2(\Omega)^3$.} \vspace{0.5cm}

\noindent
{\bf Proof.} Let $M$ be any positive number and $\bar q \in H^1(\Omega_0)$ with
$|\bar q(x)| \le M$ for a.e. $x \in \Omega_0$. First we shall show that 
$$ \lim_{n,m \to \infty} \int_{\Omega_0}({\bfv}_n-{\bfv}_m)\cdot \nabla 
   \bar q dx =0. \eqno{(4.6)}$$
For each small $\varepsilon >0$ we take $\omega_\varepsilon$ with 
$\overline{U_\varepsilon(\omega_\varepsilon)}\subset \Omega_0$
and a cut-off function $\alpha_\varepsilon$ as in Lemmas 4.2 and 4.3. 
Since $\alpha_\varepsilon q \in H^1_0(\Omega)$, we observe with the help of 
estimates in Lemmas 4.2, 4.3 and (4.4) that
\begin{eqnarray*}
&& \left| \int_{\omega_\varepsilon} ({\bfv}_n-{\bfv}_m)
     \cdot  \nabla \bar q dx \right| 
 = \left |\int_{\omega_\varepsilon} 
  ({\bfv}_n-{\bfv}_m)\cdot 
   \nabla (\alpha_\varepsilon \bar q) dx \right| \\
&=& \left|\int_\Omega ({\bfv}_n-{\bfv}_m)\cdot \nabla 
(\alpha_\varepsilon \bar q) dx \right. 
-\left.\int_{\Omega-\omega_\varepsilon} 
 ({\bfv}_n-{\bfv}_m) \cdot \nabla (\alpha_\varepsilon \bar q) dx \right | \\
&=& \left |0-\int_{\Omega-\omega_\varepsilon} ({\bfv}_n
    -{\bfv}_m)\cdot \nabla(\alpha_\varepsilon \bar q) dx \right| \\
 &=& \left |\int_{\Omega_0-\omega_\varepsilon} ({\bfv}_n
    -{\bfv}_m)\cdot \nabla(\alpha_\varepsilon \bar q) dx \right| \\
&\le& 
    \int_{\Omega_0-\omega_\varepsilon} (\gamma_n+\gamma_m)
   |\nabla \bar q|dx 
   +\int_{\Omega_0-\omega_\varepsilon} 
   (\gamma_n+\gamma_m) M \frac {C(\Omega_0)}\varepsilon dx.
 \end{eqnarray*}
In the last inequality we note from Lemma 4.2 that for any 
$x \in \Omega_0-\omega_\varepsilon$ 
there is $x' \in \Omega-\Omega_0$ such that $|x-x'| \le c_0\varepsilon$,
so that by the Lipschitz
continuity of $\gamma$ we have $\gamma(x)=|\gamma(x)-\gamma(x')| 
\le L_\gamma c_0\varepsilon$,
where $L_\gamma=L_\gamma(\kappa,t_0)$, with a small $\kappa>0$, is the 
Lipschitz constant of 
$\gamma=\gamma(\cdot,t_0)$ in a neighborhood of $\Omega-\Omega_0$.
Therefore, by (4.1)-(4.3),
\begin{eqnarray*}
 && \lim_{n,m \to \infty} \int_{\Omega_0-\omega_\varepsilon} 
   (\gamma_n+\gamma_m)\cdot |\nabla \bar q|dx \\
  &=& \int_{\Omega_0-\omega_\varepsilon} 2\gamma|\nabla \bar q|dx 
  \le  2L_\gamma c_0\varepsilon \int_{\Omega_0-\omega_\varepsilon}
    |\nabla \bar q|dx
    \to 0~{\rm as~} \varepsilon \to 0
  \end{eqnarray*}
and
  \begin{eqnarray*}
 &&\lim_{n,m \to \infty} \int_{\Omega_0-\omega_\varepsilon} (\gamma_n+\gamma_m)
  \cdot M \frac {C(\Omega_0)}\varepsilon dx \\
  &=& \int_{\Omega_0-\omega_\varepsilon} 2\gamma
         M \frac {C(\Omega_0)}\varepsilon dx 
  \le  2L_\gamma c_0M C(\Omega_0) \cdot{\rm meas}(\Omega_0-\omega_\varepsilon) 
   \to 0~{\rm as~}
  \varepsilon \to 0.
  \end{eqnarray*}
Thus we have (4.6). 

Here we note that (4.5) holds for every ${\bfz} \in {\bfH}_\sigma(\Omega_0)$,
since ${\bfW}_\sigma(\Omega_0)$ is dense in ${\bfH}_\sigma(\Omega_0)$.
Next, let ${\bfz}$ be any function in $L^2(\Omega_0)^3$. 
We denote the Helmholtz decomposition of ${\bfz}$ in 
$L^2(\Omega_0)^3$
by ${\bfz}:=\tilde{\bfz}+\nabla q$, $\tilde {\bfz}\in {\bfH}_\sigma(\Omega_0),
~q \in L^2_{loc}(\Omega_0)$ with $\nabla q \in L^2(\Omega_0)^3$. Given large 
constant $M>0$, we put
 $$ q^M(x)=(q(x)\land M)\lor (-M),~~x\in \Omega_0.$$
By Lemma 4.1, $q^M \in H^1(\Omega_0)$ and $\nabla q^M \to \nabla q$ in 
$L^2(\Omega_0(t))^3$ as $M\uparrow \infty$. Then
it follows  that
 \begin{eqnarray*}
 &&\left|\int_{\Omega_0} ({\bfv}_n-{\bfv}_m)\cdot 
{\bfz} dx \right| \\
&=& \left|\int_{\Omega_0} ({\bfv}_n-{\bfv}_m)\cdot 
(\tilde {\bfz}+\nabla q) dx \right | \\
 &\le &\left |\int_{\Omega_0} ({\bfv}_n-{\bfv}_m)
    \cdot \tilde{\bfz} dx \right | 
   + \left |\int_{\Omega_0} ({\bfv}_n-{\bfv}_m)
    \cdot \nabla q^M dx\right | \\
   & & + 2M'|\nabla q-\nabla q^M|_{L^2(\Omega_0)^3} ,
    \end{eqnarray*}
where $M'=\sup_n |{\bfv}_n|_{0,2}$. 

We derive from (4.5) and (4.6) with $\bar q=q^M$ that
   $$ \limsup_{n,m\to \infty} 
     \left|\int_{\Omega_0} ({\bfv}_n-{\bfv}_m)\cdot 
{\bfz} dx \right| \le 2M'|\nabla q-\nabla q^M|_{L^2(\Omega_0)^3}. 
\eqno{(4.7)}$$
Therefore, the right hand side of (4.7) converges to $0$ as 
$M\uparrow \infty$,  so that 
${\bfv}_n-{\bfv}_m \to 0$ weakly in 
$L^2(\Omega_0)^3$ as $n,m\to \infty$. This implies that ${\bfv}_n \to
{\bfv}$ in $L^2(\Omega)^3$ as well as in ${\bfH}_\sigma$, since 
${\bfv}_n \to 0$ uniformly on $\Omega-\Omega_0$ by (4.4).
\hfill $\diamondsuit$
\vspace{0.5cm}

\noindent
{\bf 4.2. Proof of Theorem 3.2 (Degenerate case)}\vspace{0.3cm} 

The main idea for the proof of Theorem 3.2 is found in [20],
but it will be here repeated for the
completeness under the additional condition (3.6).

In the rest of this paper we use the following notation:
  $$\Omega_\kappa(t):=\{x \in \Omega~|~\gamma(x,t)>
    \kappa\},~~\forall \kappa\geq 0,~\forall t \in [0,T],$$
  $$ E_0:=\{t \in [0,T]~|~\Omega_0(t) \ne \emptyset\};$$
$E_0:=\bigcup_{\ell=1}^\infty E_{\ell}$ is relatively open in $[0,T]$, where
$E_{\ell}$ is any connected component of $E_0$ . Also, we put
$$ \hat Q_J(\gamma>\kappa):=
   \bigcup_{t\in J}\Omega_\kappa(t)\times \{t\}, $$
for any subinterval $J$ of $[0,T]$. When $J=[0,T]$, $\hat Q_{[0,T]}(\gamma>0) 
=\hat Q(\gamma>0)$. \vspace{0.5cm}

We consider an exhaustion of the set $\hat Q(\gamma>0)$ by means of 
4-dimensional rectangulars (parallel to the (x,t)-coordinate axis) 
in $\hat Q(\gamma>0)$.

For each $\ell$ we observe that $\hat Q_{E_\ell}(\gamma>0)$ is a countable 
union of sets $\Omega^m_i \times J^m_i$ in the form:
 $$\hat Q_{E_\ell}(\gamma>0)=\bigcup_{m=1}^\infty \bigcup_{i=1}^{P_m}
   \Omega^m_i \times J^m_i,\eqno{(4.8)}$$
where
\begin{description}
\item{(a)} $E_\ell=\bigcup_{m=1}^\infty J^m$ with $J^m:= \bigcup_{i=1}^{P_m} 
J^m_i$ such that
\begin{itemize}
\item each $J^m_i$ is an 
interval of the form 
$[T_m,T'_m]$ or $(T_m,T'_m]$ or $[T_m,T'_m)$ or $(T_m,T'_m)$,
\item $J^m$ is a direct sum of $J^m_i$ and $J^m$ is
increasing in $m$, namely $J^m \subset J^{m+1}$,
\end{itemize}
\item{(b)} $\Omega^m_i$ is a smooth open set in $\Omega$ such that
\begin{itemize} 
\item $\overline{\Omega^m_i}\subset \Omega_0(t_{m,i})$ for some 
$t_{m,i} \in J^m_i$, hence ${\rm dist}(\Omega^m_i, \Omega-\Omega_0(t_{m,i}))
 >0$.
\item $\Omega^m_i$ is increasing in $m$ in the sense that $\Omega^m_i \subset
\Omega^{m+1}_j$ if $ t \in J^m_i \cap J^{m+1}_j$ for $1 \le j \le P_{m+1}$.
\end{itemize}
\end{description}
From (4.8) we see that
for each $t \in E_\ell$, 
  $${\rm meas}(\Omega_0(t) - \Omega^m_i) \to 0~~{\rm as~}m \to \infty,
     ~t \in J^m_i.$$ 

\begin{center}
\includegraphics[width=18cm, height=10cm]{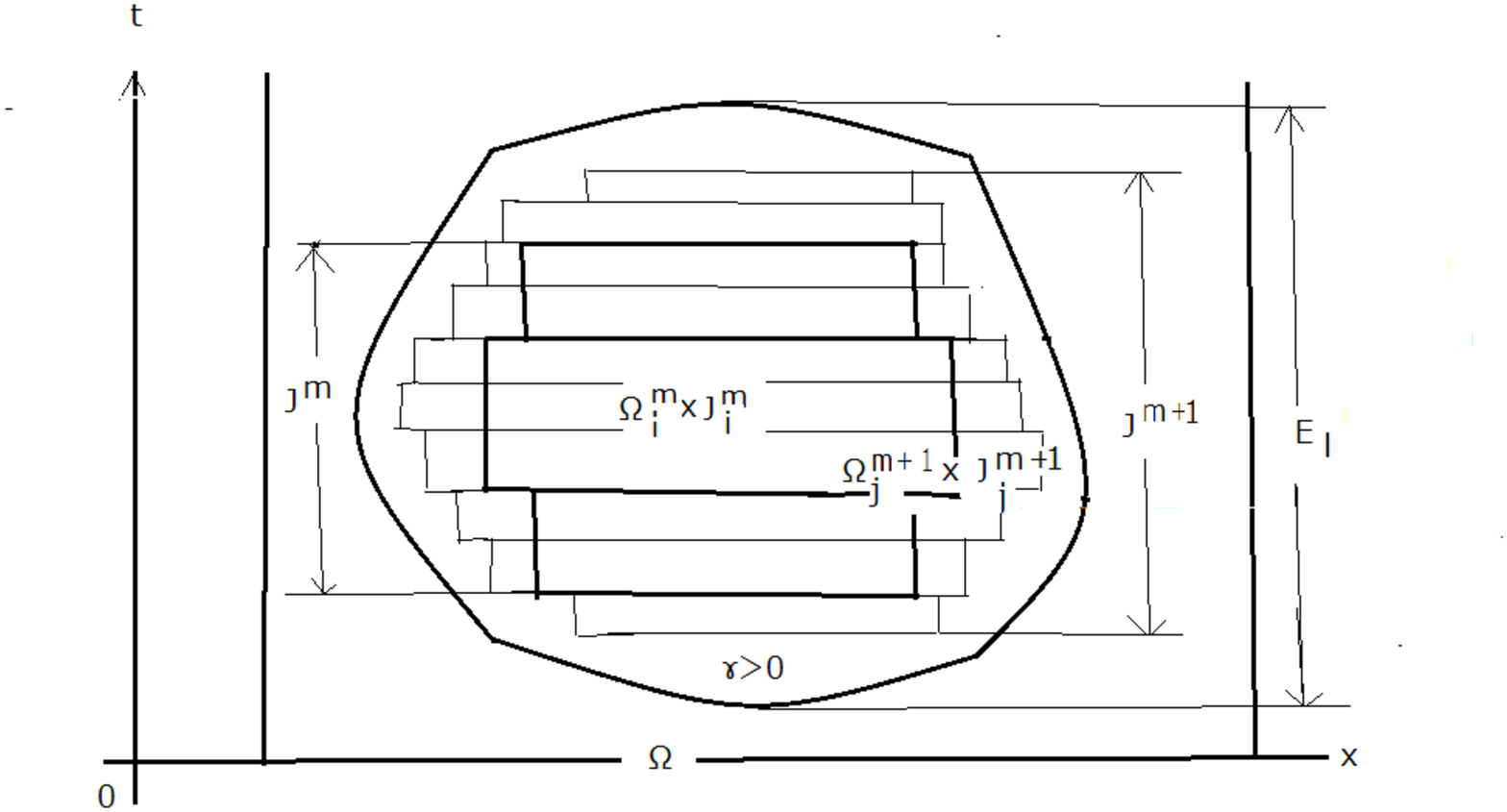}
\end{center}

Let $\gamma_n$ be the same as in the previous subsection:
 $$ \gamma_n:=(\gamma\lor \delta_n)\land N_n, ~~\delta_n\downarrow 0,~
    N_n \uparrow \infty.$$
By virtue of 
Theorem 3.1, problem $NS^1(\gamma_n;{\bfg},{\bfv}_0)$ has a solution
${\bfv}_n$ satisfying uniform estimate (cf. (3.13)):
 $$ \sup_{t \in [0,T]}|{\bfv}_n(t)|^2_{0,2}
  +\nu |{\bfv}_n|^2_{L^2(0,T;{\bfV}_\sigma)} \le M_0({\bfv}_0,{\bfg})=:M_0.
    \eqno{(4.9)}$$
From (4.9) we see that 
$\{{\bfv}_n\}$ is bounded in $L^\infty(0,T;{\bfH}_\sigma)\cap
L^2(0,T;{\bfV}_\sigma)$. Now choose a subsequence of
$\{{\bfv}_n\}$, denoted by the same notation for simplicity, such that
  $${\bfv}_n \to {\bfv}~{\rm weakly}^*~{\rm in~}
  L^\infty(0,T;{\bfH}_\sigma)~~{\rm and~weakly~ in~}L^2(0,T;{\bfV}_\sigma),
   \eqno{(4.10)}$$
for some ${\bfv} \in L^\infty(0,T;{\bfH}_\sigma)\cap L^2(0,T;{\bfV}_\sigma)$.
We note that 
${\bfv}_n \in \boldsymbol{\mathcal K}^1(\gamma_n)$ , i.e.
$|{\bfu}_n| \le \gamma_n$ a.e. on $Q$.
Since $\gamma_n \to \gamma$ uniformly on $Q$ (in the extended sense (3.7)), it 
follows that
  $$|{\bfv}| \le \gamma~{\rm a.e.~on~}Q,~~i.e.~{\bfv} \in 
     \boldsymbol{\cal K}^1(\gamma).$$
Besides, each ${\bfv}_n$ satisfies the following variational inequality:
   $$\int_0^t(\boldsymbol {\mathcal \xi}'(\tau), {\bfv}_n(\tau)
    -\boldsymbol {\mathcal \xi}(\tau))_\sigma d\tau
    +\nu \int_0^t \langle {\bfF}{\bfv}_n(\tau), {\bfv}_n(\tau)- 
    \boldsymbol {\mathcal \xi}(\tau)\rangle_\sigma d\tau ~~~~~~~~~$$
 $$  ~~+ \int_0^t \langle {\bfG}({\bfv}_n, {\bfv}_n),
    {\bfv}_n-\boldsymbol {\mathcal \xi}\rangle_\sigma d\tau 
     +\frac 12 |{\bfv}_n(t)-\boldsymbol {\mathcal \xi}(t)|^2_{0,2}  
   \eqno{(4.11)}$$
  $$ ~~~~~~~~~~~~~~~\le \int_0^t \langle{\bfg},{\bfv}_n
      -\boldsymbol {\mathcal \xi}\rangle_\sigma d\tau
     + \frac 12 |{\bfv}_0-\boldsymbol {\mathcal \xi}(0)|^2_{0,2},
~~\forall t \in [0,T],
    ~\forall \boldsymbol {\mathcal \xi} \in 
      \boldsymbol{\mathcal K}^1_0(\gamma_n).~~~~~~~~~~ $$
Putting
$${\bff}_n:={\bfg}-\nu{\bfF}{\bfv}_n-{\bfG}({\bfv}_n,{\bfv}_n),\eqno{(4.12)}$$
we observe that ${\bff}_n \in L^2(0,T;{\bfV}_\sigma^*)$ and
$\{{\bfG}({\bfv}_n,{\bfv}_n)\}$ is bounded in $L^1(0,T;{\bfW}_\sigma^*)$, 
hence $\{{\bff}_n\}$ is bounded in $L^1(0,T;{\bfW}^*_\sigma)$.
\vspace{0.5cm}

Now, for each index $\{m,i\}$ in (4.8) consider the sequence of approximate 
solutions ${\bfv}_n$ on the cylindrical open set $\Omega^m_i \times J^m_i$. 
For simplicity we write $\Omega'$ for $\Omega^m_i$ and
write $J~(=[T_1,T'_1])$ for $J^m_i$, assuming $J^m_i$ is a closed interval 
$[T_{m,i},T'_{m,i}]$; any other case can be similarly handled.
With these notations, consider the spaces ${\bfV}_\sigma(\Omega'),
~{\bfH}_\sigma(\Omega')$ and ${\bfW}_\sigma(\Omega')$ built on $\Omega'$;
  $$ {\bfW}_\sigma(\Omega') \subset {\bfV}_\sigma(\Omega'), ~~
  {\bfV}_\sigma(\Omega') \subset {\bfH}_\sigma(\Omega') 
     \subset {\bfW}^*_\sigma(\Omega') $$
with dense and compact embeddings. 

By the Helmholtz decomposition
of ${\bfv}_n(\cdot,t)$ in $L^2(\Omega')^3$,
 $$ {\bfv}_n(x, t)=\tilde{\bfv}_n(x,t)+ \nabla q_n(x,t),
 ~~\tilde{\bfv}_n(\cdot,t) \in {\bfH}_\sigma(\Omega')\cap W^{1,2}(\Omega')^3,
    ~q_n(\cdot,t) \in H^1(\Omega');$$
we know (cf. [34; Chaper 1]) that 
$|\tilde {\bfv}_n|_{W^{1,2}(\Omega')^3} \le C(\Omega')|{\bfv}_n|_
  {W^{1,2}(\Omega')^3} $ with some positive constant $C(\Omega')$
depending only on Lipschitz open set $\Omega'$. 
It is easy to see from (4.9) that
  $$ \sup_{t \in [0,T]}|\tilde{\bfv}_n(t)|^2_{{\bfH}_\sigma(\Omega')}+ 
    \nu |\tilde{\bfv}_n|^2_{L^2(0,T;W^{1,2}(\Omega')^3)} \le M_1({\bfu}_0,
    {\bfg})=:M_1. \eqno{(4.13)}$$
where $M_1$ is a positive constant independent of $n$.
We can regard ${\bff}_n(t)$ given by (4.12) as a linear continuous 
functional, denoted by $\tilde{\bff}_n$, on ${\bfW}_\sigma(\Omega')$, 
by putting
 \begin{eqnarray*}
\langle \tilde{\bff}_n(t), \boldsymbol {\mathcal \xi}\rangle_\sigma
&:=&
    \int_{\Omega'}{\bfg}(x,t) \cdot \boldsymbol {\mathcal \xi}(x)dx
  - \nu \int_{\Omega'} \nabla{\bfv}_n(x,t)\cdot 
     \nabla \boldsymbol {\mathcal \xi}(x)dx \\
   &&~~~+\int_{\Omega'} ({\bfv}_n(t)\cdot \nabla){\bfv}_n(t)\cdot 
      \boldsymbol {\mathcal \xi}(x) dx,~~~\forall 
     \boldsymbol {\mathcal \xi} \in {\bfW}_\sigma(\Omega').
 \end{eqnarray*} 
The estimate (4.13) shows that $\{\tilde{\bff}_n\}$ is
bounded in $L^1(T_1,T'_1;{\bfW}^*_\sigma(\Omega'))$ and 
  $$\int_{T_1}^{T'_1} \langle {\bff}_n, {\bfv}_n \rangle_\sigma dt\le 
\int_{T_1}^{T'_1} \langle {\bfg}, {\bfv}_n \rangle_\sigma dt $$
is bounded above; namely for a positive constant $M'_1$ it holds that
 $$ |\tilde{\bff}_n|_{L^1(T_1,T'_1;{\bfW}^*_\sigma(\Omega'))} \le M'_1,~
  ~~\forall n.
  $$

\noindent
{\bf Lemma 4.5.} {\it Let $\Omega'\times [T_1,T'_1]:=\Omega^m_i\times [T_{m,i},
T'_{m,i}]$ be as above, and let $C_1$
be a positive constant satisfying $|{\bfz}|_{C(\overline \Omega)^3} \le C_1
 |{\bfz}|_{{\bfW}_\sigma}$ for all ${\bfz}\in {\bfW}_\sigma$. Then 
$\tilde{\bfv}_n$ is of bounded variation as a function from $[T_1,T'_1]$ into
${\bfW}^*_\sigma(\Omega')$ and its total variation is estimated by:
$$ 
  {\rm Var}_{{\bfW}^*_\sigma(\Omega')}(\tilde{\bfv}_n)\\[0.3cm]
   \le M'_1 
   +\frac{C_1M_0^{\frac 12}}{\kappa_{m,i}}\int_{T_1}^{T'_1} |{\bfg}|_{0,2}
    dt
   +\frac {C_1M_0}{2\kappa_{m,i}},~~\forall {\rm large~} n,
  \eqno{(4.14)} $$
with the constant $\kappa_{m,i}$ in condition (b) and $M_0$ the same 
constant as in (3.13). }

\noindent
{\bf Proof.} We first show that there is a constant $\kappa'>0$ such that
    $$ \kappa'B_{{\bfW}_\sigma(\Omega')}(0) \subset K^1(\gamma_n;t),~\forall
     t \in [T_1,T'_1],~\forall {\rm large}~n. \eqno{(4.15)}$$
In fact, since $\overline{\Omega'}=\overline {\Omega^m_i}\subset 
\Omega_{\kappa_{m,i}}$ and $\gamma_n \to \gamma$ uniformly on $\Omega$, 
we see that $\gamma >\kappa_{m,i}$ on 
$\overline{\Omega'}$ and hence $\gamma_n > \kappa_{m,i}$ on $\overline \Omega'$
 for all large $n$. 
Therefore, if ${\bfz} \in B_{{\bfW}_\sigma(\Omega')}(0)$, then 
 $$\frac{\kappa_{m,i}}{C_1}|{\bfz}|_{C(\overline {\Omega'})^3} \le 
    \kappa_{m,i} |{\bfz}|_{{\bfW}_\sigma(\Omega')} \le
  \kappa_{m,i} \le \gamma_n(\cdot,t)~{\rm a.e. ~on~} \Omega',
   ~\forall t\in [T_1,T_1'],~\forall ~{\rm large~}n.$$ 
Thus we have (4.15) with $\kappa'=\frac{\kappa_{m,i}}{C_1}$.

Next, let $\boldsymbol {\mathcal \xi}$ be any element in 
$C^1_0(T_1,T'_1;{\bfW}_\sigma(\Omega'))$ with 
$|\boldsymbol {\mathcal \xi}|_{C([T_1,T'_1];{\bfW}_\sigma(\Omega'))} \le 1$.
 Then,
 $\hat {\boldsymbol {\mathcal \xi}}:=\pm \kappa' 
   {\boldsymbol {\mathcal \xi}}$ is a possible test function for (4.11), 
hence we get that
 $$
  \int_{T_1}^{T'_1} \langle \hat{\boldsymbol {\mathcal \xi}}', 
   \tilde{\bfv}_n - \hat{\boldsymbol {\mathcal \xi}}\rangle_\sigma dt 
   \le \int_{T_1}^{T'_1} \langle {\bff}_n, 
   {\bfv}_n -\hat {\boldsymbol {\mathcal \xi}}\rangle_\sigma dt
   +\frac 12|{\bfu}_n(T_1)|^2_{0,2}-\frac 12 |{\bfu}_n(T'_1)|^2_{0,2},
 $$
whence
 \begin{eqnarray*}
 && \int_{T_1}^{T'_1} \langle {\boldsymbol {\mathcal \xi}}', 
   \tilde{\bfv}_n \rangle_\sigma dt \\[0.3cm]
   &\le& - \int_{T_1}^{T'_1} \langle \tilde {\bff}_n, 
   {\boldsymbol {\mathcal \xi}}\rangle_\sigma dt
   +\frac 1{\kappa'}\int_{T_1}^{T'_1}({\bfg},{\bfv}_n)_\sigma dt
   +\frac 1{2\kappa'}|{\bfv}_n(T_1)|^2_{0,2}
       -\frac 1{2\kappa'} |{\bfv}_n(T'_1)|^2_{0,2} \\[0.3cm]
   &\le& \int_{T_1}^{T'_1} |\tilde {\bff}_n|_{{\bfW}^*_\sigma(\Omega')}dt
   +\frac {M_0^{\frac 12}}{\kappa'} \int_{T_1}^{T'_1}|{\bfg}|_{0,2}dt
   +\frac 1{2\kappa'}|{\bfv}_n(T_1)|^2_{0,2}
        -\frac 1{2\kappa'} |{\bfv}_n(T'_1)|^2_{0,2}.   
 \end{eqnarray*}
Consequently, we have (4.14). \hfill $\diamondsuit$
\vspace{0.5cm}

\noindent
{\bf Corollary 4.1.} {\it Under the same assumptions and notation as in 
Lemma 4.5, for each $\{m,i\}$ there are a subsequence 
$\{\tilde{\bfv}_{n_{k(m,i)}}\}_{k=1}^\infty$ of 
$\{\tilde{\bfv}_n\}$ (depending on $\Omega^m_i$) and a function  
$\tilde{\bfv}_{m,i}$ in $L^\infty(J^m_i;{\bfH}_\sigma(\Omega^m_i))
\cap L^2(J^m_i; W^{1,2}(\Omega^m_i)^3)$ such that 
  $$\tilde {\bfv}_{n_{k(m,i)}}(t) \to \tilde {\bfv}_{m,i}(t)~
  {\it weakly~in~}
     L^2(\Omega^m_i)^3,~~\forall t \in J^m_i~
     ({\it as~}k\to \infty), \eqno{(4.16)}$$
hence
  $$ \int_{\Omega^m_i} ({\bfv}_{n_{k(m,i)}}(\cdot,t)-{\bfv}_{n_{j(m,i)}}
 (\cdot,t))\cdot {\bfz}dx \to 0,~\forall {\bfz}\in {\bfH}_\sigma(\Omega^m_i),~
      \forall t \in J^m_i, ~{\it as~} k,~j \to \infty. \eqno{(4.17)}$$
 }

\noindent
{\bf Proof.} According to Lemma 4.5 and Lemma 2.1, there is a subsequence
$\{\tilde{\bfv}_{n_{k(m,i)}}\}$ of $\{\tilde{\bfv}_n\}$ and a function 
$\tilde{\bfu}_{m,i}\in L^\infty(J^m_i;{\bfH}_\sigma(\Omega^m_i))\cap 
L^2(J^m_i; W^{1,2}(\Omega^m_i)^3)$
such that
 $$\tilde {\bfv}_{n_{k(m,i)}}(t) \to \tilde {\bfv}_{m,i}(t)~{\rm weakly~in~}
     {\bfH}_\sigma(\Omega^m_i),~\forall t \in J^m_i~
     ({\rm as~}k\to \infty).\eqno{(4.18)}$$
Here, for any function ${\bfz} \in L^2(\Omega')^3$ we use the Helmholtz 
decomposition 
  $${\bfz}=\tilde{\bfz}+\nabla q,~\tilde {\bfz} \in 
{\bfH}_\sigma(\Omega^m_i),~q \in H^1(\Omega^m_i),$$
 to see by (4.18) 
 \begin{eqnarray*}
  \int_{\Omega^m_i}\tilde{\bfv}_{n_{k(m.i)}}(\cdot,t)\cdot {\bfz}dx &=&
 \int_{\Omega^m_i}\tilde{\bfv}_{n_{k(m,i)}}(\cdot,t)\cdot 
(\tilde{\bfz}+\nabla q)dx
  \\[0.3cm]
 &=& \int_{\Omega^m_i}\tilde{\bfv}_{n_{k(m,i)}}(\cdot,t)\cdot \tilde{\bfz}dx 
       \to \int_{\Omega^m_i}\tilde{\bfv}_{m,i}(\cdot,t)\cdot 
      \tilde{\bfz}dx\\[0.3cm]
 &=& \int_{\Omega^m_i}\tilde{\bfv}_{m,i}(\cdot,t)\cdot 
  (\tilde{\bfz}+\nabla q) dx
      =\int_{\Omega^m_i}\tilde{\bfv}_{m,i}(\cdot,t)\cdot {\bfz} dx;
 \end{eqnarray*}
hence (4.16) holds.
The convergence (4.17)
follows, since $\int_{\Omega^m_i} {\bfv}_{n_{k(m,i)}}\cdot {\bfz} dx= 
\int_{\Omega^m_i} \tilde{\bfv}_{n_{k(m,i)}}\cdot {\bfz} dx$ for ${\bfz} \in 
{\bfH}_\sigma(\Omega^m_i)$.
 \hfill $\diamondsuit$\vspace{0.5cm}

\noindent
{\bf Corollary 4.2.} {\it Under the same assumptions and notation as in 
Corollary 4.1, there is a subsequence $\{{\bfv}_{n_{k(m)}}\}_{k=1}^\infty$ of 
$\{{\bfv}_n\}$, depending only on $m$, such that
 $$ \int_\Omega ({\bfv}_{n_{k(m)}}(\cdot,t) - {\bfv}_{n_{j(m)}}(\cdot,t))\cdot 
       {\bfz}dx \to 0,~~\forall {\bfz} \in {\bfH}_\sigma(\Omega^m_i),~
       \forall t \in J^m_i, \eqno{(4.19)}$$
as $k,~j \to \infty$.}\vspace{0.3cm}

\noindent
{\bf Proof.}
This corollary is a direct consequence of Corollary 4.1 and the 
definition ${\bfv}_{n_{k(m,i)}}$; in fact, we extract
a subsequence $\{{\bfv}_{n_{k(m,i+1)}}\}$ from $\{{\bfv}_{n_{k(m,i)}}\}$ so as
to satisfy (4.17) with $i$ replaced by $i+1$, 
repeatedly for $i=1,2,\cdots, N_m-1$. As a result, we get a subsequence 
$\{{\bfv}_{n_{k(m)}}\}$ for which (4.19) holds. \hfill $\diamondsuit$
\vspace{0.5cm}

\noindent
{\bf Lemma 4.6.} {\it Under the same assumptions and notation as in 
Corollary 4.2,
there is a subsequence $\{{\bfv}_{n_m}\}_{m=1}^\infty$ with a function 
$\bar {\bfv}: E_\ell \to {\bfH}_\sigma$ such that
  $${\bfv}_{n_m}(t) \to \bar{\bfv}(t) ~{\it weakly~in~}L^2(\Omega)^3,~~
    \forall t \in E_\ell,~(as~m \to \infty). 
    \eqno{(4.20)}$$ }
\noindent
{\bf Proof.} We make use of the subsequences $\{{\bfv}_{n_{k(m)}}\}$ 
constructed in Corollary 4.2. In the table of these subsequences:
  \begin{eqnarray*}
 && {\bfv}_{n_{1(1)}}~~{\bfv}_{n_{2(1)}}~~{\bfv}_{n_{3(1)}}~~\cdots\cdots\cdots
    \cdots\cdots
    ~~{\bfv}_{n_{m(1)}}~~ \cdots \\ 
 && {\bfv}_{n_{1(2)}}~~{\bfv}_{n_{2(2)}}~~{\bfv}_{n_{3(2)}}~~\cdots\cdots\cdots
    \cdots\cdots
     ~~{\bfv}_{n_{m(2)}}~~\cdots \\
 && {\bfv}_{n_{1(3)}}~~{\bfv}_{n_{2(3)}}~~{\bfv}_{n_{3(3)}}~~\cdots\cdots\cdots
     \cdots\cdots
     ~~{\bfv}_{n_{m(3)}}~~\cdots \\[0.5cm]
 && ~~~~~~~~~~~~~~\cdots\cdots\cdots\cdots \\
 && ~~~~~~~~~~~~~~\cdots\cdots\cdots\cdots \\[0.5cm]
 && {\bfv}_{n_{1(m)}}~~{\bfv}_{n_{2(m)}}~~{\bfu}_{n_{3(m)}}
     ~~\cdots\cdots\cdots\cdots\cdot~~{\bfu}_{n_{m(m)}}~~\cdots
      \\[0.5cm]
 &&~~~~~~~~~~~~~~\cdots\cdots\cdots\cdots 
  \end{eqnarray*}
we pic up the diagonal functions ${\bfv}_{n_{m(m)}}=:{\bfv}_{n_m}$ and 
consider the subsequence $\{{\bfv}_{n_m}\}_{m=1}^\infty$ of $\{{\bfv}_n\}$.

We shall show below that this is a required one. 
Let $\boldsymbol {\mathcal \xi}$ be any function ${\bfW}_\sigma$ such that
$K:={\rm supp}(|\boldsymbol {\mathcal \xi}|) \subset \Omega_0(t),
~t \in E_\ell.$ Since $K$ is compact in $\Omega_0(t)$, there is a positive 
number $\kappa$ such that $\gamma(x,t) >\kappa$ for all $x \in K$. Therefore
by condition (a), $K \subset \Omega^m_i \subset \Omega_0(t)$ for a large $m$ 
and some $i$ with $t \in J^m_i$, which implies by Corollary 4.2 that
  $$ 0=\lim_{k,j\to \infty} \int_\Omega ({\bfv}_{n_{k(m)}}(\cdot,t)
     -{\bfv}_{n_{j(m)}}(\cdot,t))\cdot\boldsymbol {\mathcal \xi} dx
      =\lim_{m,m'\to \infty} \int_\Omega ({\bfv}_{n_m}(\cdot,t)
     -{\bfv}_{n_{m'}}(\cdot,t))\cdot {\boldsymbol {\mathcal \xi}} dx
   $$
Applying Lemma 4.4 to the sequence $\{{\bfv}_{n_m}\}$, we conclude that
there is $\bar{\bfv}(t) \in L^2(\Omega)^3$ such that ${\bfv}_{n_m}(t) \to
\bar{\bfv}(t)$ weakly in $L^2(\Omega)^3$ for each $t \in E_\ell$ 
as $m\to \infty$. Thus we have (4.20). 
\hfill $\diamondsuit$ \vspace{0.5cm}

\noindent
{\bf Corollary 4.3.} {\it Under the same assumptions and notation as in 
Lemma 4.6, we put 
  $$ \hat{\bfv}(x,t):=\left \{
       \begin{array}{ll}
         \bar{\bfv}(x,t),~~~& (x,t) \in \Omega\times E_\ell\\[0.2cm]
         0,   &{\it otherwise},
       \end{array} \right. $$
Then
  $${\bfv}_{n_m} \to \hat {\bfv}~~{\it in~}L^2(0,T;{\bfH}_\sigma). $$ 
}
\noindent
{\bf Proof.} From Lemmas 4.5 and 4.6 with their
corollaries adapted for every component $E_\ell$ we can construct a subsequence
$\{{\bfv}_{n_m}\}$ of $\{{\bfv}_n\}$ with a function $\bar {\bfv}: E_0 \to
{\bfH}_\sigma$ such that ${\bfv}_{n_m}(t) \to \bar {\bfv}(t)$ weakly in 
$L^2(\Omega)^3$ for all $t \in E_0$.
Besides, since ${\bfv}_{n_m}(\cdot,t) \to 0$ uniformly on 
$\Omega-\Omega_0(t)$ for 
every $t \in [0,T]-E_0$, it follows consequently that ${\bfv}_{n_m}(t) \to 
\hat {\bfv}(t)$ weakly in $L^2(\Omega)^3$ for all $t \in [0,T]$.

Now we recall a compactness lemma [27;Lemma 5.1, Chapter 1] to get
$$|{\bfv}_{n_m}(t)-{\bfv}_{n_{m'}}(t)|^2_{0,2} \le \varepsilon 
   |{\bfv}_{n_m}(t)-{\bfv}_{n_{m'}}(t)|^2_{1,2} 
      +C_\varepsilon |{\bfv}_{n_m}(t)-{\bfv}_{n_{m'}}(t)|^2_{{\bfW}^*_\sigma},
   \eqno{(4.21)}$$ 
    $$ \forall t \in [0,T],~~\forall m,~m',$$
where $\varepsilon>0$ is any positive number and $C_\varepsilon$ is a positive
constant depending only $\varepsilon$. 
Integrating the above inequality in time $t$ and letting 
$m,~m' \to \infty$ yield that
$ {\bfv}_{n_m}-{\bfv}_{n_{m'}} \to 0$ in $L^2(0,T;{\bfH}_\sigma)$ as $m,~m'\to 
\infty$, since the last term of (4.21) tends to $0$ uniformly on $[0,T]$  as 
$m,~m' \to \infty$. This shows that ${\bfv}_{n_m} \to \hat{\bfv}$ 
(strongly) in 
$L^2(0,T;{\bfH}_\sigma)$ (hence in $L^2(Q)$) as $m \to \infty$.
\hfill $\diamondsuit$ \vspace{0.5cm}

\noindent
{\bf Proof of Theorem 3.2:} Let $\hat{\bfv}$ be the same function as in 
Corpllary 4.3.
Recalling that ${\bfv}_n \to {\bfv}$ weakly in $L^2(0,T;{\bfV}_\sigma)$ and
weakly$^*$ in $L^2(0,T;{\bfH}_\sigma)$ as $n \to \infty$ (cf (4.10)), we see 
that $\hat{\bfv}={\bfv}$ a.e. on $Q$ and hence $\hat{\bfv}$ may be identified
with ${\bfv}$. In the sequel, for simplicity, let us denote ${\bfv}_{n_m}$ by 
${\bfv}_n$. 

Let $\boldsymbol {\mathcal \xi}$ be any function in 
$\boldsymbol{\cal K}^1_0(\gamma)$. Then
$\boldsymbol {\mathcal \xi}$ is a possible test function of approximate problem
$NS^1(\gamma_n; {\bfg}, {\bfv}_0)$ for all large $n$, so that
  $$\int_0^t(\boldsymbol {\mathcal \xi}'(\tau), {\bfv}_n(\tau)
   -\boldsymbol {\mathcal \xi}(\tau))_\sigma d\tau
    +\nu \int_0^t \langle {\bfF}{\bfv}_n(\tau), {\bfv}_n(\tau)- 
    \boldsymbol {\mathcal \xi}(\tau)\rangle_\sigma d\tau ~~~~~~~~~$$
 $$  ~~+ \int_0^t \langle {\bfG}({\bfv}_n, {\bfv}_n),
    {\bfv}_n-\boldsymbol {\mathcal \xi}\rangle_\sigma d\tau 
     +\frac 12 |{\bfv}_n(t)-\boldsymbol {\mathcal \xi}(t)|^2_{0,2}  
   \eqno{(4.22)}$$
  $$\le \int_0^t \langle{\bfg},{\bfv}_n-\boldsymbol {\mathcal \xi}
    \rangle_\sigma d\tau
     + \frac 12 |{\bfv}_0-\boldsymbol {\mathcal \xi}(0)|^2_{0,2},
~~\forall t \in [0,T], $$
Just as in the proof of Theorem 3.1 we obtain (cf. (3.18), (3.19))) that
   $$ \frac 12|{\bfv}(t) - \boldsymbol {\mathcal \xi}(t)|^2_{0,2}\le 
  \lim\inf_{n \to \infty}
     \frac 12|{\bfv}_n(t) - \boldsymbol {\mathcal \xi}(t)|^2_{0,2}$$
and
   $$ \lim_{n\to\infty}\int_0^t \langle {\bfG}({\bfv}_n,{\bfv}_n), 
        \boldsymbol {\mathcal \xi}\rangle d\tau
      =\int_0^t \langle {\bfG}({\bfv},{\bfv}), 
      \boldsymbol {\mathcal \xi}\rangle d\tau.$$
Hence, letting $n \to \infty$ in (4.22) gives 
$$\int_0^t(\boldsymbol {\mathcal \xi}'(\tau), {\bfv}(\tau)
    -\boldsymbol {\mathcal \xi}(\tau))_\sigma d\tau
    +\nu \int_0^t \langle {\bfF}{\bfv}(\tau), {\bfv}(\tau)- 
    \boldsymbol {\mathcal \xi}(\tau)\rangle_\sigma d\tau ~~~~~~~~~$$
 $$  ~~+ \int_0^t \langle {\bfG}({\bfv}, {\bfv}),
    {\bfv}-\boldsymbol {\mathcal \xi}\rangle_\sigma d\tau 
     +\frac 12 |{\bfv}(t)-\boldsymbol {\mathcal \xi}(t)|^2_{0,2}$$
  $$\le \int_0^t \langle{\bfg},{\bfv}-\boldsymbol {\mathcal \xi}
    \rangle_\sigma d\tau
     + \frac 12 |{\bfv}_0-\boldsymbol {\mathcal \xi}(0)|^2_{0,2},
    ~~\forall t \in [0,T]. $$
As to the other properties of ${\bfv}$, we have that
  $$ {\bfv} \in \boldsymbol{\cal K}^1(\gamma),$$
 $$ {\bfv}_n(0)={\bfv}_0 \to {\bfv}(0)~{\rm weakly~in~}{\bfH}_\sigma,
    ~{\rm hence~}{\bfv}(0)={\bfv}_0.$$
Finally we show that $t\to ({\bfv}(t),\boldsymbol {\mathcal \xi}(t))_\sigma$
is of bounded variation on $[0,T]$. It is enough to show it in the case of
${\rm supp}
(|\boldsymbol {\mathcal \xi}|) \subset \hat Q_{E_\ell}(\gamma>0)$ for some
$\ell$. In this case there is $m^*$ such 
that ${\rm supp}(|\boldsymbol {\mathcal \xi}|) \subset \bigcup_{i=1}^{P_m^*} 
\Omega^{m^*}_i \times J^{m^*}_i$. We use the Helmholtz decomposition of 
${\bfv}(t)$ in $L^2(\Omega^{m^*}_i)$ which is of the form
$ {\bfv}(t)=\tilde{\bfv}(t)+\nabla q(t),~\tilde{\bfv}(t)\in 
{\bfH}_\sigma(\Omega^{m^*}_i),~q(t)\in H^1(\Omega^{m^*}_i),$ for each 
$t \in J^{m^*}_i$. 
By the estimate (4.14) in Lemma 4.5 we have
  $$ {\rm Var}_{{\bfW}_\sigma^*(\Omega^{m^*}_i)}(\tilde{\bfv}) \le M'_1 
   +\frac{C_1M_0^{\frac 12}}{\kappa_{m,i}}\int_{T_1}^{T'_1} |{\bfg}|_{0,2}
    dt
   +\frac {C_1M_0}{2\kappa_{m,i}}.
    \eqno{(4.23)} $$
In fact, since ${\bfv}_n \to {\bfv}$ weakly in $L^2(Q)^3$, it follows that 
$\tilde{\bfv}_n \to \tilde{\bfv}$ 
in $L^2(0,T; {\bfH}_\sigma(\Omega^{m^*}_i)$. By the lower semicontinuity of
the total variation functional, (4.14) implies (4.23). Moreover, for any $i$
 we observe that
 \begin{eqnarray*}
  &&|({\bfv}(t),\boldsymbol {\mathcal \xi}(t))_\sigma 
        -({\bfv}(s),\boldsymbol {\mathcal \xi}(s))_\sigma|
  = |(\tilde {\bfv}(t),\boldsymbol {\mathcal \xi}(t))_\sigma 
        -(\tilde {\bfv}(s),\boldsymbol {\mathcal \xi}(s))_\sigma|  \\
  &\le& |\tilde{\bfv}(t)-\tilde{\bfv}(s)|_{0,2}|
    |\boldsymbol {\mathcal \xi}(t)|_{{\bfW}_\sigma}
       + |{\bfv}(s)|_{0,2}|\boldsymbol {\mathcal \xi}(s)
        -\boldsymbol {\mathcal \xi}(t)|_{{\bfW}_\sigma},
      ~~\forall s,~t \in J^{m^*}_i.
 \end{eqnarray*}
From this we see easily that the total variation of 
$t\to ({\bfv}(t),\boldsymbol {\mathcal \xi}(t))_\sigma$
on $J^{m^*}$ is bounded by 
 $${\rm const.} \left(\sum_{i=1}^{P_{m^*}}{\rm Var}_{{\bfW}^*_\sigma
   (\Omega^{m^*}_i)}(\tilde{\bfv})
         +|\boldsymbol {\mathcal \xi}'|_{L^1(0,T;{\bfW}_\sigma)}\right),
\eqno{(4.24)}$$
and (4.24) is valid for every connected component $E_\ell$. Hence 
$t\to ({\bfv}(t),\boldsymbol {\mathcal \xi}(t))_\sigma$ is of bounded
variation on $[0,T]$.
Thus ${\bfv}$ is a weak solution of 
$NS^1(\gamma;{\bfg},{\bfv}_0)$. \hfill $\diamondsuit$ \vspace{0.5cm}

\noindent
{\bf Remark 4.1.} We suppose that 
\begin{itemize}
\item for each $t \in [0,T]$ the closure of
any connected component of $\Omega-\overline{\Omega_0(t)}$
meets the boundary $\Gamma$.
\end{itemize}
In this case, the same type of lemma as Lemma 4.4
can be proved under gradient constraint $|\nabla{\bfv}|\le \gamma$. Therefore
with some modification in the proof of Theorem 3.2, an existence result can be
shown for $NS^2(\gamma;{\bfg}, {\bfv}_0)$. However, this additional geometric
assumption is too restrictive to apply this result to Stefan/Navier-Stokes 
problem discussed in the next section.\vspace{1cm}

\noindent
{\large\bf 5. Stefan/Navier-Stokes problems}\vspace{0.5cm}

In this section, let us consider Stefan/Navier-Stokes problem.
As was mentioned in the introduction, it is a system of the enthalpy
 formulation of solid-liquid phase change in fluid flows with freezing and 
melting effect.

We begin with the precise formulations in non-degenerate and
degenerate cases of obstacle function, postulating that
the region $\Omega$ is divided into three unknown time-dependent regions,
 $$\Omega= \Omega_s(t)\cup \Omega_\ell(t) \cup \Omega_m(t),$$ 
which are respectively 
called the solid, liquid and mixture (mussy) regions, and velocity constraint
depends on the order parameter of phase, namely the velocity field is 
independent of space variable $x$ on $\Omega_s(t)$, governed by Navier-Stokes equation in $\Omega_\ell(t)$ and constrained by an 
order parameter dependent obstacle function in $\Omega_m(t)$. In our model, one
of main ideas is that the dynamics of velocity field is described as a 
quasi-variational inequality of Navier-Stokes type.
\vspace{0.5cm}

\noindent
{\bf 5.1. Variational formulation of Stefan/Navier-Stokes problem} 
\vspace{0.3cm} 

First of all we give the weak variational formulations of Stefan/Navier-Stokes 
problems with constraints on the velocity and its gradient.

Let $\beta=\beta(r)$ be a non-decreasing and Lipschitz continuous function with
Lipschitz constant $L_\beta$ from ${\bf R}$ into ${\bf R}$ such that
  $$ \left\{
      \begin{array}{l}
    \displaystyle{
      \liminf_{|r|\to \infty} \frac {\beta(r)}{|r|} >0,
       ~~\beta(r)=0~{\rm for ~}r \in [0,1],}\\[0.5cm]
     \displaystyle{ \beta(r) ~{\rm is~strictly~increasing~for~}r<0~
      ~{\rm and~for~}r>1. }
      \end{array} \right. \eqno{(5.1)}$$
Let $\hat \beta:=\hat \beta(r)$ be the primitive of $\beta$ given by
$\hat \beta(r):= \int_0^r \beta(s) ds$ for all $r \in {\bf R}$. Clearly
 $$\hat \beta(0)=0,~~\beta(r)= \hat\beta'(r)
    \left(= \frac {d\hat\beta(r)}{dr}\right),~
 \forall r \in{\bf R}.$$
Next, let $\gamma=\gamma(r)$ be a non-negative and non-decreasing continuous 
function from ${\bf R}$ into $[0,\infty]$
and consider the following two cases:
\vspace{0.5cm}

\noindent
{\bf (Non-degenerate case)}
  $$ \left\{
     \begin{array}{l}
  \displaystyle{\gamma(r)\geq c_*>0 ~{\rm for~}r\in {\bf R},~
   {\rm where~} c_*~{\rm is~ a~ positive~ constant},~~~~~}\\[0.3cm]
   \displaystyle{\gamma~{\rm is~strictly~increasing~on~}[0,1),}\\[0.3cm]
    \displaystyle{ \gamma(r) \uparrow \infty~{\rm as ~}r\uparrow 1,
   ~~\gamma(r)=\infty ~{\rm for~}r \geq 1.}
     \end{array} \right. \eqno{(5.2)}$$

\noindent
{\bf (Degenerate case)}
  $$ \left\{
     \begin{array}{l}
  \displaystyle{\gamma(r)=0 ~{\rm for~}r \le 0,~~\gamma~{\rm is~strictly
      ~increasing~on~}[0,1)} \\[0.3cm]
    \displaystyle{ \gamma(r) \uparrow \infty~{\rm as ~}r\uparrow 1,
   ~~\gamma(r)=\infty  ~{\rm for~}r \geq 1,}\\[0.3cm]
  \displaystyle{\gamma(r)~{\rm is~locally~Lipschitz~continuous~in~a 
   ~neighborhood~ of~}r=0. }
     \end{array} \right. \eqno{(5.3)}$$

For the enthalpy formulation of the Stefan problem we introduce some function
spaces. Let $V:=H^1(\Omega)$ with norm
   $$ |z|_V := \left\{\int_\Omega |\nabla z|^2dx + n_0 \int_\Gamma |z|^2d\Gamma
     \right \}^{\frac 12},~~\forall z \in V,$$
where $n_0$ is a fixed positive number and $d\Gamma$ is the usual surface 
measure on $\Gamma$.
The dual space $V^*$ is equipped with the dual norm of $V$. In this case
 $$V\subset L^2(\Omega) \subset V^*~~{\rm with~ dense~ and~ compact~ 
 embeddings}$$
and the duality mapping $F: V \to V^*$ is given by
  $$ \langle Fu, z \rangle :=\int_\Omega \nabla u\cdot \nabla z dx
     + n_0\int_\Gamma u z d\Gamma,~~
      \forall u,~z \in V,\eqno{(5.4)}$$
where $\langle \cdot,\cdot \rangle=\langle \cdot,\cdot \rangle_{V^*,V}$.
We know that $F$ is linear, continuous and uniformly monotone from $V$ onto 
$V^*$. Now we set up an inner product $(\cdot,\cdot)_*$ in $V^*$ by
 $$ (u^*,z^*)_* := \langle u^*, F^{-1}z^* \rangle,~~\forall u^*,~z^* \in V^*;
      $$
note that $V^*$ is a Hilbert space with this inner product $(\cdot,\cdot)_*$.
Moreover, we define a proper, l.s.c. and convex function 
$\varphi(\cdot)$ on $V^*$ by
    $$ \varphi(z):= \left\{
        \begin{array}{ll}
  \displaystyle{\int_\Omega \hat \beta(z(x))dx,
            ~~}&\displaystyle{{\rm for~} z \in L^2(\Omega),}\\[0.3cm]
  \displaystyle{\infty, ~~}&\displaystyle{{\rm for~}z \in V^*-L^2(\Omega).}
        \end{array} \right. $$
We know (cf. [11, 12]) that the subdifferential $\partial_*\varphi(\cdot)$ of 
$\varphi(\cdot)$ in the Hilbert space $V^*$ is a singlevalued mapping in 
$V^*$ such that
 $$ D(\partial_*\varphi):=\{z \in L^2(\Omega)~|~\beta(z) \in V\}~{\rm and~}
    \partial_*\varphi(z)=F \beta(z),
     ~~\forall z \in D(\partial_*\varphi). \eqno{(5.5)}$$

The following evolution equation is considered in the space $V^*$ as the 
enthalpy formulation of the Stefan problem:
  $$ \begin{array}{l}
    w'(t)+\partial_*\varphi(w(t))
     +{\rm div}\hspace{0.05cm}(w(t){\bfv}(t)) =h(t)~~{\rm in}~V^*,
     ~~t \in (0,T), \\[0.3cm]
     w(0)=w_0, 
     \end{array}\eqno{(5.6)} $$
where ${\bfv}:={\bfv}(x,t)$ is a vector field on $Q$ in 
$L^2(0,T;{\bfV}_\sigma)$ 
as well as $w_0 \in L^2(\Omega)$ and $h\in L^2(Q)$.

On account of the results in [11, 12], the problem (5.6) with
${\bfv}\equiv 0$ admits one and only one solution $w \in W^{1,2}(0,T;V^*)
\cap C_w([0,T];L^2(\Omega))$ such that 
$t \mapsto \int_\Omega 
\hat\beta(w(\cdot,t))dx $ is absolutely continuous on $[0,T]$, where
$C_w([0,T];L^2(\Omega))$ is the space of all
weakly continuous functions from $[0,T]$ into $L^2(\Omega)$. 
We notice from (5.4) and (5.5) that (5.6) is 
written in the form
 $$ w'(t) +F(\beta(w(t))+{\rm div}\hspace{0.05cm}(w(t){\bfv}(t)) =h(t)~~
   {\rm in~}V^*, ~{\rm a.e.~}t \in (0,T),~w(0)=w_0, \eqno{(5.7)}$$
or 
 $$ w_t  -\Delta \beta(w)+{\bfv} \cdot \nabla w =h~~{\rm in~}Q,~
    \frac{\partial{\beta(w)}}{\partial n}+n_0\beta(w)=0~{\rm on ~}\Sigma,
    ~w(\cdot,0)=w_0~{\rm on~}\Omega,
       \eqno{(5.8)}$$
 $$({\rm in~the~distribution~sense}); $$
note here that ${\bfv}\cdot \nabla w={\rm div}\hspace{0.05cm}(w{\bfv})$ is 
the convective term for the enthalpy
$w$ caused by the fluid flow. \vspace{0.5cm}

We are now consider a coupled system of (5.7) (hence (5.8)) and the 
variational inequality
of Navier-Stokes type discussed in sections 2$\thicksim$4:
 $$ {\bfv} \in \boldsymbol{\cal K}^i(\gamma(w^{\varepsilon_0})),~~{\bfv}(0)=
   {\bfv}_0, $$
$$\int_0^t(\boldsymbol {\mathcal \xi}'(\tau), {\bfv}(\tau)-\boldsymbol {\mathcal \xi} (\tau))_\sigma d\tau
    +\nu \int_0^t \langle {\bfF}{\bfv}(\tau), {\bfv}(\tau)- 
    \boldsymbol {\mathcal \xi}(\tau)\rangle_\sigma d\tau ~~~~~~~~~$$
 $$  ~~+ \int_0^t \langle {\bfG}({\bfv}, {\bfv}),
    {\bfv}-\boldsymbol {\mathcal \xi}\rangle_\sigma d\tau 
     +\frac 12 |{\bfv}(t)-\boldsymbol {\mathcal \xi}(t)|^2_{0,2} 
   \eqno{(5.9)}$$
  $$\le \int_0^t \langle{\bfg},{\bfv}-\boldsymbol {\mathcal \xi}
    \rangle_\sigma d\tau
     + \frac 12 |{\bfv}_0-\boldsymbol {\mathcal \xi}(0)|^2_{0,2},
    ~~\forall t \in [0,T],~\forall \boldsymbol {\mathcal \xi} \in
    \boldsymbol{\cal K}_0^i(\gamma(w^{\varepsilon_0})). $$
In the above formulation, for $i=1,2$ the constraint sets 
$K^i(\gamma(w^{\varepsilon_0}); t)$, the classes 
$\boldsymbol {\cal K}^i(\gamma(w^{\varepsilon_0}))$ of test functions
and $\boldsymbol{\cal K}^i_0(\gamma(w^{\varepsilon_0}))$ of smooth test 
functions are defined by as follows: 
   $$w^{\varepsilon_0}(x,t) =
\int_\Omega\rho_{\varepsilon_0}(x-y) w(y,t)dy,~x \in \Omega, $$ 
 $$K^1(\gamma(w^{\varepsilon_0});t):=\{{\bfz} \in {\bfV}_\sigma~|~
  |{\bfz}| \le \gamma(w^{\varepsilon_0}(\cdot,t))~{\rm a.e.~on~}\Omega\},~
   t \in [0,T], \eqno{(5.10)}$$
 $$ \boldsymbol{\cal K}^1(\gamma(w^{\varepsilon_0})):=
    \{\boldsymbol {\mathcal \xi}\in
  L^2(0,T; {\bfV}_\sigma)~|~\boldsymbol {\mathcal \xi}(t)\in K^1
     (\gamma(w^{\varepsilon_0};t)~{\rm a.e.~}t \in (0,T)\},\eqno{(5.11)}$$
$$ \boldsymbol {\cal K}^1_0(\gamma(w^{\varepsilon_0})):=\left\{
     \boldsymbol {\mathcal \xi} \in
   C^1([0,T];{\bfW}_\sigma)~\left|
      \begin{array}{l}
      |\boldsymbol {\mathcal \xi}| \le
      \gamma(w^{\varepsilon_0}) ~{\rm on~}Q,\\[0.3cm]
      {\rm supp}\hspace{0.05cm}(|\boldsymbol {\mathcal \xi}|)
     \subset \hat Q(\gamma(w^{\varepsilon_0}>0))
      \end{array}
      \right. \right\} \eqno{(5.12)}$$
and
 $$K^2(\gamma(w^{\varepsilon_0});t):=\{{\bfz} \in {\bfV}_\sigma~|~
  |\nabla {\bfz}| \le \gamma(w^{\varepsilon_0}(\cdot,t))~{\rm a.e.~on~}
   \Omega\}, ~ t \in [0,T],\eqno{(5.13)}$$
$$ \boldsymbol {\cal K}^2(\gamma(w^{\varepsilon_0})):=
\{\boldsymbol {\mathcal \xi}\in
  L^2(0,T; {\bfV}_\sigma)~|~\boldsymbol {\mathcal \xi}(t)\in K^2
  (\gamma(w^{\varepsilon_0};t)~{\rm a.e.~}t \in (0,T)\},\eqno{(5.14)}$$
  $$ \boldsymbol {\cal K}^2_0(\gamma(w^{\varepsilon_0})):=\left\{
     \boldsymbol {\mathcal \xi} \in
   C^1([0,T];{\bfW}_\sigma)~\left|
      \begin{array}{l}
      |\nabla \boldsymbol {\mathcal \xi}| \le
      \gamma(w^{\varepsilon_0}) ~{\rm on~}Q,\\[0.3cm]
      {\rm supp}\hspace{0.05cm}(|\nabla \boldsymbol {\mathcal \xi}|)
     \subset \hat Q(\gamma(w^{\varepsilon_0}>0))
      \end{array}
      \right. \right\}. \eqno{(5.15)}$$
The initial datum ${\bfv}_0$ and source ${\bfg}$ are respectively prescribed
in ${\bfV}_\sigma$ and in $L^2(0,T;{\bfH}_\sigma)$.
 \vspace{0.5cm}

\noindent
{\bf Definition 5.1.} 
For given data
 $$ w_0 \in L^2(\Omega),~h \in L^2(Q),~
     ~ {\boldsymbol g}\in L^2(0,T; {\boldsymbol H}_\sigma),~
   {\boldsymbol v}_0 \in {\boldsymbol H}_\sigma(\Omega),
 $$
our problem, referred to $SNS^i(\beta,\gamma; h, {\boldsymbol g},w_0,
{\bfv}_0),~i=1,2$, 
is to find a pair of functions $\{w, {\bfv}\}$
from $[0,T]$ into $V^*\times {\bfH}_\sigma$ satisfying the following 
(i), (ii), (iii) and (iv):
\begin{description}
\item{(i)} $w \in W^{1,2}(0,T;V^*)\cap L^\infty(Q)
\cap C_w([0,T];L^2(\Omega))$, $w(0)=w_0$ and
$t \mapsto \int_\Omega \hat\beta(x,t)dx$ is absolutely continuous on $[0,T]$,
\item{(ii)} ${\bfv} \in L^\infty(0,T;{\bfH}_\sigma)\cap 
{\cal K}^i(\gamma(w^{\varepsilon_0}))$,
${\bfv}(0)={\bfv}_0$ and $t\mapsto ({\bfv}(t), 
\boldsymbol {\mathcal \xi}(t))_\sigma$ is
of bounded variation on $[0,T]$ for every $\boldsymbol {\mathcal \xi}\in
C^1([0,T];{\bfW}_\sigma)$ with 
 $$ 
     {\rm supp}\hspace{0.05cm}
(|\boldsymbol {\mathcal \xi}|) \subset \hat Q(\gamma(w^{\varepsilon_0})>0)~
    {\rm for~}i=1,~~
   {\rm supp}\hspace{0.05cm}
(|\nabla\boldsymbol {\mathcal \xi}|) \subset 
   \hat Q(\gamma(w^{\varepsilon_0})>0)~
    {\rm for~}i=2,$$
\item{(iii)} $w$ satisfies the evolution equation of (5.6) (hence (5.7)).
\item{(iv)} ${\bfv}$ satisfies (5.9).
\end{description}
Such a pair of functions $\{w, {\bfv}\}$ is called a weak solution of 
$SNS^i(\beta, \gamma;h, {\bfg},w_0, {\bfv}_0)$.\vspace{0.5cm}

The existence results are stated in the following two theorems.\vspace{0.5cm}

\noindent
{\bf Theorem 5.1. (Non-degenerate case)} {\it Assume that the data $w_0,~h,
~{\bfv}_0$ and ${\bfg}$ satisfy 
 $$ w_0 \in L^\infty(\Omega),~~h\in L^\infty(Q), \eqno{(5.16)}$$
and
   $$ \begin{array}{l}
   {\bfv}_0 \in {\bfV}_\sigma,~ {\bfg} \in L^2(0,T;{\bfH}_\sigma),\\
 |{\bfv}_0| \in L^\infty(\Omega),~|{\bfv}_0|\le \gamma(w^{\varepsilon_0}_0)
    ~{\it a.e.~on~}\Omega~{\it in~ the~ case~of~}i=1, \\
  |\nabla{\bfv}_0| \in L^\infty(\Omega),~
   |\nabla {\bfv}_0|\le \gamma(w^{\varepsilon_0}_0)
    ~{\it a.e.~on~}\Omega~{\it in~ the~ case~of~}i=2,
      \end{array}  \eqno{(5.17)}$$
where $w_0^{\varepsilon_0}(x)= \int_\Omega \rho_{\varepsilon_0}(x-y)w_0(y)dy,
~x \in \Omega.$ Moreover, assume that $\beta$ and $\gamma$ satisfy (5.1) and
 (5.2), respectively. 
Then $SNS^i(\beta, \gamma; h, {\bfg}, w_0, {\bfv}_0),~i=1,2,$ admits at least 
one weak solution $\{w,{\bfv}\}$.} \vspace{0.5cm}

\noindent
{\bf Remark 5.1.} In the non-degenerate case of $\gamma$ the classes
$\boldsymbol {\cal K}^i_0(\gamma(w^{\varepsilon_0})),~i=1,2,$ of admissible test functions 
are given by
 $$ \boldsymbol {\cal K}^1_0(\gamma(w^{\varepsilon_0})):=\{
     \boldsymbol {\mathcal \xi} \in
   C^1([0,T];{\bfW}_\sigma)~|~
      |\boldsymbol {\mathcal \xi}| \le
      \gamma(w^{\varepsilon_0}) ~{\rm on~}Q\},$$
and
  $$ \boldsymbol {\cal K}^2_0(\gamma(w^{\varepsilon_0})):=\{
     \boldsymbol {\mathcal \xi} \in
   C^1([0,T];{\bfW}_\sigma)~|~
      |\nabla \boldsymbol {\mathcal \xi}| \le
      \gamma(w^{\varepsilon_0}) ~{\rm on~}Q\}.$$
Therefore, just as the case of Theorem 3.1, the mathematical treatment of
non-degenerate case is much easier than the degenerate one.\vspace{0.5cm}

\noindent
{\bf Theorem 5.2. (Degenerate case)} {\it Assume that (5.1) and (5.3) are
satisfied. For the data $w_0,~h,~{\bfv}_0$ and
${\bfg}$, in addition to (5.16) and (5.17) for $i=1$, suppose that
  $$ {\rm supp}\hspace{0.05cm}(|{\bfv}_0|) \subset 
  \{x \in \Omega~|~\gamma(w_0^{\varepsilon_0}(x))>0\}.\eqno{(5.18)}$$
Then $SNS^1(\beta,\gamma;h, {\bfg}, w_0, {\bfv}_0)$ admits at least one weak solution
$\{w,{\bfv}\}$ in the sense of Definition 5.1 for $i=1$.  
}\vspace{0.5cm}

The solvability of $SNS^2(\beta,\gamma;h, {\bfg}, w_0, {\bfv}_0)$ under 
gradient constraint is an open question
in the degenerate case in which some difficulties arise just as in the case of
$NS^2(\gamma;{\bfg},{\bfv}_0)$. \vspace{0.5cm}

\noindent
{\bf 5.2. Approximate problems}\vspace{0.3cm}

We approximate the functions $\beta$ and $\gamma$ by $\beta_\delta$ and
$\gamma_{\delta,N}$ for each small $\delta>0$ and large $N>0$,
which are defined by
  $$ \beta_\delta(r):=\beta(r)+\delta r,~~\forall r \in {\bf R}$$
and 
  $$ \gamma_{\delta,N}(r):=(\gamma(r)\lor \delta)\land N, ~~\forall r
     \in {\bf R}.$$
The proper, l.s.c. and convex function $\varphi_\delta(\cdot)$ 
is defined by
 $$ \varphi_{\delta}(z):= \left\{
        \begin{array}{ll}
  \displaystyle{\int_\Omega \hat \beta_\delta (z(x))dx,
            ~~}&\displaystyle{{\rm for~} z \in L^2(\Omega),}\\[0.3cm]
  \displaystyle{\infty, ~~}&\displaystyle{{\rm for~}z \in V^*-L^2(\Omega),}
        \end{array} \right. $$
where $\hat \beta_\delta(r):=\int_0^r \beta_\delta(s)ds$ for $r \in {\bf R}$.
It is easy to see that $\varphi_{\delta}$ converges to $\varphi$ in the 
sense of Mosco [28] as $\delta \downarrow 0$. Also, 
the evolution equation
$$ w'(t)+\partial_*\varphi_\delta(w(t))+{\rm div}\hspace{0.05cm}(w(t){\bfv}(t))
     =h(t)~~{\it in~}V^*,~{\rm a.e.~}t \in(0,T),~~w(0)=w_0, 
  \eqno{(5.19)}$$
is formulated similarly in the case (5.6). This is equivalently written in the 
form:
 $$ w'(t)+F\beta_\delta(w(t))+{\rm div}\hspace{0.05cm}(w(t){\bfv}(t))
     =h(t)~~{\it in~}V^*,~{\rm a.e.~}t \in(0,T),~~w(0)=w_0.
\eqno{(5.20)}$$

\vspace{0.5cm}

\noindent
{\bf Lemma 5.1.} {\it Let $\{w_n\}$ and $\{{\bfv}_n\}$ be bounded 
sequences in $L^\infty(Q)$ and $L^2(0,T;{\bfV}_\sigma)$,
respectively, and let $w \in L^\infty(Q)$ and
${\bfv} \in L^2(0,T;{\bfV}_\sigma)$ such that 
  $$w_n \to w~{\it weakly~ in~} L^2(Q),~~
{\bfv}_n \to {\bfv} ~{\it in~} L^2(0,T;{\bfH}_\sigma).$$
Then 
  $${\rm div}\hspace{0.05cm}(w_n{\bfv}_n)
   \to {\rm div}\hspace{0.05cm}(w{\bfv})~{\it~ weakly~in~} L^2(0,T;V^*).
    \eqno{(5.21)}$$
Moreover, if $w_n\to w$ in $L^2(Q)$, then the convergence of (5.21) holds 
in $L^2(0,T;V^*)$.}\vspace{0.3cm}

\noindent
{\bf Proof.} First we note by the boundedness of $\{w_n\}$ in $L^\infty(Q)$
that $w_n \to w$ weakly$^*$ in $L^\infty(Q)$. 
For any function $z^* \in L^2(0,T;V^*)$ we have 

 \begin{eqnarray*}
  && \left|\int_0^T({\rm div}\hspace{0.05cm}(w_n{\bfv}_n-w{\bfv}), z^*)_* dt
       \right| \\
  &=& \left|\int_Q {\rm div}\hspace{0.05cm}(w_n{\bfv}_n-w{\bfv})
    F^{-1} z^* dx dt \right|\\
  &=& \left|\int_Q ((w_n-w){\bfv}+w_n({\bfv}_n-{\bfv}))\cdot
    \nabla F^{-1} z^* dxdt\right|\\
  &\le& \left| \int_Q (w_n-w){\bfv}\cdot \nabla F^{-1} z^* dxdt\right|
      + \int_Q |w_n||{\bfv}_n-{\bfv}| |\nabla F^{-1} z^*| dxdt \\
  &\le& \left| \int_Q (w_n-w){\bfv}\cdot \nabla F^{-1} z^* dxdt\right|
      + C |{\bfv}_n-{\bfv}|_{L^2(0,T;{\bfH}_\sigma)}| z^*|_{L^2(0,T;V^*)},\\
  &\to& 0~~{\it as ~}n \to \infty,
 \end{eqnarray*}
where $C=\sup_{n\in {\bf N}}|w_n|_{L^\infty(Q)}$. Hence (5.21) is obtained. 
In particular, if $w_n \to w$ in $L^2(Q)$, then the last 
convergence holds in the strong topology of
$L^2(0,T;V^*)$.
 \hfill $\diamondsuit$ \vspace{0.5cm}

\noindent
{\bf Lemma 5.2.} {\it Assume that (5.18) holds and let ${\bfv}$ be any function
in $L^2(0,T;{\bfV}_\sigma)$. Then we have:
\begin{description}
\item{(a)} For each $\delta>0$, problem (5.6) (hence (5.7))
has a unique solution $w$ in $W^{1,2}(0,T;V^*) \cap 
L^2(0,T;V)\cap C([0,T];L^2(\Omega))$ such that
$\sqrt{t}w' \in L^2(0,T; L^2(\Omega))$ and
$\sqrt{t} |\nabla w|_{L^2(\Omega)} \in L^\infty(0,T)$.
\item{(b)} For each $\delta>0$, the solution $w$ of (5.6) 
satisfies that
 $$ \sup_{t \in [0,T]}|w(t)|^2_{L^2(\Omega)}
      + \delta|w|^2_{L^2(0,T;V)} \le
    |w_0|^2_{L^2(\Omega)}+\frac{C_2^2}\delta |h|^2_{L^2(Q)},\eqno{(5.22)}$$
 $$ \sup_{t \in [0,T]}\int_\Omega \hat\beta_\delta(w(t))dx 
     + \frac 12 |\beta_\delta (w)|^2_{L^2(0,T;V)} \le 
     \int_\Omega \hat\beta(w_0)dx+|w_0|^2_{L^2(\Omega)} 
    +\frac {C_2^2}2 |h|^2_{L^2(Q)},
  \eqno{(5.23)}$$ 
$$ |w|_{L^\infty(Q)} \le R_0,~~\forall \delta \in (0,1], \eqno{(5.24)}$$
where $C_2$ is a positive constant satisfying 
$|z|_{L^2(\Omega)} \le C_2|z|_{V}$ for all $z \in V$
and $R_0$ is a positive constant independent of $\delta \in (0,1]$ 
and ${\bfv}$.
\end{description}} 

\noindent
{\bf Proof.} There are several approaches to the solvability of the approximate problem (5.7). Here we are going to use the $L$-pseudomonotone theory (cf.
[8, 9]) for it.
 
First we consider the case of ${\bfv} \in L^\infty(Q)^3$. In this case we 
observe that 
the operator $Aw:=-\Delta (\beta_\delta(w))+{\rm div}\hspace{0.05cm}
(w{\bfv})$ is coercive and $L$-pseudomotone from $L^2(0,T;V)$
into $L^2(0,T;V^*)$ in case $L=L_{w_0}$ is the time-derivative $\frac {d}{dt}$
associated for constraint $K(t)=V$ (see section 2). Therefore, applying a 
result in [8, 9], the range of $L_{w_0}+A$ is the
whole of $L^2(0,T;V^*)$; namely, given $h \in L^2(Q)$, there is 
$w \in L^2(0,T;V)$ 
such that $L_{w_0}w+Aw=h$ in $L^2(0,T;V^*)$. Moreover, since 
$\beta_\delta(\cdot)$
is bi-Lipschitz continuous on ${\bf R}$, it follows from the general theory 
in [24; Chapter 2] that the solution $w$ of (5.20) has the required regularity 
properties, and 
as to the uniqueness of solutions, for two solutions $w_k,~k=1,2,$
of (5.20) associated for the initial datum $w_{k0}\in L^2(\Omega)$ and
the source $h_k \in L^2(Q)$ we have for all $t \in [0,T]$
  $$ |(w_1(t)-w_2(t))^+|_{L^1(\Omega)} \le |(w_{10}-w_{20})^+|_{L^1(\Omega)}
      +\int_0^t |(h_1(\tau)-h_2(\tau))^+|_{L^1(\Omega)}d \tau,
      $$
where $(\cdot)^+$ denotes the positive part of $(\cdot)$. In particular, if
 $w_{10}=w_{20}$ and $h_1=h_2$, then this shows $w_1=w_2$ on $Q$, namely the
uniqueness of the solution of (5.20). Thus we have (a).

Next we show (b).
Multiplying (5.20) by $w$ and noting that 
   \begin{eqnarray*}
   \int_\Omega {\rm dvi}\hspace{0.05cm}(w(x,t) {\bfv}(x,t))w(x,t)dx &=&
    -\int_\Omega w(x,t){\bfv}(x,t)\cdot \nabla w(x,t) dx \\
    &=&-\int_\Omega {\bfv}(x,t)\cdot \nabla\left(\frac 12 w(x,t)^2\right)dx\\
    &=&\frac 12\int_\Omega {\rm div}\hspace{0.05cm}{\bfv}(x,t)w(x,t)^2dx=0
   \end{eqnarray*}
we have that
  $$\frac 12 \frac d{dt} |w(t)|^2_{L^2(\Omega)} +\delta \int_\Omega
   \{\beta'(w(x,t))|\nabla w(x,t)|^2 +\delta|\nabla w(x,t)|^2\} dx $$ 
  $$ + n_0\int_\Gamma \{\beta(w(x,t))w(x,t)+\delta|w(x,t)|^2\}d\Gamma
    \le \int_\Omega|h(x,t)||w(x,t)|dx  $$
for all $t \in [0,T]$. By integrating it over $[0,T]$ 
in time, we have an estimate of the form (5.22).
The estimate (5.23) is similarly obtained. In fact, multiply (5.20)
by $\beta_\delta(w)$. Then
  $$\frac d{dt} \int_\Omega \hat \beta_\delta(x,t)dx + \int_\Omega 
  |\nabla \beta_\delta(w(x,t))|^2dx + n_0\int_\Gamma |\beta_\delta(w)|^2d\Gamma
  \le \int_\Omega h \beta_\delta(w(x,t))dx.$$
Integrating this over $[0,t]$ in time, we easily get an estimate of the form
(5.23).

Next, let $R_1$ be a positive constant such that
$|w_0|_{L^\infty(\Omega)} \le R_1$ and $|h|_{L^\infty(Q)} \le R_1$. Note that
the solution $w$ of (5.20) satisfies
  $$w_t -\Delta \beta_\delta(w)+{\bfv}\cdot \nabla w
    = h~~{\rm a. e.~in~}Q, ~\frac {\partial \beta_\delta(w)}{\partial n}+n_0
      \beta_\delta(w)=0~{\rm a.e.~on~}\Sigma. $$ 
Multiplying this  by $(w -R_1-R_1t)^+$, we obtain that
  $$\int_\Omega(w_t-R_1)(w-R_1-R_1t)^+ dx+\int_\Omega 
   \beta'_\delta(w)|\nabla (w-R_1-R_1t)^+|^2 dx
   ~~~~~~~~~~~~~~~$$
   $$+\int_\Omega \frac 12 {\bfv}\cdot \nabla [(w-R_1-R_1t)^+]^2 dx
   = \int_\Omega (h-R_1)(w-R_1-R_1t)^+ dx \le 0.$$
In the above inequality the second integral is non-negative  and third integral is zero becaose of the divergencefreeness of ${\bfv}$, so that
  $$ \frac d{dt} \int_\Omega |(w(t)-R_1-R_1t)^+|^2dx \le 0~~{\rm a.e.~}
    t \in [0,T].$$
Hence
 $$\int_\Omega |(w(t)-R_1-R_1t)^+|^2dx \le \int_\Omega |(w_0-R_1)^+|^2dx
    =0,~\forall t \in [0,T],$$
namely $w \le R_1+R_1T=:R_0$ a.e. on $Q$ and (5.24) is obtained.

In the general case of ${\bfv} \in L^2(0,T; {\bfV}_\sigma)$, we approximate 
${\bfv}$ by a bounded solenoidal function ${\bfv}_\varepsilon, ~\varepsilon>0,$
so that
$|{\bfv}-{\bfv}_\varepsilon|_{L^2(0,T;{\bfV}_\sigma)} \le \varepsilon$. For
each $\varepsilon>0$ we denote by $w_\varepsilon$ the solution of (5.18) 
associated for convection vector ${\bfv}_\varepsilon$. Then, by the above
results, estimates (5.22), (5.23) and (5.24) hold true and they are independent
of $\varepsilon$. By using these uniform estimates and Lemma 5.1, it is not 
difficult to show
that $w_\varepsilon$ converges to the solution $w$ of (5.19)  
as $\varepsilon \to 0$, and 
we have the conclusion of the lemma. 
    \hfill $\diamondsuit$ \vspace{0.5cm}

According to Lemma 5.2, we see that the solution $w$ of (5.6) (hence (5.7)) is 
constructed in the class
  $$ {\cal Y}:= \left\{w \in W^{1,2}(0,T;V^*)~ \left |~
     \begin{array}{l}
       w(0)=w_0,~~\beta_\delta(w) \in L^2(0,T;V),\\[0.2cm]
     (5.22),~(5.23)~ {\rm and~} (5.24)~{\rm hold}~{\rm for~all~}
       \delta \in (0,1]
     \end{array} \right. \right\}.
  $$
With the functions $\beta_\delta$ and $\gamma_{\delta,N}$
approximate problem $SNS^i(\beta, \gamma_{\delta,N};h, {\bfg},w_0,{\bfv}_0)$
is formulated as follows.\vspace{0.5cm}

We use the notation $K^i(\gamma_{\delta,N};t),~\boldsymbol {\cal K}^i
(\gamma_{\delta,N}
(w^{\varepsilon_0}))$ and $\boldsymbol {\cal K}^i_0(\gamma_{\delta,N}
(w^{\varepsilon_0}))$,
etc., which are defined by (5.12)-(5.17) with $\gamma$ replaced by 
$\gamma_{\delta,N}$. \vspace{0.5cm}

\noindent
{\bf Definition 5.2.} We denote by $SNS^i(\beta_\delta,\gamma_{\delta,N}; 
h, {\bfg}, w_0, {\bfv}_0)$ the problem to find a pair of
functions $\{w,{\bfv}\}$ which satisfies:
\begin{description}
\item{(i)} $w \in W^{1,2}(0,T;V^*)\cap L^2(0,T;V) 
\cap C([0,T];L^2(\Omega)),~w(0)=w_0$ and
  $$w'(t)+F(\beta_\delta(w(t)))+{\rm div}\hspace{0.05cm}(w(t){\bfv}(t))
      =h(t)~~{\rm in ~}V^*,~{\rm a.e.~}t \in (0,T), $$
\item{(ii)} ${\bfv} \in C([0,T];{\bfH}_\sigma)\cap 
\boldsymbol {\cal K}^i(\gamma_{\delta,N}(w^{\varepsilon_0}))$, 
  ${\bfv}(0)={\bfv}_0$ and
  $$\int_0^t(\boldsymbol {\mathcal \xi}'(\tau), {\bfv}(\tau)-
 \boldsymbol {\mathcal \xi} (\tau))_\sigma d\tau
    +\nu \int_0^t \langle {\bfF}{\bfv}(\tau), {\bfv}(\tau)- 
    \boldsymbol {\mathcal \xi}(\tau)\rangle_\sigma d\tau ~~~~~~~~~$$
 $$  ~~+ \int_0^t \langle {\bfG}({\bfv}, {\bfv}),
    {\bfv}-\boldsymbol {\mathcal \xi}\rangle_\sigma d\tau 
     +\frac 12 |{\bfv}(t)-\boldsymbol {\mathcal \xi}(t)|^2_{0,2} 
   $$
  $$\le \int_0^t \langle{\bfg},{\bfv}-\boldsymbol {\mathcal \xi}
    \rangle_\sigma d\tau
     + \frac 12 |{\bfv}_0-\boldsymbol {\mathcal \xi}(0)|^2_{0,2},
    ~~\forall t \in [0,T],~\forall \boldsymbol {\mathcal \xi} \in
    \boldsymbol {\cal K}_0^i(\gamma_{\delta,N}(w^{\varepsilon_0})), $$
\end{description}

As to the approximate problem $SNS^i(\beta_\delta, \gamma_{\delta,N};
h, {\bfg}, w_0, {\bfv}_0)$
we prove:\vspace{0.5cm}

\noindent
{\bf Proposition 5.1.} {\it Assume (5.16) and (5.17) hold for the data. Then
for each $\delta >0$ and $N>0$ there is at least one solution 
$\{w_{\delta,N}, {\bfv}_{\delta,N}\}$ of $SNS^i(\beta_\delta,\gamma_{\delta,N};
h, {\bfg},w_0, {\bfv}_0)$ in the sense of Definition 5.2 such that
   $$ |w_{\delta,N}|_{L^\infty(Q)} \le R_0 \eqno{(5.25)} $$
and 
   $$\sup_{t \in [0,T]}|{\bfv}_{\delta,N}(t)|^2_{0,2}
    +\nu |{\bfv}_{\delta,N}|^2_{L^2(0,T;{\bfV}_\sigma)} \le |{\bfv}_0|^2_{0,2}
     +\frac {C_0^2}\nu |{\bfg}|^2_{L^2(0,T;{\bfH}_\sigma)}, \eqno{(5.26)}$$
where $R_0$ is the same one as in (b) of Lemma 5.1 and $C_0$ is a positive
constant satisfying $|{\bfz}|_{0,2}\le C_0|{\bfz}|_{1,2}$ for all 
${\bfz} \in {\bfV}_\sigma$. }
\vspace{0.3cm}

\noindent
{\bf Proof.} We shall prove the proposition by the Schauder fixed point
argument.

\noindent 
{\bf (Step 1)} Let $\bar {\bfv}$ be any function in $L^2(0,T;{\bfV}_\sigma)
\cap L^\infty(Q)^3$ and denote by $w$ the unique solution of
 $$ w'(t)+F\beta_\delta(w(t))+{\rm div}\hspace{0.05cm}(w(t)\bar{\bfv}(t))
    =h(t)~{\rm in~}V^*,~{\rm a.e.~}t \in (0,T),~~w(0)=w_0. \eqno{(5.27)}$$
By Lemma 5.1, $w \in W^{1,2}(0,T;V^* )\cap L^2(0,T;V)$ and
$|w|_{L^\infty(Q)} \le R_0.$ Since 
$\gamma_{\delta,N} \geq \delta$ on ${\bf R}$, we have
    $$ \frac \delta{C_1} B_{{\bfW}^*_\sigma}(0) \le 
      K^i(\gamma_{\delta,N}(w^{\varepsilon_0});t),~~\forall t \in [0,T],
       \eqno{(5.28)}$$
where $C_1$ is a positive constant satisfying
$$|{\bfz}|_{C(\overline \Omega)^3}\le C_1|{\bfz}|_{{\bfW}_\sigma}
       ~{\rm in~the~ case~of~}i=1;~~
   |\nabla{\bfz}|_{C(\overline \Omega)^3}\le C_1|{\bfz}|_{{\bfW}_\sigma}
       ~{\rm in~the~ case~of~}i=2.$$
In particular, if $\gamma$ is non-degenerate, i.e. $\gamma\geq c_*$, then
  $$ \frac {c_*}{C_1} B_{{\bfW}^*_\sigma}(0) \le 
   K^i(\gamma_{\delta,N}(w^{\varepsilon_0});t),~~\forall t \in [0,T],~\forall
       \delta \in (0,c_*],~\forall {\rm large~}N>0.
       $$

Now, consider the variational inequality:
 $$ {\bfv} \in C([0,T];{\bfH}_\sigma)\cap  \boldsymbol {\cal K}^i
  (\gamma_{\delta,N}(w^{\varepsilon_0})),~{\bfv}(0)={\bfv}_0,
          $$
$$\int_0^t(\boldsymbol {\mathcal \xi}'(\tau), {\bfv}(\tau)-\boldsymbol {\mathcal \xi} (\tau))_\sigma d\tau
    +\nu \int_0^t \langle {\bfF}{\bfv}(\tau), {\bfv}(\tau)- 
    \boldsymbol {\mathcal \xi}(\tau)\rangle_\sigma d\tau ~~~~~~~~~$$
 $$  ~~+ \int_0^t \langle {\bfG}({\bfv}, {\bfv}),
    {\bfv}-\boldsymbol {\mathcal \xi}\rangle_\sigma d\tau 
     +\frac 12 |{\bfv}(t)-\boldsymbol {\mathcal \xi}(t)|^2_{0,2} 
   \eqno{(5.29)}$$
  $$\le \int_0^t \langle{\bfg},{\bfv}-\boldsymbol {\mathcal \xi}
    \rangle_\sigma d\tau
     + \frac 12 |{\bfv}_0-\boldsymbol {\mathcal \xi}(0)|^2_{0,2},
    ~~\forall t \in [0,T],~\forall \boldsymbol {\mathcal \xi} \in
    \boldsymbol {\cal K}_0^i(\gamma_{\delta,N}(w^{\varepsilon_0})), $$
Then, by (5.28), Proposition 3.1 gives that this variational inequality 
possesses one and only one solution ${\bfv}$ together with estimates:
 $$ \sup_{t \in [0,T]}|{\bfv}(t)|^2_{0,2} \le |{\bfv}_0|^2_{0,2} +
         \frac {C_0}\nu |{\bfg}|^2_{L^2(0,T;{\bfH}_\sigma)}=:N_1,~
    \nu |{\bfv}|^2_{L^2(0,T;{\bfV}_\sigma)} \le N_1. \eqno{(5.30)}$$
Furthermore, as is easily seen, the weak solution ${\bfv}$ of (5.29) belongs 
to the set $Z_2(\frac \delta{C_1}, N'_1, {\bfv}_0)$, 
where $N'_1$ is a positive constant depending only on $N_1$ and 
$\delta>0$. 
Therefore, by Lemma 2.2, 
 $$  {\rm Var}_{{\bfW}^*_\sigma}({\bfv}) \le N'_1+\frac {C_0}\delta N'_1
      +\frac {C_0}{2\delta}|{\bfv}_0|^2_{0,2}=:C^*(\delta). \eqno{(5.31)}$$
In particular, when $\gamma$ is non-degenerate, we have
$$  {\rm Var}_{{\bfW}^*_\sigma}({\bfv}) \le N'_1+\frac {C_0}{c_*} N'_1
      +\frac {C_0}{2c_*}|{\bfv}_0|^2_{0,2}=:C^*. $$

Now, we consider the set ${\cal X}_\delta$ given by:
   $$ {\cal X}_\delta:=\left\{{\bfv}\in L^2(0,T;{\bfV}_\sigma)~\left|~
           \begin{array}{l}
           \nu |{\bfv}|^2_{L^2(0,T;{\bfV}_\sigma)} \le N_1, \\
        |{\bfv}|^2_{L^\infty(0,T;{\bfH}_\sigma)} \le N_1, \\
        {\rm Var}_{{\bfW}^*_\sigma}({\bfv}) \le C^*(\delta)
           \end{array} \right. \right\}. \eqno{(5.32)}$$
By Lemma 2.1, ${\cal X}_\delta$ is a compact convex subset of 
$L^2(0,T;{\bfH}_\sigma)$. Also, in the non-degenerate case of $\gamma$,
we consider 
  $$ {\cal X}_{c_*}:=\left\{{\bfv}\in L^2(0,T;{\bfV}_\sigma)~\left|~
           \begin{array}{l}
           \nu |{\bfv}|^2_{L^2(0,T;{\bfV}_\sigma)} \le N_1, \\
        |{\bfv}|^2_{L^\infty(0,T;{\bfH}_\sigma)} \le N_1, \\
        {\rm Var}_{{\bfW}^*_\sigma}({\bfv}) \le C^*
           \end{array} \right. \right\}. \eqno{(5.33)}$$

We denote by ${\cal S}_1$ the mapping which assigns to $\bar{\bfv} \in 
{\cal X}_\delta$ the 
solution $w$ of (5.27), namely $w={\cal S}_1\bar {\bfv}$, and by ${\cal S}_2$
the mapping which assigns to $w$ the weak solution ${\bfv}$ of (5.29), 
${\bfv}={\cal S}_2w$.
The composition of ${\cal S}={\cal S}_2{\cal S}_1$ of ${\cal S}_1$ and 
${\cal S}_2$ is a mapping which maps ${\cal X}_\delta$ into itself, and
 ${\bfv}={\cal S}\bar{\bfv}$.
\vspace{0.3cm}

\noindent
{\bf (Step 2)} We now show the continuity of ${\cal S}$ in ${\cal X}_\delta$
with respect to the topology of $L^2(0,T;{\bfH}_\sigma)$. 
Let $\{\bar{\bfv}_n\}$ be a 
sequence in ${\cal X}_\delta$ such that $\bar{\bfv}_n \to \bar{\bfv}$ in
$L^2(0,T;{\bfH}_\sigma)$ and put $w_n:={\cal S}_1\bar {\bfv}_n$, which
is the solution of
 $$w'_n +F\beta_\delta(w_n) +{\rm div}\hspace{0.05cm}(w_n\bar{\bfv}_n)=h~
  {\rm in~}V^*,~{\rm a.e.~on~}(0,T),~w_n(0)=w_0. \eqno{(5.34)}$$
Note  (cf. [24; Chapter 2]) that $\{w_n\}$ is bounded in $W^{1,2}(0,T;V^*)\cap 
C([0,T];L^2(\Omega))\cap L^2(0,T;V)$ with $|w_n|_{L^\infty(Q)}\le R_0$, 
whence it follows from the Aubin compactness theorem [3, 27] that $\{w_n\}$
is relatively compact in $L^2(Q)$. Now, choose any subsequence $\{w_{n_k}\}$
of $\{w_n\}$ so that $w_{n_k} \to w$ in $L^2(Q)$; for this subsequence we have
that $w_{n_k} \to w$ weakly in $W^{1,2}(0,T;V^*)\cap L^2(0,T;V)$ 
(as $k \to \infty$) and 
$\nabla \beta_\delta(w_n)=\beta'_\delta(w_{n_k})\nabla w_{n_k} \to 
\beta'_\delta(w)\nabla w=\nabla \beta_\delta(w)$ weakly in $L^2(Q)^3$.
Besides, by Lemma 5.2, ${\rm div}\hspace{0.05cm}(w_{n_k}\bar{\bfv}_{n_k}) \to
{\rm div}\hspace{0.05cm}(w\bar{\bfv})$ in $L^2(0,T;V^*)$.
Therefore, letting $k \to \infty$ in (5.34) with $n=n_k$ yields that
 $$w' +F\beta_\delta(w) +{\rm div}\hspace{0.05cm}(w\bar{\bfv})=h~
  {\rm in~}V^*,~{\rm a.e.~on~}(0,T),~w(0)=w_0. \eqno{(5.35)}$$
Since the solution of (5.35) is unique, it follows that the above argument
holds true without extracting any subsequence from $\{w_n\}$. Namely we have
shown that $w_n={\cal S}_1\bar{\bfv}_n$ converges to $w={\cal S}_1\bar{\bfv}$
in weakly in $W^{1,2}(0,T;V^*)\cap L^2(0,T;V)$ and in $L^2(Q)$.

Next, put ${\bfv}_n:={\cal S}_2w_n$ which is the weak solution of
$$ {\bfv}_n \in C([0,T];{\bfH}_\sigma)\cap \boldsymbol {\cal K}^i
(\gamma_{\delta,N}
(w_n^{\varepsilon_0})),~{\bfv}_n(0)={\bfv}_0,$$
$$\int_0^t(\boldsymbol {\mathcal \xi}'(\tau), {\bfv}_n(\tau)
   -\boldsymbol {\mathcal \xi} (\tau))_\sigma d\tau
    +\nu \int_0^t \langle {\bfF}{\bfv}_n(\tau), {\bfv}_n(\tau)- 
    \boldsymbol {\mathcal \xi}(\tau)\rangle_\sigma d\tau ~~~~~~~~~$$
 $$  ~~+ \int_0^t \langle {\bfG}({\bfv}_n, {\bfv}_n),
    {\bfv}_n-\boldsymbol {\mathcal \xi}\rangle_\sigma d\tau 
     +\frac 12 |{\bfv}_n(t)-\boldsymbol {\mathcal \xi}(t)|^2_{0,2} 
   \eqno{(5.36)}$$
  $$\le \int_0^t \langle{\bfg},{\bfv}_n-\boldsymbol {\mathcal \xi}
    \rangle_\sigma d\tau
     + \frac 12 |{\bfv}_0-\boldsymbol {\mathcal \xi}(0)|^2_{0,2},
    ~~\forall t \in [0,T],~\forall \boldsymbol {\mathcal \xi} \in
    \boldsymbol {\cal K}_0^i(\gamma_{\delta,N}(w_n^{\varepsilon_0})). $$
Since $w_n \to w$ in weakly in $W^{1,2}(0,T;V^*)\cap L^2(0,T;V)$, 
strongly in $L^2(Q)$ and
$\{w_n\}$ is uniformly bounded on $\overline Q$, it follows that
$w_n^{\varepsilon_0} \to w^{\varepsilon_0}$ in 
$C(\overline Q)$ as well as $\gamma_{\delta,N}(w_n^{\varepsilon_0}) \to
\gamma_{\delta,N}(w^{\varepsilon_0})$
in $C(\overline Q)$. This convergence shows that 
$K^i(\gamma_{\delta,N}(w^{\varepsilon_0}_n);t) \Longrightarrow
 K^i(\gamma_{\delta,N}(w^{\varepsilon_0});t)$ on $[0,T]$ (see section 2).
Therefore, by virtue of the convergence result in [17; Theorem 2.2], 
${\bfv}_n$ converges to the weak solution ${\bfv} ~(\in {\cal X}_\delta)$ of 
(5.36) in $C([0,T];{\bfH}_\sigma)$, satisfying (5.30) and (5.31).
Moreover,  by Lemma 5.2, 
${\rm div}\hspace{0.05cm}(w_n{\bfv}_n)\to {\rm div}\hspace{0.05cm}(w{\bfv})$
in $L^2(0,T;V^*)$, so that it follows from the general convergence result
on evolution equations generated by subdifferentials (cf. [24; Chapter 1])
that $w$ is the solution of 
  $$w'+F\beta_\delta(w)+{\rm div}\hspace*{0.05cm}(w{\bfv})=h~{\rm in~}V^*,
    ~{\rm a.e.~on~}(0,T),~w(0)=w_0.$$
Hence ${\cal S}_2w_n={\bfv}_n \to {\bfv}={\cal S}_2w$ in 
$L^2(0,T;{\bfH}_\sigma)$, which shows that ${\cal S}\bar{\bfv}_n \to {\cal S}
\bar{\bfv}$ in $L^2(0,T; {\bfH}_\sigma)$. 
Thus ${\cal S}$ is continuous in 
${\cal X}_\delta$ with respect to the topology of $L^2(Q)$.

Now, applying the fixed-point theorem
to the continuous mapping ${\cal S}: {\cal X}_\delta \to {\cal X}_\delta$
we can find a fixed point ${\bfv}_{\delta,N}$ of ${\cal S}$, i.e. 
${\bfv}_{\delta,N}={\cal S}{\bfv}_{\delta,N}$, which gives a weak 
solution $\{w_{\delta,N},{\bfv}_{\delta,N}\}$ with $w_{\delta,N}=
{\cal S}_1{\bfv}_{\delta,N}$.
\hfill $\diamondsuit$ \vspace{1cm}

\noindent
{\bf 5.3. Proof of Theorems 5.1 and 5.2}\vspace{0.3cm}

\noindent
{\bf (Non-degenerate case)} \vspace{0.3cm}

In the non-degenerate case we note that
    $$\gamma_{\delta,N}(r)=\gamma(r)\land N=:\gamma_N(r),~~\forall r \in 
     {\bf R},~\forall \delta \in (0,c_*],~\forall {\rm large~} N.$$
By virtue of Proposition 5.1 and Lemma 5.2, for every $0<\delta \le c_*$ and 
large $N$ the approximate problem $SNS^i(\beta_\delta,\gamma_{N}; h,
{\bfg},w_0,{\bfv}_0)$ admits a pair of solutions 
$\{w_{\delta,N},{\bfv}_{\delta,N}\}$ in the class ${\cal Y}\times 
{\cal X}_{c_*} $. Hence there are a sequence $\{\delta_n,N_n\}$
with $\delta_n \downarrow 0$ and $N_n \uparrow \infty$ (as $n \to \infty$) 
and a pair of function $\{w,{\bfv}\}$ with $\tilde \beta \in L^2(0,T;V)$ such 
that
  $$ \left \{\begin{array}{l} 
     w_n:=w_{\delta_n,N_n} \to w ~{\rm weakly}^*~{\rm in~}L^\infty(Q),
      ~{\rm weakly~in~}W^{1,2}(0,T;V^*),\\[0.3cm]
   \beta_{\delta_n}(w_n) \to \tilde\beta~{\rm weakly~in~}L^2(0,T;V)
     \end{array} \right. \eqno{(5.37)} $$
and
  $$ \left \{
     \begin{array}{l} 
     {\bfv}_n:={\bfv}_{\delta_n,N_n} \to {\bfv}~{\rm weakly~in~}
     L^2(0,T;{\bfV}_\sigma),\\[0.3cm]
     ~~~~~~~~~~~~~~~~~~~~~~~~{\rm weakly}^*~{\rm in~}
     L^\infty(0,T;{\bfH}_\sigma) ~{\rm and ~in~}L^2(0,T;{\bfH}_\sigma).
     \end{array} \right. \eqno{(5.38)}$$
Moreover, by (5.37), $w_n^{\varepsilon_0}(t) \in C^\infty(\overline\Omega)$
and $\{w_n^{\varepsilon_0}\}$ is bounded in 
$W^{1,2}(0,T;C(\overline\Omega))$. 
This shows that $w_n^{\varepsilon_0} \to w^{\varepsilon_0}$ in 
$C(\overline Q)$, so that $\gamma_{N_n}(w_n^{\varepsilon_0}) \to
\gamma(w^{\varepsilon_0})$ in the following sense:
 $$ \left\{ \begin{array}{l}
  \forall \kappa >0, \gamma_{N_n}(w_n^{\varepsilon_0}) \to 
      \gamma(w^{\varepsilon_0})
   ~{\rm uniformly~on~}\overline Q(\gamma(w^{\varepsilon_0})\le k),\\[0.3cm]
  \forall {\rm large~} M>0, ~\exists M'~(>M),~\exists n_M~{\rm such~that~}\\
    ~~~~~~~~~~~~\gamma_{N_n}(w_n^{\varepsilon_0}) > M~{\rm on~}
    \overline Q(\gamma >M'),~ \forall n \geq n_M. 
    \end{array} \right.\eqno{(5.39)}$$

First we prove:\vspace{0.5cm}

\noindent
{\bf Lemma 5.3.} {\it $\tilde \beta=\beta(w)$, and 
$\beta_{\delta_n}(w_n) \to \beta(w)$ in $L^2(Q)$ and
weakly in $L^2(0,T;V)$ as $n \to \infty$.}

\noindent
{\bf Proof.} Taking the inner product of the difference of two equations
$w'_n-w'_m +F(\beta_{\delta_n}(w_n)-\beta_{\delta_m}(w_m))+{\rm div}
    \hspace{0.05cm}(w_n{\bfv}_n-w_m{\bfv}_m)=0$ and $w_n-w_m$ in $V^*$
and integrating the resultant in time, we have
  $$\frac 12 |w_n(t)-w_m(t)|^2_{V^*}
   +\int_0^t \int_\Omega (\beta(w_n)-\beta(w_m))(w_n-w_m)dxd\tau $$
 $$~~~~~~~~~~ +\int_\Omega(\delta_nw_n-\delta_mw_m)(w_n-w_m) dx d\tau$$ 
  $$ = \int_0^t \int_\Omega (w_n{\bfv}_n-w_m{\bfv}_m)\cdot \nabla
     F^{-1}(w_n-w_m)dxd\tau, ~~\forall t \in [0,T],$$
so that
 $$\frac 12 |w_n(t)-w_m(t)|^2_{V^*}
   +\frac 1{L_\beta}\int_0^t \int_\Omega |\beta(w_n)-\beta(w_m)|^2 dxd\tau $$
 $$\le \int_0^t \int_\Omega|\delta_nw_n-\delta_mw_m||w_n-w_m| dx d\tau~~~~~
  ~~~~~~~~~~~~~~~~~~~~\eqno{(5.40)}$$ 
  $$ + \int_0^t \int_\Omega (w_n{\bfv}_n-w_m{\bfv}_m)\cdot \nabla
     F^{-1}(w_n-w_m)dxd\tau,$$
for all $t \in [0,T]$, where $L_\beta$ is a (positive) Lipschitz constant of 
$\beta$. We note here that $w_n(t)-w_m(t) \to 0$ weakly in $L^2(\Omega)$ 
for every $t \in [0,T]$ and $F^{-1}$ is linear and bounded from $L^2(\Omega)$ 
into $H^2(\Omega)$ and compact from $L^2(\Omega)$ into $V$. This implies that
$F^{-1}(w_n(t)-w_m(t)) \to 0$ in $V$ for all $t \in [0,T]$ and the last 
integral in (5.40) converges to $0$ as $n,~m\to \infty$. As a consequence,
$\beta(w_n)$ converges in $L^2(Q)$ as well as $\beta_{\delta_n}(w_n)$ in 
$L^2(Q)$ as $n \to \infty$. Since $w_n \to w$ weakly in $L^2(Q)$ and 
$\beta(\cdot)$ is a maximal monotone mapping in $L^2(Q)$. It follows from the
demiclosedness of maximal monotone mappings that
$\beta(w_n)\to \beta(w)$ in $L^2(Q)$, hence $\beta_{\delta_n}(w_n)=
\beta(w_n)+\delta_n w_n\to \beta(w)$ in $L^2(Q)$. 

Also, by (5.22), $\{\beta_{\delta_n}(w_n)\}$ is bounded in $L^2(0,T;V)$ which
implies that $\beta_{\delta_n}(w_n) \to \beta(w)$ weakly in $L^2(0,T;V)$.
\hfill $\diamondsuit$ \vspace{0.5cm}

\noindent
{\bf Proof of Theorem 5.1:} We denote $SNS^i(\beta_{\delta_n},
\gamma_{\delta_n, N_n}; h, {\bfg},w_0,{\bfv}_0)$ by $SNS^i_{\delta_n,N_n}$.
In the non-degenerate case, the approximate solution $\{w_n,{\bfv}_n\}$
of $SNS^i_{\delta_n,N_n}$ 
is constructed in
${\cal Y}\times {\cal X}_{c_*}$. We already have seen
(5.37), (5.38) with $\tilde \beta=\beta(w)$ and (5.39). 
Since the estimates in the class ${\cal Y}\times {\cal X}_{c_*}$
are independent of parameter $N_n$,
we can show by using convergence property (5.39) of obstacle functions
$\gamma_{N_n}(w_n^{\varepsilon_0})$
(as mentioned in Remark 3.4)
that the convergence of $\{w_n,{\bfv}_n\}$
are obtained just as in the proof of Proposition 5.1 and
$\{w,{\bfv}\} \in {\cal Y}\times {\cal X}_{c^*}$ 
is a weak solution of $SNS^i(\beta,\gamma; h, {\bfg}, w_0,{\bfv}_0)$.
The limit ${\bfv}$ of $\{{\bfv}_n\}$ may lose the continuity in ${\bfH}_\sigma$
in time, but it still stays in the class ${\cal X}_{c_*}$,
whence $t\mapsto ({\bfv}(t),\boldsymbol {\mathcal \xi}(t))_\sigma$ is
of bounded variation on $[0,T]$ for every $\boldsymbol {\mathcal \xi}\in
C^1([0,T];{\bfW}_\sigma)$.
\hfill $\diamondsuit$
\vspace{0.5cm}

\noindent
{\bf Proof of Theorem 5.2:} 
We consider only the problem $SNS^1(\beta,\gamma; h,{\bfg},w_0,{\bfv}_0)$.
In the degenerate case, by Proposition 5.1, for given parameters 
$\delta_n \downarrow 0$ and $N_n \uparrow \infty$ our approximate solution 
$\{w_n,{\bfv}_n\}:=\{w_{\delta_n,N_n},{\bfv}_{\delta_n,N_n}\}$ was constructed
 in ${\cal Y}\times {\cal X}_{\delta_n}$
and we may assume by the uniform estimates (5.22), (5.23), (5.25) and (5.26)
that $\{w_n,{\bfv}_n\}$ satisfies 
   $$ \begin{array}{l} 
  w_n \to w ~{\rm weakly}^*~{\rm in~}L^\infty(Q),
      ~{\rm weakly~in~}W^{1,2}(0,T;V^*),\\[0.2cm]
   \beta_{\delta_n}(w_n) \to \beta(w)~{\rm in~}L^2(Q)~{\rm and~weakly~in~}
    L^2(0,T;V)
   \end{array} \eqno{(5.41)}$$
and
  $$ {\bfv}_n  \to {\bfv}~{\rm weakly~in~}
     L^2(0,T;{\bfV}_\sigma)~{\rm and~weakly}^*~{\rm in~}
     L^\infty(0,T;{\bfH}_\sigma). \eqno{(5.42)}$$
From (5.41) it follows that 
 $$ \left\{ \begin{array}{l}
  \forall \kappa >0, \gamma_{\delta_n, N_n}(w_n^{\varepsilon_0}) \to 
      \gamma(w^{\varepsilon_0})
   ~{\rm uniformly~on~}\overline Q(\gamma(w^{\varepsilon_0})\le k),\\[0.3cm]
  \forall {\rm large~} M>0, \exists M'~(>M),~\exists n_M~{\rm such~that~}\\
    ~~~~~~~~~~~~\gamma_{\delta_n,N_n}(w_n^{\varepsilon_0}) > M~{\rm on~}
    \overline Q(\gamma >M'),~ \forall n \geq n_M. 
    \end{array} \right.\eqno{(5.43)}$$
Now, applying Lemma 4.6 with Corollary 4.3 to $\{{\bfv}_n\}$
and the sequence of obstacle functions 
$\{\gamma_{\delta_n,N_n}(w_n^{\varepsilon_0})\}$ satisfying (5.43), we 
conclude that
  $$ {\bfv}_n \to {\bfv}~~{\rm in~}L^2(0,T;{\bfH}_\sigma). \eqno{(5.44)}$$
Hence, by Lemma 5.1, ${\rm div}\hspace{0.05cm}(w_n{\bfv}_n) \to
{\rm div}\hspace{0.05cm}(w{\bfv})$ weakly in $L^2(0,T;V^*)$ and $w_n$
converges to $w$ weakly in $W^{1,2}(0,T;V^*)$, which is the solution of
 $$ w'+F\beta(w)+{\rm div}\hspace{0.05cm}(w{\bfv})=h~~{\rm in~}V^*,~
    {\rm a.e.~on~}(0,T),~w(0)=w_0.$$
Besides, just as in the proof of Theorem 3.2, it is obtained that
${\bfv}_n$ converges to ${\bfv}$ in the sense of (5.42) and (5.44) and the
limit ${\bfv}$ satisfies the variational inequality (5.9) and (ii) for $i=1$
in Definition 5.1. Thus $\{w,{\bfv}\}$ is a weak solution of 
$SNS^1(\beta,\gamma;h,{\bfg},w_0,{\bfv}_0)$. \hfill $\diamondsuit$
\vspace{1cm}

 \begin{center}
{\large \bf References}
 \end{center}

\begin{enumerate}
\item H. Abels, Longtime behavior of solutions of a Navier-Stokes/Cahn-Hilliard
 system, pp. 9--19 in {\it Nonlocal and abstract parabolic equations and their
applications}, Banach Center Publ. {\bf Vol. 86}, Polish Acad. Sci., Inst. 
Math., Warsaw, 2009.
\item H.W. Alt, {\it Linear Functional Analysis}, Springer, Berlin Heidelberg, 
2012.
\item J. P Aubin, Un th\'eor\`eme de
compacit\'e, C. R. Acad. Sci. Paris,
{\bf 256}(1963), 5042-5044.
\item M. Biroli, Sur l'in\'equation d'{\'e}volution de Navier--Stokes.
C. R. Acad. Sci. Paris Ser. A-B, {\bf 275}(1972), A365--A367.
\item M. Biroli, Sur in\'equation d'{\'e}volution de Navier--Stokes.
Nota I, Atti Accad. Naz. Lincei Rend. Cl. Sci. Fis. Mat. Natur. (8), {\bf 52}
(1972), 457--460.
\item M. Biroli, Sur in\'equation d'{\'e}volution de Navier--Stokes.
Nota II. Atti Accad. Naz. Lincei Rend. Cl. Sci. Fis. Mat. Natur. (8), {\bf 52}
(1972), 591--598.
\item M. Biroli, Sur in\'equation d'{\'e}volution de Navier--Stokes.
Nota III, Atti Accad. Naz. Lincei Rend. Cl. Sci. Fis. Mat. Natur. (8), {\bf 52}
(1972), 811--820
\item H. Br\'ezis, Perturbations nonlin\'eaires d'op\'erateurs
maximaux monotones, C. R. Acad. Sci. Paris, {\bf 269}(1969), 566-569.
\item H. Br\'ezis, Non linear perturbations of monotone
operators, Technical Report 25, Univ. Kansas, 1972.
\item H. Br\'ezis, {\it Op\'eratuers Maximaux Monotones et 
Semi-groupes de 
Contractions dans les espaces de Hilbert}, Math.~Studies 5, North-Holland,
Amsterdam, 1973.
\item A. Damlamian, Some results on the multi-phase Stefan problem, Comm.
Partial Differential Equations, {\bf 2} (1977), 1017--1044.
\item A. Damlamian and N. Kenmochi, Evolution equations generated by
  subdifferentials in the dual space of $H^1(\Omega)$,
  Discrete Contin. Dyn. Syst., {\bf 5}(1999), 269-278.
\item J. A. Dubinskii, Weak convergence in nonlinear elliptic
and parabolic equations, Matematicheskii Sbornik {\bf 109}(1965),
609-642 and  Amer. Math. Soc. Translations Series 2, {\bf 67}(1968),
226-258.
\item  H.J. Eberl, D.F. Parker and M.C.M. van Loosdrecht, 
A new deterministic
spatio-temporal continuum model for biofilm development, 
Computational and Mathematical Methods in Medicine, {\bf 3} (2001), 161--175.
\item L. C. Evans and R. F. Gariepy, {\it Measure Theory and Fine Properties of Functions }, CRC Press, Boca Raton--London--New York--Washington, D.C., 1992.
\item T. Fukao and N. Kenmochi, Variational inequality for the Navier-Stokes
equations with time-dependent constraint, Gakuto Intern. Ser. Math. Sci. Appl., {\bf Vol. 34} (2011), 87--102.
\item T. Fukao and N. Kenmochi, 
Parabolic variational inequalities with weakly time-dependent constraints, 
 Adv. Math. Sci. Appl., {\bf 23}(2013), 365--395. 
\item T. Fukao and N. Kenmochi, Quasi-variational inequalities 
approach to heat convection problems with temperature dependent velocity 
constraint,  Discrete Contin. Dyn. Syst., 
{\bf 35}(2015), 2523--2538.
\item G. P. Galdi, {\it An Introduction to the Mathematical Theory of the 
Navier-Stokes Equations}, Second Edition, Springer, 2011.
\item M, Gokieli, N. Kenmochi and M. Niezg\'odka, Variational inequalities of
Navier--Stokes type with time dependent constraints, J. Math. Anal. Appl., 
{\bf 449}(2017), 1229-1247.
\item M, Gokieli, N. Kenmochi and M. Niezg\'odka, Mathematical modeling for
biofilm development, Nonlinear Anal.: Real World Applications, {\bf 42}(2018), 
422-447.
\item M, Gokieli, N. Kenmochi and M. Niezg\'odka, A new compactness 
theorem for variational inequalities of parabolic type, Houston J. Math.,
 {\bf 44}
(2018), 319-350.
\item N. Kenmochi, R\'esolution de compacit\'e dans les
espaces de Banach d\'ependant du temps,
S\'eminaires d'analyse convexe, Montpellier 1979,
Expos\'e 1, 1-26.     .
\item N. Kenmochi, Solvability of nonlinear evolution equations with 
time-dependent constraints and applications, Bull. Fac. Edu., Chiba Univ., {\bf 30} (1981), 1--87.
\item N. Kenmochi, Monotonicity and compactness methods for
nonlinear variational inequalities, pp. 203-298 in { Handbook of
Differential Equations: Stationary Partial Differential Equations}
{Vol. 4}, Elsevier, Amsterdam, 2007.

\item N. Kenmochi and M. Niezg\'odka, Weak
solvability for parabolic variational inclusions
and applications to quasi-variational problems,
Adv. Math. Sci. Appl., {\bf 25}(2016), 62-97.
\item J. L. Lions, {\it Quelques m\'ethodes de r\'esolution des probl\`emes aux
limites non lin\'eaires}, Dunod Gauthier--Villars, Paris, 1969.
\item U. Mosco, Convergence of convex sets and of solutions of variational
inequalities, Advances Math., {\bf 3}(1969), 510-585.
\item Y. Murase and A. Ito, Mathematical model for the process of brewing 
Japanese sake and its analysis, Adv. Math. Sci. Appl., 
{\bf 23}(2013), 297--317.
\item  M. Peszy\'nska, A. Trykozko, G, Iltis and S. Schlueter, 
Biofilm growth in
porous media: Experiments, computational modeling at the porescale, 
and upscaling, Advances in Water Resources, 2015, 1--14.
\item G. Prouse, On an inequality related to the motion, in
any dimension, of viscous, incompressible fluids. Note I, Atti Accad. Naz. Lincei Rend. Cl. Sci. Fis. Mat. Natur. (8), {\bf 67}(1979), 191--196.
\item G. Prouse, On an inequality related to the motion,
in any dimension, of viscous, incompressible fluids. Note II, Atti Accad. Naz. Lincei Rend. Cl. Sci. Fis. Mat. Natur. (8), {\bf 67}(1979), 282--288.
\item J. Simon, Compact sets in the space of $L^p(0,T;B)$,
Ann. Mat. pura appl., {\bf 146}(1986), 65-96.
\item R. Temam, {\it Navier-Stokes Equations, Theory and Numerical Analysis},
North-Holland, Amsterdam, 1984.

\end{enumerate}
\vspace{1cm}

\label{page:e}

\end{document}